\documentclass{article}
\usepackage{amssymb,amsmath}
\usepackage{anysize}
\usepackage{psfrag}
\newtheorem{thm}{Theorem}[section]
\newtheorem{def.}{Definition}[section]

\newtheorem{prop}{Proposition}[section]
\newtheorem{cor}{Corollary}[section]

\newtheorem{lem}{Lemma}[section]

\numberwithin{table}{section}

\begin{document}
\title{On computational aspects\\
 of two classical knot invariants}
        \author{Pedro Lopes\\
        Department of Mathematics\\
        Instituto Superior T\'ecnico\\
        Technical University of Lisbon\\
        Av. Rovisco Pais\\
        1049-001 Lisbon\\
        Portugal\\
        \texttt{pelopes@math.ist.utl.pt}\\
}
\date{December 31, 2006}
\maketitle

\begin{abstract}
We look into computational aspects of two classical knot
invariants. We look for ways of simplifying the computation of the
coloring invariant and of the Alexander module. We support our
ideas with explicit computations on pretzel knots.
\end{abstract}

\section{Introduction} \label{sect:intro}

\noindent

Knots (\cite{BZ, CF, Fox, lhKauffman, Lickorish, Murasugi,
Rolfsen}) are embeddings or placements of the standard circle,
$S\sp{1}$, into $3$-space. Roughly speaking, the idea is to break
$S\sp{1}$ at one point, make the line segment so obtained go over
and under itself a number of times and finally to connect the two
ends. We consider two such embeddings to represent the same knot
if one of the embeddings can be deformed into the other. In the
sequel, {\bf knot} will stand for a class of these embeddings that
are deformable into one another or for one individual embedding,
the context will make the choice clear.

Although Knot Theory has varied and interesting applications and
connections to other fields of study, we believe that one of its
basic goals is the classification of the knots modulo
deformations, or at least an attempt to do so. Ideally, an
invariant of knots would yield different outcomes when applied to
knots which are not deformable into each other. It is still an
open question whether such an invariant exists or not. On the
other hand, invariants should be tested for their efficiency. This
can be accomplished by using classes of knots and testing the
efficiency of a given invariant in telling apart the elements of
the given class. The next step is to look into the possible
simplifications in the calculation of the invariant. In this
article we will focus on these issues of testing and simplifying
on two classical invariants of knots: the coloring invariant and
the Alexander module. Moreover, we will be using the class of
pretzel knots as our working example.

\begin{figure}[h!]
    \psfrag{a}{\huge $a$}
    \psfrag{a'}{\huge $a\sb{1}$}
    \centerline{\scalebox{.50}{\includegraphics{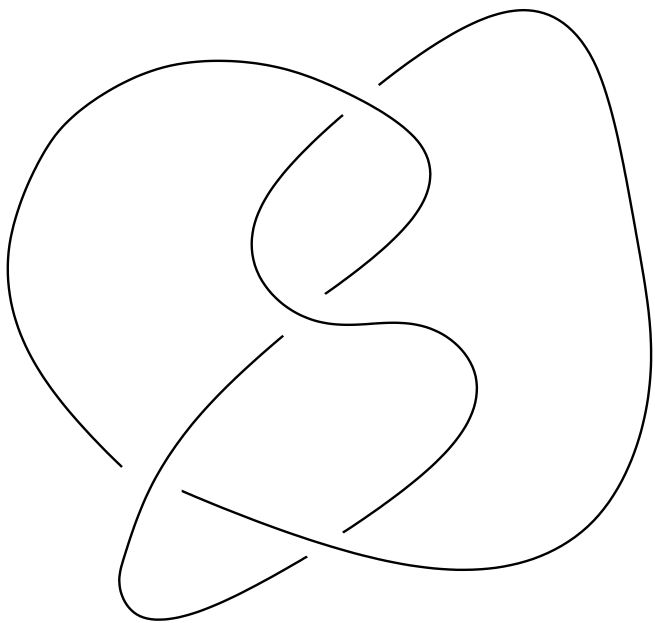}}}
    \caption{A diagram of the Figure Eight knot}\label{Fi:f8}
\end{figure}

Knots are usually represented by drawing a diagram on a plane.
This begins with a projection of the embedding on a plane,
possibly after some deformation of the embedding so that each
point of intersection in the projection has the following
property. There is a neighborhood of each such point such that the
intersection of this neighborhood with the projection of the knot
is formed by exactly two arcs meeting transversally. These points
of intersection are called crossings. At crossings, in order to
make sense of what goes over and what goes under in the embedding,
the line that goes under in the embedding is broken in the
projection. The result so obtained is called a {\bf knot diagram},
see Figure \ref{Fi:f8}.

\begin{figure}[h!]
    \psfrag{a}{\huge $a$}
    \psfrag{a'}{\huge $a'$}
    \centerline{\scalebox{.50}{\includegraphics{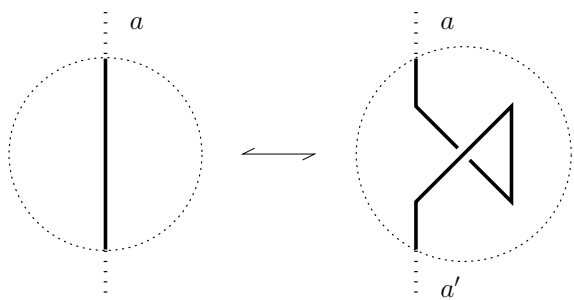}}}
    \caption{Type I Reidemeister move}\label{Fi:r1}
\end{figure}

\begin{figure}[h!]
    \psfrag{a}{\huge $a$}
    \psfrag{b}{\huge $b$}
    \psfrag{a'}{\huge $a\sb{1}'$}
    \psfrag{a''}{\huge $a\sb{2}'$}
    \centerline{\scalebox{.50}{\includegraphics{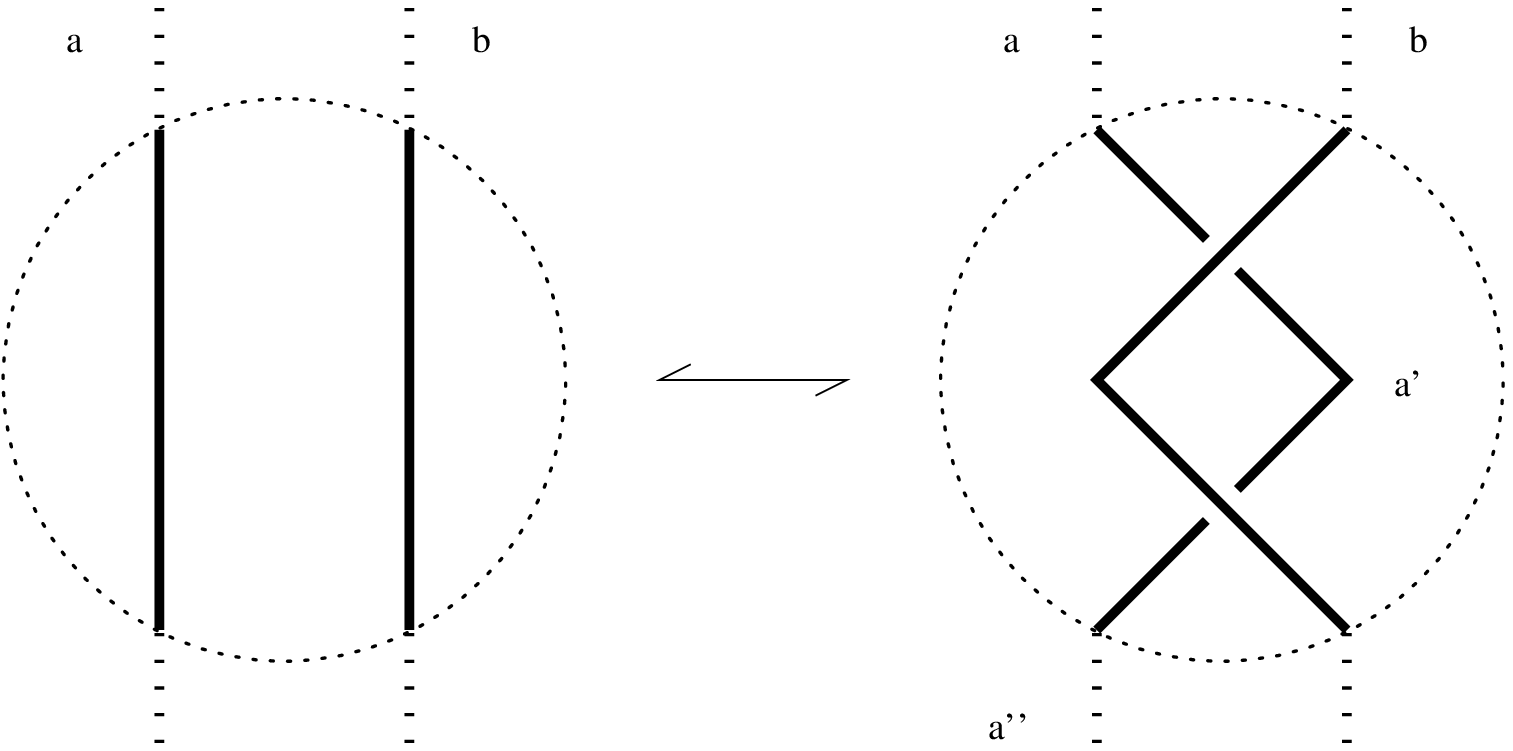}}}
    \caption{Type II Reidemeister move}\label{Fi:r2}
\end{figure}

\begin{figure}[h!]
    \psfrag{a}{\huge $a$}
    \psfrag{b}{\huge $b$}
    \psfrag{c}{\huge $c$}
    \psfrag{a1}{\huge $a\sb{1}$}
    \psfrag{a2}{\huge $a\sb{2}$}
    \psfrag{a1'}{\huge $a\sb{1}'$}
    \psfrag{a2'}{\huge $a\sb{2}'$}
    \psfrag{b1}{\huge $b\sb{1}$}
    \psfrag{b1'}{\huge $b\sb{1}'$}
    \centerline{\scalebox{.50}{\includegraphics{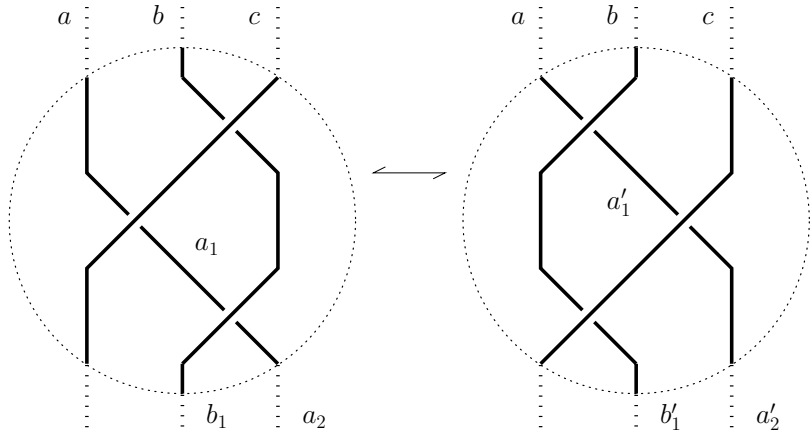}}}
    \caption{Type III Reidemeister move}\label{Fi:r3}
\end{figure}

Knot diagrams are useful in characterizing and visualizing knots.
Another important aspect of knot diagrams stems from the so-called
Reidemeister moves. These are transformations on knot diagrams
known as Reidemeister moves of type I, II, and III as shown in
Figures \ref{Fi:r1}, \ref{Fi:r2}, and \ref{Fi:r3}. Their interplay
with knot diagrams is stated in Theorem \ref{thm:reid}.

\begin{thm}[Reidemeister-Alexander]\label{thm:reid} If two knot
diagrams are related by a finite number of Reidemeister moves,
then the corresponding knots are deformable into each other.
Conversely, if two knots are deformable into each other, then any
diagram of one of these two knots is related to any diagram of the
other knot by a finite sequence of Reidemeister moves.
\end{thm} $\hfill\blacksquare$

See \cite{Reidemeister} or \cite{lhKauffman00} for a proof. The
Reidemeister moves are local i.e., the changes they make on the
diagrams occur only inside a given neighbourhood of the diagram.
Outside this neighbourhood the diagram remains the same.

A possible way of obtaining a knot invariant is to construct an
assignment that uses data from the knot diagrams and whose outcome
(essentially) does not change upon the performance of the
Reidemeister moves. In Section \ref{sect:backcols} we will see
examples of such  assignments.

\bigbreak

\subsection{Organization and acknowledgements}

\noindent

This article is organized as follows. In Section
\ref{sect:backcols} we develop the background material on
colorings, filling in what we believe to be a gap in the
literature. In Section \ref{sect:pretzelcols} we calculate the
coloring matrices for pretzel knots. In Section \ref{sect:seif} we
calculate the Alexander polynomial for two subclasses of pretzel
knots partially recovering, but also extending, previous work of
Parris (\cite{Parris}).

\bigbreak

The author acknowledges partial support from {\em Programa
Operacional ``Ci\^{e}ncia, Tecnologia, Inova\c{c}\~{a}o''} (POCTI)
of the {\em Funda\c{c}\~{a}o para a Ci\^{e}ncia e a Tecnologia}
(FCT) cofinanced by the European Community fund FEDER.

\bigbreak

\section{Background material on colorings}\label{sect:backcols}

\noindent

\bigbreak

We begin this Section by defining {\bf Integer colorings of
diagrams} and {\bf Coloring matrix of a diagram}.

\begin{figure}[h!]
    \psfrag{a}{\huge $a$}
    \psfrag{b}{\huge $b$}
    \psfrag{c}{\huge $c$}
    \psfrag{2b=a+c}{\huge $2b=a+c$}
    \centerline{\scalebox{.50}{\includegraphics{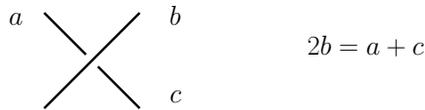}}}
    \caption{Coloring at a crossing}\label{Fi:coloratx}
\end{figure}

\begin{def.}[Integer coloring of a diagram]\label{def : integercoloring}
Given a knot $K$, consider a diagram ${\cal D}\sb{K}$ of it. An
{\bf integer coloring of } $\mathbf{{\cal D}\sb{K}}$ is an
assignment of integers to the arcs of this diagram such that, at
each crossing, twice the integer assigned to the over-arc equals
the sum of the integers assigned to the under-arcs $($see Figure
\ref{Fi:coloratx}$)$. There are always the so-called {\bf trivial
colorings}. In a {\bf trivial coloring of a diagram} the same
integer is assigned to each and every arc of the diagram.
\end{def.}

\begin{def.}[Coloring matrix of a diagram]\label{def : integercoloringmat}
We keep the notation and the terminology of Definition \ref{def :
integercoloring}. From each crossing of a knot diagram stems an
arc. In this way, there are as many arcs as crossings. Suppose
there are $\mathbf n$ arcs $($resp., crossings$)$. Enumerating the
arcs $($resp., the crossings$)$ as we go along the knot diagram
under study starting at a given reference point, we say that the
$i$-th arc stems from the $i$-th crossing. We then set up a system
of $\mathbf n$ linear equations with $\mathbf n$ variables as
follows. Each arc stands for a variable and each crossing stands
for an equation. Each of these equations states that the sum of
the variables corresponding to the under-arcs minus twice the
variable corresponding to the over-arc at the crossing at issue
equals zero. The integer solutions of this system of equations are
the {\bf Integer colorings} of Definition \ref{def :
integercoloring}. The matrix of integer coefficients of this
system of equations is the {\bf Coloring matrix of the diagram}.
\end{def.}

\bigbreak

We now define {\bf Elementary transformations on Integer matrices}
and prove that Reidemeister moves on diagrams give rise to {\bf
Coloring matrices of diagrams} which are related by a finite
number of {\bf Elementary transformations}. This being related by
a finite number of Elementary transformations is an equivalence
relation on the set of integer matrices. As a consequence, the
equivalence class of a coloring matrix of a diagram of a given
knot is a topological invariant of that knot.

\begin{def.} [Elementary Transformations on Integer Matrices
\cite{CF, Jacobson, Lickorish}] \label{def:elementransf} Given an
integer matrix $M$, an {\bf elementary transformation} on $M$ is
obtained by performing a finite sequence of the following
operations $($and/or their inverses$)$ on $M$:
\begin{enumerate}
\item Permutation of rows or columns; \item Replacement of the
matrix $M$ by \begin{equation*}\left(
\begin{matrix} $M$ & 0\\
0 & 1
\end{matrix}\right)
\end{equation*}
\item Addition of a scalar multiple of a row $($resp., column$)$
to another row $($resp., column$)$
\end{enumerate}
\end{def.}

\bigbreak

\begin{prop}\label{prop:equivcolmat}Consider a knot diagram of a knot.
 Upon performance of  Reidemeister moves on this diagram, the Coloring
matrix of the resultant diagram is obtained by performing a finite
sequence of elementary transformations on the Coloring matrix of
the original diagram.
\end{prop} Proof: We split the proof into three parts, each one
concerning one type of Reidemeister move.
\begin{itemize}
\item Type I Reidemeister move, see Figure
\ref{Fi:r1}.\\
Upon performance of this type of move (in one direction) the new
diagram has an extra arc and crossing with respect to the original
diagram, see Figure \ref{Fi:r1}. Let this extra arc be denoted
$a'$ and let the remaining arcs keep the notation of the previous
diagram. The equation at the new crossing reads
\[
a+a'-2a'=0 \quad \Leftrightarrow \quad a'=a \quad \Leftrightarrow
\quad a'+a-2a=0
\]
which states that the new variable equals one of the old
variables. Let us now see what happens with the Coloring matrices.
We will start from the Coloring matrix corresponding to the
diagram on the right of Figure \ref{Fi:r1}. This matrix has an
extra row and an extra column due to the extra arc $a'$. We write
down only the entries in the rows and columns that have directly
to do with $a'$. In this way,  there is a row formed by $1$ and
$1-2$ which corresponds to the equation read at the extra crossing
$a+a'-2a'=0$. The row with a sole $1$ corresponds to the
contribution of the other end of arc $a'$ to another equation
where it enters as an under-arc (hence the $1$). This portion of
the Coloring matrix is the one down to the left. The other
matrices were obtained by performing Elementary transformations
which we hope are clear to the reader. In the matrix to the far
right we will then perform a Transformation $2$. to get rid of the
row and column with the sole $-1$. The matrix so obtained
corresponds to the diagram on the left of Figure \ref{Fi:r1}.

\begin{equation*}
\left(
\begin{matrix}
1 & 1-2 & \dots \\
0 & 1 & \dots
\end{matrix}
\right)
\longleftrightarrow
\left(
\begin{matrix}
0 & -1 & \dots \\
1 & 1 & \dots
\end{matrix}
\right)
\longleftrightarrow
\left(
\begin{matrix}
0 & -1 & \dots \\
1 & 0 & \dots
\end{matrix}
\right)
\end{equation*}

\item Type II Reidemeister move, see Figure
\ref{Fi:r2}. \\
Upon performance of this type of Reidemeister move (in one
direction) the new diagram has two extra arcs and two extra
crossings, with respect to the original diagram, see Figure
\ref{Fi:r2}. Let these extra arcs be denoted $a\sb{1}'$ and
$a\sb{2}'$. The equations at the new crossings read
\[
\begin{cases}
a+a\sb{1}'-2b=0 \\
a\sb{1}'-a\sb{2}'-2b=0
\end{cases}
 \Leftrightarrow\quad
\begin{cases}
a\sb{1}'=2b-a \\
a\sb{2}'=a
\end{cases}
\]

We write down the progression of the portions of the matrices
relevant to this move, under the Elementary Transformations. In
the matrix to the left, the top row corresponds to the equation
read off at the top crossing of the diagram on the right of Figure
\ref{Fi:r2}, the second row corresponds to the bottom crossing on
the same Figure, and the last row corresponds to an equation to
which the bottom left arc contributes also.

\begin{equation*}
\left(
\begin{matrix}
1 & 1 & 0 & -2\\
0 & 1 & 1 & -2\\
\dots & \dots & \dots & \dots \\
0 & 0 & 1 & \dots
\end{matrix}
\right) \longleftrightarrow \left(
\begin{matrix}
1 & 1 & 0 & 0\\
0 & 1 & 1 & 0\\
\dots & \dots & \dots & \dots \\
0 & 0 & 1 & \dots
\end{matrix}
\right)
\longleftrightarrow
\left(
\begin{matrix}
1 & 1 & 0 & 0\\
0 & 1 & 1 & 0\\
\dots & \dots & \dots & \dots \\
0 & -1 & 0 & \dots
\end{matrix}
\right) \longleftrightarrow
\end{equation*}

\begin{equation*}
\longleftrightarrow \left(
\begin{matrix}
1 & 1 & 0 & 0\\
0 & 1 & 1 & 0\\
\dots & \dots & \dots & \dots \\
1 & 0 & 0 & \dots
\end{matrix}
\right)
\longleftrightarrow \left(
\begin{matrix}
0 & 1 & 0 & 0\\
0 & 0 & 1 & 0\\
\dots & \dots & \dots & \dots \\
1 & 0 & 0 & \dots
\end{matrix}
\right)
\end{equation*}

\item Type III Reidemeister move, see Figure
\ref{Fi:r3}. \\
Without further remarks, we write down the progression of the
relevant portions of the matrices under the Elementary
Transformations.

\begin{equation*}
\left(
\begin{matrix}
1 & 1 & 0 & \dots & -2 & \dots & \dots & \dots \\
0 & 1 & 1 & \dots & \dots & \dots & \dots & -2\\
\dots & \dots & \dots & \dots & \dots & \dots & \dots & \dots \\
\dots & \dots & \dots & \dots & 1 & 1 & \dots & -2
\end{matrix}
\right)
\longleftrightarrow
\left(
\begin{matrix}
1 & 1 & 0 & \dots & -2 & \dots & \dots & \dots \\
0 & 1 & 1 & \dots & -1 & -1 & \dots & 0\\
\dots & \dots & \dots & \dots & \dots & \dots & \dots & \dots \\
\dots & \dots & \dots & \dots & 1 & 1 & \dots & -2
\end{matrix}
\right) \longleftrightarrow
\end{equation*}

\begin{equation*}
\longleftrightarrow \left(
\begin{matrix}
1 & 1 & 0 & \dots & -1 & 1 & \dots & -2 \\
0 & 1 & 1 & \dots & -1 & -1 & \dots & 0\\
\dots & \dots & \dots & \dots & \dots & \dots & \dots & \dots \\
\dots & \dots & \dots & \dots & 1 & 1 & \dots & -2
\end{matrix}
\right) \longleftrightarrow \left(
\begin{matrix}
1 & 1 & 0 & \dots & 0 & 1 & \dots & -2 \\
0 & 1 & 1 & \dots & 0 & -1 & \dots & 0\\
\dots & \dots & \dots & \dots & \dots & \dots & \dots & \dots \\
\dots & \dots & \dots & \dots & 1 & 1 & \dots & -2
\end{matrix}
\right)
\end{equation*}

\begin{equation*}
\longleftrightarrow \left(
\begin{matrix}
1 & 1 & 0 & \dots & 0 & 0 & \dots & -2 \\
0 & 1 & 1 & \dots & 0 & -2 & \dots & 0\\
\dots & \dots & \dots & \dots & \dots & \dots & \dots & \dots \\
\dots & \dots & \dots & \dots & 1 & 1 & \dots & -2
\end{matrix}
\right)
\end{equation*}

\end{itemize}
 $\hfill \blacksquare $

\bigbreak

\begin{cor}\label{cor : class} Given a knot $K$ along with one of its diagrams,
consider the coloring matrix of this diagram. The equivalence
class of the coloring matrix of this diagram under elementary
transformations as described in Definition \ref{def:elementransf},
is a topological invariant of the knot $K$.
\end{cor} Proof: Omitted.
 $\hfill \blacksquare $

\bigbreak

\begin{def.} [Coloring matrix of a knot] Given a knot $K$ we call
{\bf  Coloring matrix of $K$} any matrix which is obtained from a
{\bf  Coloring matrix of a diagram of $K$} by performing a finite
number of transformations as described in Definition
\ref{def:elementransf}.
\end{def.}

\bigbreak

In view of Corollary \ref{cor : class}, given a Coloring matrix we
would like to diagonalize it, using the elementary
transformations. That this is possible is a result of Smith on
integer matrices (\cite{Smith}) which further states that the
diagonalized matrix, say $diag (d\sb{1}, d\sb{2}, \dots , d\sb{s},
0, \dots , 0)$, can be obtained so that $d\sb{i}|d\sb{i+1}$ (see
\cite{Jacobson, {brMcDonald}}). This is called the {\bf Smith
normal form} of the matrix under consideration. It is unique
modulo multiplication of the $d\sb{i}$'s by units. The $d\sb{i}$'s
are called the invariant factors.

Using elementary transformation $2$., $\pm 1$'s can be introduced
in or removed from the diagonal thus increasing or decreasing the
diagonal's length. In this way, two Smith normal forms of the same
matrix may differ on the length of the diagonal, and the
difference of the lengths is the surplus of $\pm 1$'s one has with
respect to the other. We distinguish the Smith normal forms
without $\pm 1$'s, in our work.

\begin{def.}[Normal form of the coloring
matrix]\label{def:colnormal} Consider a knot $K$. Consider the
Smith normal form of any one of its Coloring matrices. Eliminate
the $\pm 1$'s from the diagonal of this matrix along with the
corresponding rows and columns. We call the matrix so obtained
{\bf Normal form of the Coloring matrix of $K$}.

In case the diagonal of the Smith normal form consists only of
$\pm 1$'s, we call {\bf Normal form of the Coloring matrix of $K$}
the $1\times 1$ matrix formed by the entry $1$.
\end{def.}

It is also a known result (see \cite{Jacobson}) that the invariant
factors can be calculated in the following way. Calculate the
greatest common divisors (g.c.d.'s) of the $i$-rowed minors i.e.,
the g.c.d.'s of the determinants of all $i\times i$ submatrices of
the matrix under study, denoted $\Delta\sb{i}'s$. Then, modulo
units,
\[
d\sb{1}=\Delta\sb{1}, \qquad d\sb{2}=\Delta\sb{2}\cdot
\Delta\sb{1}\sp{-1}, \qquad d\sb{3}=\Delta\sb{3}\cdot
\Delta\sb{2}\sp{-1}, \qquad \dots
\]

\begin{cor}\label{cor:irowedminors} Given  a knot $K$, the
invariant factors $($resp., the greatest common divisors of the
$i$-rowed minors$)$ of any Coloring matrix of $K$ form a set of
topological invariants of $K$. \end{cor} Proof: Omitted. $\hfill
\blacksquare $

\bigbreak

\begin{def.} [Torsion invariants] The {\bf torsion invariants} of a
knot $K$ are the invariant factors of any Coloring matrix of $K$.
\end{def.}

\bigbreak

\begin{cor}\label{cor : det0} Given a knot $K$, the determinant of
any of its coloring matrices is zero.
\end{cor} Proof: Consider a coloring matrix of $K$ obtained by
associating a system of equations to a given diagram of $K$ in the
way indicated in \ref{def : integercoloringmat}. Since the
equations read at each crossing have the form
$a\sb{i}+a\sb{i+1}-2a\sb{j}=0$, then the entries of the Coloring
matrix along each row add up to zero. If we add all columns but
the last one to the last one we obtain a new matrix which is
identical to the preceding one except at the last column which is
now made of zero's. Since the passage from the former to the
latter matrices is obtained via a finite sequence of elementary
transformations, their determinants are equal. Since the latter
has a column made of zeros the determinant is zero. $\hfill
\blacksquare $

\bigbreak

\begin{def.}[The determinant of the knot] Consider an $n\times n$
coloring matrix of a given knot $K$. The {\bf determinant of $K$}
is the greatest common divisor of the $(n-1)\times (n-1)$-rowed
minors of this coloring matrix. It is clearly a topological
invariant of $K$.
\end{def.}

\bigbreak

\begin{cor} Given a knot $K$ and a coloring matrix of a diagram of $K$, its
determinant can be computed by eliminating any row and any column
of this matrix and computing the determinant of the matrix so
obtained.\end{cor} Proof: We prove that certain linear
combinations of the rows (respect., columns) yield zero.
Furthermore the coefficients of these linear combinations are $\pm
1$. We thus prove that any row (respect., column) is a linear
combination of the remaining rows (respect., columns) hence
proving the statement in this Corollary.

In the proof of Corollary \ref{cor : det0} it was shown how to
obtain a column of $0$'s by adding all columns but the one in
stake to it. We now consider how to obtain a row of zero's. Assume
the diagram $D$ that is being used is an alternating diagram i.e.,
each arc goes over exactly one crossing. The contribution of each
arc to the coloring matrix of $D$ is then as follows. It
contributes a $2$ at the crossing of which it is the over-arc; it
contributes a $-1$ at each of the two crossings of which it is an
under-arc. In this way, each column contains one $2$ and two
$-1$'s. Then, adding all but a given row to this row changes it to
a new row made of $0$'s.

\begin{figure}[h!]
    \psfrag{-(-ai)}{\huge $-(a\sb{i}+ \dots )$}
    \psfrag{(-ai)}{\huge $(a\sb{i}+ \dots )$}
    \psfrag{-(2ai)}{\huge $-(\dots - 2a\sb{i})$}
    \psfrag{(2ai)}{\huge $(\dots - 2a\sb{i})$}
    \psfrag{ai}{\huge $a\sb{i}$}
    \psfrag{x}{\huge $+$}
    \centerline{\scalebox{.50}{\includegraphics{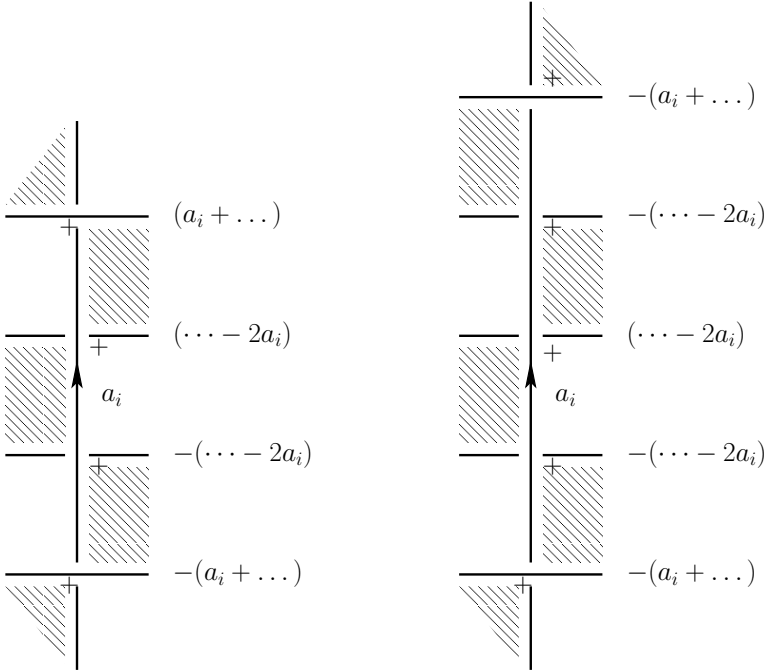}}}
    \caption{Illustrative cases of the row cancellation}\label{Fi:check1}
\end{figure}

Here is the algorithm for the general case. Consider a knot
diagram $D$ and introduce a checkerboard shading on it. This is a
shading of the faces of the diagram using two tones (grey and
white) such that at each crossing the faces that share a common
arc receive distinct tones, as (partially) illustrated in Figure
\ref{Fi:check1}.

Now, go along the diagram starting at a given point. As you go
over a crossing you print a mark to the right of the crossing and
just before it. Do this for all crossings. Note that, at each
crossing the markings either fall on grey faces or on white faces.
Also note that, at each crossing, this falling of the markings on
grey or white faces of the diagram is independent of the
orientation of the components of the diagram. As a matter of fact,
should the orientation of any of these components be reversed,
then at each crossing of this component, the marking would now
fall on the opposite side of this crossing. But shadings on
opposite sides of crossings are the same.

When you associate the $a\sb{i-1}+a\sb{i}-2a\sb{j}=0$ equation to
the $i$-th crossing, you multiply this equation by $-1$ if the
mark corresponding to this crossing is in a grey face, and you
leave the equation as it is otherwise. It is now easy to see that
adding all rows, a row of $0$'s is obtained. In Figure
\ref{Fi:check1} we show illustrative instances of this row
cancellation. The left case is intended to illustrate the case
when an $a\sb{i}$ is an over-arc for an even number of crossings
and the right case for the odd number of crossings. Finally, note
that the markings associated to the other arcs meeting at the
crossing $a\sb{i}$ stems from, comply with the markings associated
to $a\sb{i}$, so that the overall effect is to obtain a row of
zeros.
 $\hfill \blacksquare $

\bigbreak

We now know that the {\bf Normal form of the Coloring matrix} of
any knot has at least one $0$ along its diagonal. We define:

\begin{def.} [Reduced coloring matrix] A {\bf Reduced coloring
matrix} is a matrix obtained from a {\bf Coloring matrix} by
eliminating one row and one column - or any matrix equivalent to
it. The determinant of the knot is the determinant of the reduced
coloring matrix.
\end{def.}

To some extent and from the topological point of view, the {\bf
Reduced coloring matrix} is the {\bf Coloring matrix} with some
redundant information removed. Note that the procedure explained
above leading to the elimination of a row corresponds to changing
all but the last variable to a new set of variables each of which
equals the old one minus the variable associated with the given
row. Suppose the determinant of the knot is $d$. If $d$ is a unit
i.e., an element such that there is a $\bar{d}$ with
$\bar{d}d=1=d\bar{d}$ then we can use Cramer's rule (see
\cite{brMcDonald}) to obtain a unique solution for this reduced
system. According to what was previously said, this unique
solution corresponds to trivial colorings. If $d=0$ then this knot
would have more solutions then just the trivial ones.  On the
other hand since we are working with matrices over the integers,
the units are $\pm 1$. Also, there are always infinitely many
solutions even if these are only the trivial ones. Note that if
the ring we were using were finite then the number of solutions
would constitute a topological invariant of the knot. It would
also constitute a computable topological invariant of the knot.

This is the situation when we replace the ring of integers,
$\mathbb{Z}$, by the ring of residues modulo a given $r$,
$\mathbb{Z}\sb{r}$. Clearly, these rings comply with the set up
developed above.

\begin{def.} [$r$-colorings of a knot $K$] The $r$-colorings of a knot
$K$ are the solutions of the system of equations whose coefficient
matrix is a Coloring matrix of $K$ over the ring
$\mathbb{Z}\sb{r}$.
\end{def.}

\begin{cor} The number of $r$-colorings of a knot $K$ is a
topological invariant of $K$.
\end{cor}Proof: Omitted. $\hfill\blacksquare$

\bigbreak

We distinguish two types of situations, as far as $r$ is
concerned. If $r$ is prime, then $\mathbb{Z}\sb{r}$ is a field. In
this case, after diagonalizing the Coloring Matrix we count the
number of $0$'s modulo $r$ in the diagonal, say $k(\geq 1)$. This
$k$ is then the dimension of the linear subspace of the solutions;
the number of the solutions is $r\sp{k}$. $k$ is at least $1$ on
account of the trivial solutions. Moreover if $r$ divides the
determinant of the knot, then there are non-trivial solutions.

If $r$ is not prime then $\mathbb{Z}\sb{r}$ is a ring which is not
a field. In this case, if the determinant of the knot is a unit
(say $d\in \mathbb{Z}\sb{r}$ such that there is $\bar{d}\in
\mathbb{Z}\sb{r}$ and $d\bar{d}=1=\bar{d}d$) then we can use
Cramer's rule to obtain a unique solution for the reduced matrix
which means as discussed above, there are only trivial solutions.
If the determinant of the knot is not a unit (which in these rings
is equivalent to saying that $r$ is a zero divisor) then there are
more than just the trivial colorings.

The prime $r$ case is relatively known since it comes down to
linear algebra over fields, as $\mathbb{Z}\sb{r}$ is a field. The
case for non-prime $r$ is still being studied (see \cite{klopes}).

\bigbreak

There is another aspect which makes this use of rings of residues
worthwhile. As a matter of fact, given any diagram of a knot we
can set up the Coloring matrix of this diagram and proceed to
counting solutions in a given modulus $r$ by brute force. Of
course, once the diagonalization of the Coloring matrix of the
diagram is achieved, the counting of solutions in any given
modulus is very easy. But, given the known difficulties with
diagonalization of integer matrices, the method above is
particularly interesting from the practical point of view.

This point of view is explored in \cite{DL}. In this article the
coloring equations are set up for each knot from their minimal
diagrams. Then the candidates to solutions are tried out on this
system of equations. The candidates that comply with the system of
equations are the solutions and are thus counted. The knots
considered are all prime knots up to ten crossings.

\bigbreak

Our stand point in the current article is more on the theoretical
side as mentioned before. In this way, we also diagonalize the
coloring matrices of each knot of a class of knots, the pretzel
knots. This is done in Section \ref{sect:pretzelcols}.

\section{Coloring matrices of pretzel knots}\label{sect:pretzelcols}

\subsection{Preliminary formulas}

\noindent

\bigbreak

In this subsection we set up formulas for the colorings of twists.

\begin{figure}[h!]
    \psfrag{b1}{\huge $2b-a=b+(b-a)$}
    \psfrag{b2}{\huge $2(2b-a)-b=b+2(b-a)$}
    \psfrag{b3}{\huge $b+3(b-a)$}
    \psfrag{tla}{\huge $a$}
    \psfrag{tra}{\huge $b$}
    \centerline{\scalebox{.50}{\includegraphics{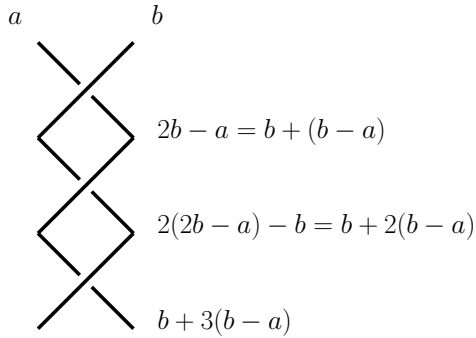}}}
    \caption{Coloring of an $n$-twist at the $n=3$ instance}\label{Fi:c3t}
\end{figure}

\begin{def.}
An {\bf $n$-twist} is obtained by producing $n$ half-twists on two
line segments.  {\bf Tassel on $n$ crossings}, and {\bf twist on
$n$ crossings} are synonyms to {\bf $n$-twist}. We refer to top
arcs and bottom arcs and order crossings as shown for the $n=3$
instance in Figure \ref{Fi:3t}. We note that this is a different
sort of ordering from that presented in Definition \ref{def :
integercoloringmat}.

\begin{figure}[h!]
    \psfrag{1st}{\huge $1\text{st crossing}$}
    \psfrag{2nd}{\huge $2\text{nd crossing}$}
    \psfrag{3rd}{\huge $3\text{rd crossing}$}
    \psfrag{tla}{\huge $\text{top (left) arc}$}
    \psfrag{tra}{\huge $\text{top (right) arc}$}
    \psfrag{bla}{\huge $\text{bottom (left) arc}$}
    \psfrag{bra}{\huge $\text{bottom (right) arc}$}
    \centerline{\scalebox{.50}{\includegraphics{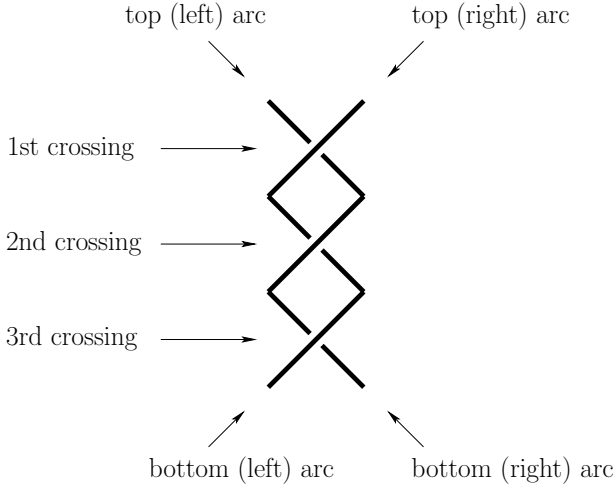}}}
    \caption{Aspect and terminology of an $n$-twist at the $n=3$ instance}\label{Fi:3t}
\end{figure}
A {\bf coloring} of an $n$-twist is the assignment of integers to
the arcs of the twist so that at each crossing the color on the
emergent under-arc is $2b-a$, where $b$ is the color on the
over-arc and $a$ is the color on the incoming under-arc (see
Figure \ref{Fi:c3t}). We will also use expressions as {\bf colored
$n$-twist} to refer to an $n$-twist endowed with a coloring.
\end{def.}

\bigbreak

Typically, such a coloring is established once the $a$ and $b$
colors on the top arcs of the twist are specified, as shown in
Proposition \ref{prop : coltwist}.

\bigbreak

\begin{figure}[h!]
    \psfrag{b1}{\huge $2b-a=b+(b-a)$}
    \psfrag{b2}{\huge $b+(n-1)(b-a)$}
    \psfrag{b3}{\huge $b+n(b-a)$}
    \psfrag{b4}{\huge $2\bigl( b+n(b-a)\bigr) - \bigl( b+(n-1)(b-a)\bigr)  =  $}
    \psfrag{b41}{\huge $ = b + (n+1)(b-a) $}
    \psfrag{tla}{\huge $a$}
    \psfrag{tra}{\huge $b$}
    \psfrag{n}{\huge $n$}
    \psfrag{...}{\huge $ \vdots $}
    \centerline{\scalebox{.50}{\includegraphics{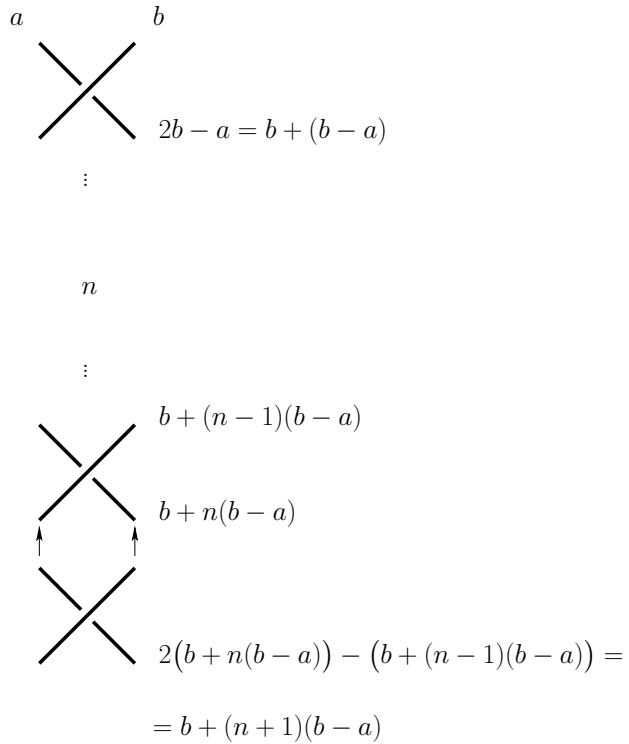}}}
    \caption{Juxtaposing the $(n+1)$-th crossing to a colored $n$-twist}\label{Fi:c3tn2}
\end{figure}

\begin{prop}\label{prop : coltwist}
Consider an $n$-twist and assign integers $a$ and $b$ to its top
arcs as shown in Figure \ref{Fi:c3t} for $n=3$. These integers
induce a coloring of the twist such that the color on the arc
emerging from the $i$-th crossing is
\[
b+i(b-a)
\]
for each $1\leq i \leq n$.
\end{prop} Proof: By induction on $n$. The $n=3$ instance is
illustrated in Figure \ref{Fi:c3t}. Assume the statement is true
for an $n>2$. Upon juxtaposition of the $(n+1)$-th crossing to a
colored $n$-twist (see Figure \ref{Fi:c3tn2}) we obtain on the arc
emerging from the $(n+1)$-th crossing:
\[
2\biggl( b+n(b-a)\biggr) - \biggl(   b+(n-1)(b-a) \biggr) = b +
(n+1)(b-a)
\]
This concludes the proof.  $\hfill \blacksquare $

\bigbreak

\bigbreak

We now introduce simpler symbols for both an $n$-twist and a
colored $n$-twist.

\bigbreak

\begin{figure}[h!]
    \psfrag{n}{\huge $n$}
    \psfrag{a}{\huge $a$}
    \psfrag{b}{\huge $b$}
    \psfrag{an}{\huge $b+(n-1)(b-a)$}
    \psfrag{bn}{\huge $b+n(b-a)$}
    \centerline{\scalebox{.50}{\includegraphics{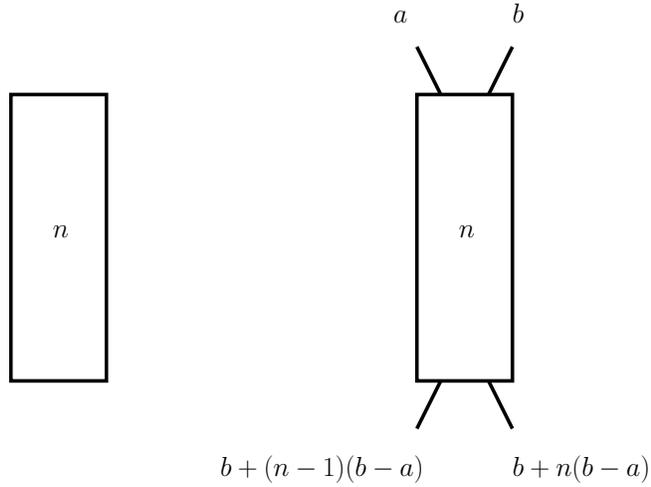}}}
    \caption{Symbolic representation of an $n$-twist and a colored $n$-twist}\label{Fi:n}
\end{figure}

\begin{def.}\label{def. : symbtwist}
Given a positive integer $n$, the symbol to the left in Figure
\ref{Fi:n} represents an $n$-twist. The symbol to the right in the
same figure represents a colored $n$-twist, where the top arcs are
assigned colors $a$ and $b$ (from left to right). Consequently,
according to Proposition \ref{prop : coltwist} the bottom left arc
bears color $b+(n-1)(b-a)$ and the bottom right arc bears color
$b+n(b-a)$.
\end{def.}

\bigbreak

We remark that, in Definition \ref{def. : symbtwist} above, $n$
does not have to be positive. In fact, given a negative $n$ the
meaning of the symbols in this Definition is analogous taken in
consideration that now the relative position of the arcs at each
crossing is changed. The under-arcs merge into one over-arc; the
over-arc splits into two under-arcs.

\subsection{The coloring matrices of pretzel knots}

\noindent

\bigbreak

We start by defining a pretzel knot.

\bigbreak

\begin{def.}
Given a positive integer $N$ and $n\sb{1}, n\sb{2}, \dots ,
n\sb{N}$ integers, a $P(n\sb{1}, n\sb{2}, \dots , n\sb{N})$
pretzel knot is represented by the knot diagram in Figure
\ref{Fi:pn} $($disregard the $a\sb{i}'s$ and $b\sb{i}'s$ for
now$)$.
\begin{figure}[h!]
    \psfrag{n1}{\huge $n\sb{1}$}
    \psfrag{n2}{\huge $n\sb{2}$}
    \psfrag{n3}{\huge $n\sb{3}$}
    \psfrag{nN}{\huge $n\sb{N}$}
    \psfrag{aN}{\huge $a\sb{1}$}
    \psfrag{a1}{\huge $a\sb{2}$}
    \psfrag{a2}{\huge $a\sb{3}$}
    \psfrag{a3}{\huge $a\sb{4}$}
    \psfrag{aN1}{\huge $a\sb{N}$}
    \psfrag{bN}{\huge $b\sb{1}$}
    \psfrag{b1}{\huge $b\sb{2}$}
    \psfrag{b2}{\huge $b\sb{3}$}
    \psfrag{b3}{\huge $b\sb{4}$}
    \psfrag{bN1}{\huge $b\sb{N}$}
    \psfrag{...}{\huge $\dots$}
    \centerline{\scalebox{.50}{\includegraphics{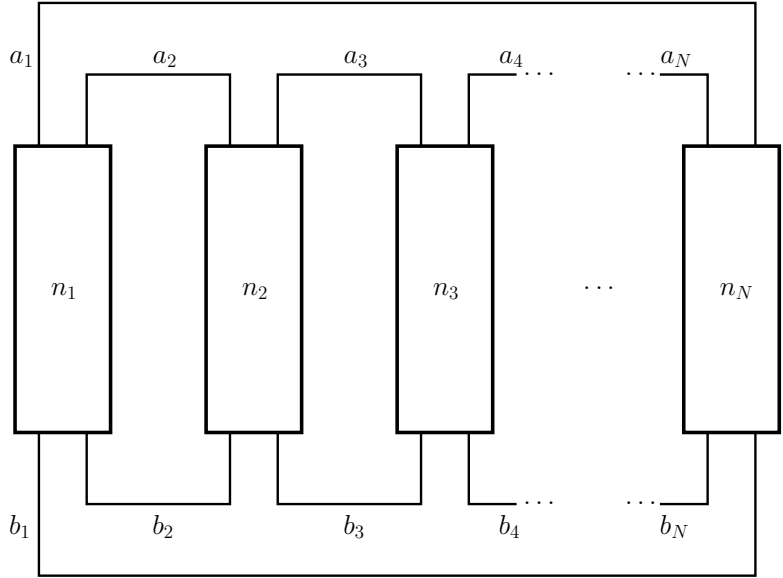}}}
    \caption{Knot diagram of a $P(n\sb{1}, n\sb{2}, n\sb{3}, \dots , n\sb{N})$
 pretzel knot}\label{Fi:pn}
\end{figure}
\end{def.}

\bigbreak

\begin{figure}[h!]
    \psfrag{n1}{\huge $n\sb{1}$}
    \psfrag{n2}{\huge $n\sb{2}$}
    \psfrag{n3}{\huge $n\sb{3}$}
    \psfrag{a3}{\huge $a\sb{1}$}
    \psfrag{a1}{\huge $a\sb{2}$}
    \psfrag{a2}{\huge $a\sb{3}$}
    \psfrag{b3}{\huge $b\sb{1}$}
    \psfrag{b1}{\huge $b\sb{2}$}
    \psfrag{b2}{\huge $b\sb{3}$}
    \psfrag{...}{\huge $\dots$}
    \centerline{\scalebox{.50}{\includegraphics{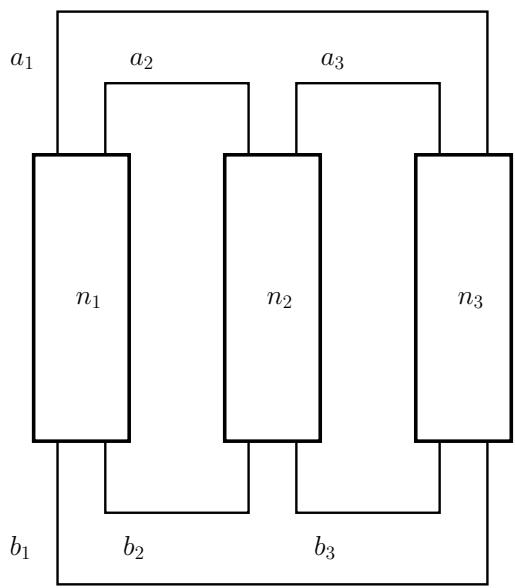}}}
    \caption{Knot diagram of a $P(n\sb{1}, n\sb{2}, n\sb{3})$ pretzel knot}\label{Fi:3n}
\end{figure}

We will now consider the $N=3$ situation with three positive
integers, $n\sb{1}, n\sb{2}, n\sb{3}$, and write down its coloring
system of equations. Consider Figure \ref{Fi:3n}.

\bigbreak

\bigbreak

From Proposition \ref{prop : coltwist} we know that once the
$a\sb{i}$ colors have been assigned to the tops of the twists that
compose the pretzel knot, the colors on the arcs along each of
these twists are uniquely expressed in terms of the $a\sb{i}$'s.
Thus, we just have to equate the two different ways of writing
each $b\sb{i}$. In this way the coloring system of equations for
this pretzel knot is
\[
\begin{cases}
b\sb{2}: & \qquad
a\sb{2}+n\sb{1}(a\sb{2}-a\sb{1})=a\sb{2}+n\sb{2}(a\sb{3}-a\sb{2})\\
b\sb{3}: & \qquad
a\sb{3}+n\sb{2}(a\sb{3}-a\sb{2})=a\sb{3}+n\sb{3}(a\sb{1}-a\sb{3})\\
b\sb{1}: & \qquad
a\sb{1}+n\sb{3}(a\sb{1}-a\sb{3})=a\sb{1}+n\sb{1}(a\sb{2}-a\sb{1})\\
\end{cases}
\]
This is equivalent to
\[
\begin{cases}
-n\sb{1}a\sb{1}+(n\sb{1}+n\sb{2}) a\sb{2}-
n\sb{2}a\sb{3}=0\\
-n\sb{3}a\sb{1}-n\sb{2}a\sb{2}+(
n\sb{2}+n\sb{3}) a\sb{3}=0\\
( n\sb{3}+n\sb{1}) a\sb{1}-n\sb{1}a\sb{2}-n\sb{3}a\sb{3}=0
\end{cases}
\]

The Coloring matrix of this diagram is then
\begin{equation*}
\left(
\begin{matrix}
-n\sb{1}  &  n\sb{1}+n\sb{2} &
 -n\sb{2}\\
-n\sb{3}  & -n\sb{2}  &
n\sb{2}+n\sb{3}\\
n\sb{3}+n\sb{1}  & -n\sb{1}  & -n\sb{3}
\end{matrix}
\right)
\end{equation*}

\bigbreak

As we saw before, along any row the coefficients add up to zero
and analogously along any column. Adding the first two columns to
the third one, and adding the first two rows to the third one, we
obtain
\begin{equation*}
\left(
\begin{matrix}
-n\sb{1}  &  n\sb{1}+n\sb{2} &
 0\\
-n\sb{3}  & -n\sb{2}  &
0\\
0  & 0  & 0
\end{matrix}
\right)
\end{equation*}
In order to keep up with what we will do later in the general
case, we add the first column to the second one to obtain
\begin{equation*}
\left(
\begin{matrix}
-n\sb{1}  &  n\sb{2} &
 0\\
-n\sb{3}  & -n\sb{3}-n\sb{2}  &
0\\
0  & 0  & 0
\end{matrix}
\right)
\end{equation*}

\bigbreak

The greatest common  divisors  of the $i$-rowed minors are now
easy to calculate.
\begin{align}\notag
\Delta\sb{1}&= ( n\sb{1}, n\sb{2}, n\sb{3}  ) \\ \notag
\Delta\sb{2}&= n\sb{1}n\sb{2}+ n\sb{1}n\sb{3}+ n\sb{2}n\sb{3}
  \\ \notag \Delta\sb{3}&=0
\end{align}
where we use the following

\begin{def.}\label{def:gcd's} Given $I$ integers, $a\sb{1}, a\sb{2}, \dots ,
a\sb{I}$, we denote their greatest common divisor $(g.c.d.)$ by
\[
(a\sb{1}, a\sb{2}, \dots , a\sb{I})  \qquad \text{ or } \qquad
(a\sb{i})\sb{i=1, \dots , I}
\]
\end{def.}

\bigbreak

With this notation, the diagonalized Coloring matrix is:

\begin{equation*}
\left(
\begin{matrix}
( n\sb{1}, n\sb{2}, n\sb{3} )   &  0 &
 0\\
0  & \frac{n\sb{1}n\sb{2}+ n\sb{1}n\sb{3}+ n\sb{2}n\sb{3}
 }{( n\sb{1}, n\sb{2},
n\sb{3} ) }  &
0\\
0  & 0  & 0
\end{matrix}
\right)
\end{equation*}

\bigbreak

Now for the general case of a pretzel knot with $N$ tassels. The
result describing the diagonalized Coloring matrix is Theorem
\ref{thm : coefmat}. It is the corollary to Propositions \ref{prop
: coefmat} and \ref{prop:gcdcolmat} and Lemma \ref{lem :
quasidiag}.

\bigbreak

\begin{prop}\label{prop : coefmat}
Given positive integers $N>2$ and $n\sb{1}, \dots , n\sb{N}$,
consider the Pretzel knot $P(n\sb{1}, \dots , n\sb{N})$ (see
Figure \ref{Fi:pn}). Its coloring matrix is equivalent to
\begin{equation*}
\left(
\begin{matrix}
-n\sb{1} &  n\sb{2}   &          0        & 0
&   \cdots     &   0       &   0   \\
    0             &  -n\sb{2}  & n\sb{3}  &
    0      &  \cdots      &   0       &   0   \\
    0      &     0     &     -n\sb{3}   &
    n\sb{4}  &      \cdots         &   0      &   0    \\
\vdots        &      \vdots    &       \vdots         & \ddots &
\vdots     &     \vdots         &   0     \\
    0    &     0    &     0     &  \cdots  &  -n\sb{N-2}   &           n\sb{N-1}  &       0\\
-n\sb{N} & -n\sb{N}  & -n\sb{N}  &  \dots  & -n\sb{N}   &   -n\sb{N}-n\sb{N-1}      &   0   \\
    0    &     0    &     0     &  \cdots      &           0  &       0     &   0
\end{matrix}
\right)
\end{equation*}
\end{prop} Proof: We remark that the $N=3$ instance is the calculation we worked
out before Proposition \ref{prop : coefmat}. For any $N\geq 3$ the
coloring matrix is

\begin{equation*}
\left(
\begin{matrix}
-n\sb{1} &  n\sb{1}+n\sb{2}   & -n\sb{2} & 0
&   \cdots     &   0     \\
    0             &  -n\sb{2}  & n\sb{2}+n\sb{3}  & -n\sb{3}
    &  \cdots      &   0   \\
\vdots        &      \vdots    &       \vdots         & \ddots &
\vdots     &     \vdots        \\
 -n\sb{N} & 0 & 0 &  \dots &        -n\sb{N-1} &
n\sb{N-1}+n\sb{N}     \\
n\sb{N}+n\sb{1}    & -n\sb{1}  &    0  & \cdots &    0 &
  -n\sb{N}
\end{matrix}
\right)
\end{equation*}

In each row (resp., column) the coefficients add up to zero. We
then add all but the last column (resp., row) to the last column
(resp., row). In this way we obtain a column  (resp., row) of
zero's for the new last column  (resp., row). We now add the first
column to the second one, the resultant second column to the third
one and so on and so forth to obtain

\begin{equation*}
\left(
\begin{matrix}
-n\sb{1} &  n\sb{2}   & 0 & 0
&   \cdots     &   0     &   0      &   0  \\
    0             &  -n\sb{2}  & n\sb{3}  & 0
    &  \cdots      &   0    &   0     &   0  \\
\vdots        &      \vdots    &       \vdots         & \ddots &
\vdots     &     \vdots     &     \vdots   &     \vdots      \\0
&    0    &    0   &   \cdots &   -n\sb{N-3} &  n\sb{N-2}   &
0 & 0 \\
0    &    0    &    0   &   \cdots &       0             &  -n\sb{N-2}  & n\sb{N-1}  &    0   \\
 -n\sb{N} &  -n\sb{N} &  -n\sb{N} &  \dots &       -n\sb{N}    &     -n\sb{N}   &         -n\sb{N}-n\sb{N-1} &
0    \\
0    & 0  &    0  & \cdots &    0    &    0   &      0 &  0
\end{matrix}
\right)
\end{equation*}
$\hfill \blacksquare $

\bigbreak

\begin{lem}\label{lem : quasidiag} For an integer $N>2$ consider
the $(N-2)\times (N-1)$ matrix
\begin{equation*}
\left(
\begin{matrix}
-n\sb{1} &  n\sb{2}   & 0 & 0
&   \cdots     &   0     &   0        \\
    0             &  -n\sb{2}  & n\sb{3}  & 0
    &  \cdots      &   0    &   0       \\
\vdots        &      \vdots    &       \vdots         & \ddots &
\vdots     &     \vdots     &     \vdots        \\0    &    0    &
0   &   \cdots &   -n\sb{N-3} &  n\sb{N-2}   &
0  \\
0    &    0    &    0   &   \cdots &       0             &
-n\sb{N-2}  & n\sb{N-1}     \end{matrix} \right)
\end{equation*}
Its $i$-rowed minors are either zero or the product of the $i$
distinct $n\sb{s}$'s from the corresponding square submatrix.
Moreover, for each $1\leq i \leq N-2$, and for each $1\leq s\sb{1}
< s\sb{2} < s\sb{3} < \cdots < s\sb{i} \leq N-1 $, there is at
least one $i$-rowed minor of the indicated matrix equal to
\[
n\sb{s\sb{1}}n\sb{s\sb{2}}n\sb{s\sb{3}} \cdots n\sb{s\sb{i}}
\]
modulo sign.
\end{lem} Proof: By induction on $i$. Let $i=2$. Any $2\times 2$
submatrix of the indicated matrix has at least one zero entry for
otherwise there would be two consecutive rows with non-zero
entries along two consecutive columns. Thus, any $2$-rowed minor
of the indicated matrix is either zero or the product of two
$n\sb{s}$'s with distinct indices, for those with the same indices
lie along the same column. Moreover, for each $1\leq s<s'\leq
N-1$, each $n\sb{s}n\sb{s'}$ equals at least one of the $2$-rowed
minors, e.g.,
\begin{equation*}
\det \left(
\begin{matrix}
-n\sb{s} &  \ast \\
    0             &  -n\sb{s'}
\end{matrix}
\right)
\end{equation*}
where
\[
\ast =
\begin{cases}
n\sb{s'}, & \text{ if } | s-s'| =1\\
0, & \text{ if } | s-s'| >1
\end{cases}
\]
Now suppose that for some $i$, for each $1 \leq j \leq i$, any
$j$-rowed minor is either zero or the product of $j$ distinct
$n\sb{s}$'s and, moreover, for any $1\leq s\sb{1}<s\sb{2}< \cdots
< s\sb{j}\leq N-1$, the product $n\sb{s\sb{1}}n\sb{s\sb{2}}\cdots
n\sb{s\sb{j}}$ is realized. Consider an $(i+1)\times (i+1)$
submatrix of the indicated matrix. If its determinant is not zero
then it has at least one column with only one non-zero entry. Then
this determinant equals this non-zero entry times an $i$-rowed
minor. Further this $i$-rowed minor is the product of $i$ distinct
$n\sb{s}$'s by the induction hypothesis and none of these
$n\sb{s}$'s equals the indicated non-zero entry - for it was the
only non-zero entry in its column.

Consider a sequence of $i+1$ integers $1\leq s\sb{1}<s\sb{2} <
\cdots <s\sb{i} < s\sb{i+1}\leq N-1$. Then
\[
n\sb{s\sb{2}} \cdots n\sb{s\sb{i}}n\sb{s\sb{i+1}}
\]
can be realized as a minor of some $i\times i$ submatrix say $M$
by the induction hypothesis. Then
\begin{equation*}
\left|
\begin{matrix}
-n\sb{s\sb{1}} &  \ast & 0 & \cdots & 0 \\
    0             &  & &  \\
\vdots &  & M & \\
0 &  & &
\end{matrix}
\right| =n\sb{s\sb{1}} n\sb{s\sb{2}} \cdots
n\sb{s\sb{i}}n\sb{s\sb{i+1}}
\end{equation*}
modulo sign, where $\ast $ is either $0$ or $n\sb{s\sb{1}}$. This
completes the proof.

$\hfill \blacksquare $

\bigbreak

\begin{prop}\label{prop:gcdcolmat}
Given an integer $N>2$ and integers $n\sb{1}, n\sb{2}, \dots ,
n\sb{N}$ consider the $N\times N$ matrix
\begin{equation*}
\left(
\begin{matrix}
-n\sb{1} &  n\sb{2}   & 0 & 0
&   \cdots     &   0     &   0      &   0  \\
    0             &  -n\sb{2}  & n\sb{3}  & 0
    &  \cdots      &   0    &   0     &   0  \\
\vdots        &      \vdots    &       \vdots         & \ddots &
\vdots     &     \vdots     &     \vdots   &     \vdots      \\0
&    0    &    0   &   \cdots &   -n\sb{N-3} &  n\sb{N-2}   &
0 & 0 \\
0    &    0    &    0   &   \cdots &       0             &  -n\sb{N-2}  & n\sb{N-1}  &    0   \\
 -n\sb{N} &  -n\sb{N} &  -n\sb{N} &  \dots &       -n\sb{N}    &     -n\sb{N}   &         -n\sb{N}-n\sb{N-1} &
0    \\
0    & 0  &    0  & \cdots &    0    &    0   &      0 &  0
\end{matrix}
\right)
\end{equation*}
Then,
\begin{align}\notag
\Delta\sb{N}&=0 \\ \notag \Delta\sb{N-1}&=\sum\sb{1\leq s\sb{1}<
s\sb{2}< \cdots < s\sb{N-1}\leq N} n\sb{s\sb{1}} n\sb{s\sb{2}}
\dots  n\sb{s\sb{N-1}} \\ \notag \Delta\sb{i}&=(n\sb{s\sb{1}}
n\sb{s\sb{2}} \dots  n\sb{s\sb{i}})\sb{1\leq s\sb{1}< s\sb{2}<
\cdots < s\sb{i}\leq N}, \qquad \text{ for } 1\leq i \leq N-2
\end{align}
\end{prop} Proof: $\Delta\sb{N}=0$ since there is one row of
zeroes (and one column of zeroes) in the indicated $N\times N$
matrix.

There is only one $(N-1)\times (N-1)$ submatrix, call it $M$, with
non-zero determinant. We calculate the determinant of $M$ by
Laplace expansion on the $(N-1)$-th row ($|M\sb{}|\sb{N-1, s}$
stands for the $N-1, s$ minor of $M$ i.e., the determinant of the
submatrix of $M$ obtained by removing the $(N-1)$-th row and the
$s$-th column from $M$):
\begin{align}\notag
\sum&\sb{s=1}\sp{N-2}(-n\sb{N})(-1)\sp{N-1+s}|M\sb{}|\sb{N-1, s} +
(-n\sb{N}-n\sb{N-1})(-1)\sp{N-1+N-1}|M\sb{}|\sb{N-1, N-1} =\\
\notag =
\sum&\sb{s=1}\sp{N-2}n\sb{N}(-1)\sp{N+s}(-n\sb{1})(-n\sb{2})
\cdots \widehat{(-n\sb{s})}n\sb{s+1}\cdots n\sb{N-1}+
(-1)\sp{1}(n\sb{N}+n\sb{N-1})(-1)\sp{N-2}n\sb{1}n\sb{2} \cdots
n\sb{N-2} =\\ \notag
=\sum&\sb{s=1}\sp{N-2}(-1)\sp{N+s}(-1)\sp{s-1}n\sb{1} \cdots
\widehat{n\sb{s}}\cdots n\sb{N-1}n\sb{N}+(-1)\sp{N-1}\biggl(
n\sb{1} \cdots  n\sb{N-2}n\sb{N} + n\sb{1} \cdots
n\sb{N-2}n\sb{N-1} \biggr) =\\ \notag =
&(-1)\sp{N-1}\sum\sb{s=1}\sp{N}n\sb{1}n\sb{2} \cdots
\widehat{n\sb{s}}n\sb{s+1}\cdots n\sb{N-1}n\sb{N}
\end{align}

We will now investigate the possible contributions of the
$i$-rowed minors to $\Delta\sb{i}$, for an otherwise arbitrary
$1\leq i \leq N-2$. Note that, in the original matrix, $n\sb{s}$
is found only along column $s$, for each $1\leq s \leq N-1$,
$n\sb{s}$, while $n\sb{N}$ is found only along row $N-1$. In this
way, $i$-rowed minors are either zero, the product of $i$ $n$'s no
two of them with the same index, or sums of such products. Thus,
these non-zero $i$-rowed minors are generated by products of $i$
$n$'s no two of them with the same index. We will now prove that
either such a product equals a specific $i$-rowed minor or is a
linear combination of some of these $i$-rowed minors, in this way
concluding the proof. Note that products not involving $n\sb{N}$
are dealt with in Lemma \ref{lem : quasidiag}. So here we will
just consider products involving $n\sb{N}$.

We first consider products of $i$ $n$'s no two of them having the
same index such that one of them has index $N$ and none of them
has index $1$ i.e., consider the sequence of indices
\[
2\leq s\sb{1} < s\sb{2} < \cdots < s\sb{i-1} < s\sb{i} = N
\]
Then
\[
n\sb{s\sb{1}}n\sb{s\sb{2}}\cdots n\sb{s\sb{i-1}}n\sb{N}
\]
is realized, modulo sign, by the determinant of the following
$i\times i$ submatrix of the indicated matrix
\begin{equation*}
\left(
\begin{matrix}
0  &  -n\sb{s\sb{1}} &  \ast &  \cdots & 0 \\
0  &   0        &  -n\sb{s\sb{2}}         &  \cdots &  0  \\
0  &   0        & 0  &  \ddots   & \vdots \\
0  &   0              &  0  &  \cdots  &  -n\sb{s\sb{i-1}} \\
-n\sb{N} & -n\sb{N} & -n\sb{N} &  \cdots & -n\sb{N}
\end{matrix}
\right)
\end{equation*}
where $\ast $ is a possible non-zero entry.

Assume now that each product of $i$ $n$'s, no two of them with the
same index, involve the factors $n\sb{1}$ and $n\sb{N}$ but do not
involve the factor $n\sb{N-1}$. Then
\[
n\sb{1}n\sb{s\sb{2}}\cdots n\sb{s\sb{i-1}}n\sb{N}
\]
(where $1<s\sb{2}< \dots < s\sb{i-1} < N-1$) is realized modulo
sign and modulo addition or subtraction of a product obtained in
Lemma \ref{lem : quasidiag}, by the determinant of the following
$i\times i$ submatrix of the indicated matrix
\begin{equation*}
\left(
\begin{matrix}
-n\sb{1} &  \ast & 0  &  \cdots & 0 \\
0        &  -n\sb{s\sb{2}}  &  \ast         &  \cdots &  0  \\
0        & 0  &  \ddots  & & \vdots \\
0  &   0              &  \cdots  &  -n\sb{s\sb{i-1}} & 0  \\
-n\sb{N} & -n\sb{N} &   \cdots & -n\sb{N}  & -n\sb{N}-n\sb{N-1}
\end{matrix}
\right)
\end{equation*}

This proof is now complete except for products of $i$ $n$'s, no
two of them with distinct indices, such that this product includes
 factors $n\sb{N-1}$ and $n\sb{N}$. If there were a formally distinct expression
for these products, this would indicate that there would be
something special about the $(N-1)$-th twist in the pretzel knot.
We see though from the symmetry of this knot that there is no
preferred twist. This concludes the proof.

$\hfill \blacksquare $

\bigbreak

\begin{thm}\label{thm : coefmat}
Given positive integers $N>2$ and $n\sb{1}, \dots , n\sb{N}$,
consider the Pretzel knot $P(n\sb{1},  \dots , n\sb{N})$ - see
Figure \ref{Fi:pn}. Its coloring matrix is equivalent to an $N
\times N$ diagonal matrix whose $i$-th entry along the diagonal is
\[
\begin{cases}
\Delta\sb{1}, & i=1\\
\Delta\sb{i}/\Delta\sb{i-1}, & 2\leq i \leq N
\end{cases}
\]
where
\begin{align}\notag
\Delta\sb{i}&=(  n\sb{s\sb{1}}n\sb{s\sb{2}} \cdots n\sb{s\sb{i}} )
\sb{1\leq s\sb{1} < s\sb{2} < \dots < s\sb{i} \leq N }    \qquad
\qquad 1\leq i \leq N-2 \\ \notag \Delta\sb{N-1}&=\sum\sb{1\leq
s\sb{1} < s\sb{2} < \dots < s\sb{N-1} \leq N}
n\sb{s\sb{1}}n\sb{s\sb{2}} \cdots n\sb{s\sb{N-1}}  \\
\notag \Delta\sb{N}&=0
\end{align}
\end{thm}Proof: This is a straightforward consequence of Propositions
\ref{prop : coefmat} and \ref{prop:gcdcolmat} and Lemma \ref{lem :
quasidiag}. $\hfill \blacksquare $

\bigbreak

Pretzel knots constitute an interesting class of knots when it
comes to coloring matrices. Rational knots are poor in this sense.
As a matter of fact, the coloring matrices of rational knots
reduce to $2\times 2$ matrices i.e., their reduced coloring
matrices have only one entry, the determinant of the knot (see
\cite{klopes2}).

\bigbreak

\section{Seifert surfaces and the Alexander module}\label{sect:seif}

\noindent

\bigbreak

In this Section we address another invariant of knots, the
so-called Alexander module of a knot. We outline the basic facts
pertaining to the Alexander module below; we elaborate on them in
the next subsection.

\bigbreak

The Alexander module is presented by a(n equivalence class of)
matrix(ces) of the form $tS-S\sp{T}$, where $S$ is a Seifert
matrix and $t$ is an indeterminate. A Seifert matrix is obtained
from a Seifert surface of the knot under study. A Seifert surface
of a knot is any connected orientable surface whose boundary is
this knot. Consider then the generators of the first homology
group of such a surface. Consider also their translates by pushing
off each generator slightly along the direction of the normal to
the surface. The Seifert matrix, $S$, is then the matrix whose
elements are the linking numbers between generators and their
translates. It is a theorem (\cite{Lickorish}) that the matrix
$tS-S\sp{T}$ is a presentation matrix for the Alexander module,
which is, further, a knot invariant. In particular, the
determinant of this matrix is a knot invariant, the so-called
Alexander polynomial.

\bigbreak

In this Section we concentrate on computational aspects of this
presentation matrix. Our working examples will be subclasses of
pretzel knots. We will partially recover work of Parris
(\cite{Parris}) on pretzel knots with an odd number of tassels,
 each tassel with an odd number of crossings.  On the
other hand we will also look into a subclass of pretzel knots not
considered by Parris.

\bigbreak

In \cite{lhKauffman0}, Kauffman outlined a way of simplifying the
calculation of the Seifert matrix. We have the same goal in mind
in the current Section although we believe we came up with a
different approach.  In order to achieve our goal, we look for a
simpler form of the Seifert surface so that the generators of the
first homology group are easily identified and the linking with
their translates easily calculated.

\bigbreak

\subsection{Definitions of Seifert surface and Seifert
matrix}\label{subsect:seifdefs}

\noindent

\bigbreak

 We start with the Definition of Seifert surface
followed by the proof of its existence, for any knot.

\bigbreak

\begin{def.}[Seifert surface of a knot] A Seifert surface of a
knot is a connected, orientable surface whose boundary is the
given knot.
\end{def.}

\bigbreak

\begin{figure}[h!]
    \centerline{\scalebox{0.65}{\includegraphics{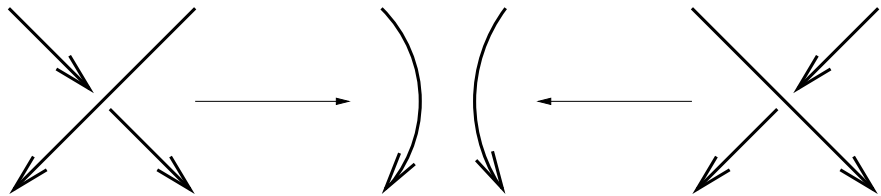}}}
    \caption{Smoothing of crossings}\label{Fi:smoothcross}
\end{figure}

\begin{prop}\label{prop:seifertalgo} There is at least one Seifert surface for each knot.
\end{prop} Proof: We consider a knot diagram of the given knot and
orient it. We smooth each crossing as shown in Figure
\ref{Fi:smoothcross}.
Upon the performance of these smoothings, the knot diagram becomes
a finite sequence of closed curves. We cap off each of these
closed curves with discs and connect them at the crossings with
half-twists. Each of these half-twists has to be consistent with
the corresponding crossing so the final result is an orientable
surface whose boundary is the given knot. If this surface is not
connected, connect any two distinct components by cutting one
small disc on each of them and gluing the boundaries of a thin
tube to the boundaries of the removed discs. This completes the
proof. $\hfill \blacksquare $

\bigbreak

\begin{figure}[h!]
    \psfrag{...}{\LARGE $\mathbf{\vdots}$}
    \psfrag{2i1+1}{\LARGE $\mathbf{2i\sb{1}+1}$}
    \psfrag{2i2+1}{\LARGE $\mathbf{2i\sb{2}+1}$}
    \psfrag{2i3+1}{\LARGE $\mathbf{2i\sb{3}+1}$}
    \psfrag{a4}{\LARGE $4b-3a$}
    \psfrag{+1}{\LARGE $+1$}
    \psfrag{-1}{\LARGE $-1$}
    \centerline{\scalebox{0.39}{\includegraphics{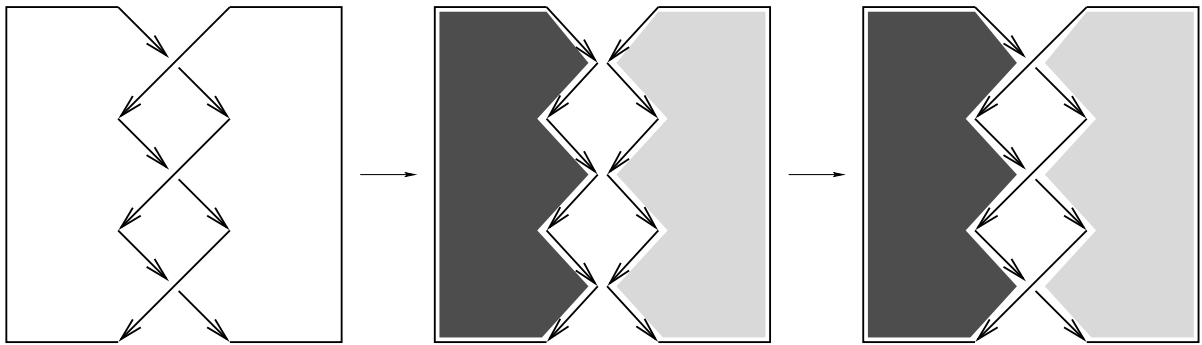}}}
    \caption{Towards a Seifert surface of the trefoil}\label{Fi:trefseifert}
\end{figure}

Figure \ref{Fi:trefseifert} illustrates this procedure for the
trefoil knot. The different shadings represent the different sides
of the surface. As a matter of fact, note that the orientation of
the knot induces an orientation of the normal to the Seifert
surface by the right-hand rule. It is this orientation of the
normal that we will consider in the Seifert surfaces, in the
sequel.

\bigbreak

\begin{figure}[h!]
    \psfrag{n}{\LARGE $\text{Normal}$}
    \psfrag{l1}{\LARGE $\mathbf{l\sb{1}}$}
    \psfrag{l2}{\LARGE $\mathbf{l\sb{2}}$}
    \centerline{\scalebox{0.6}{\includegraphics{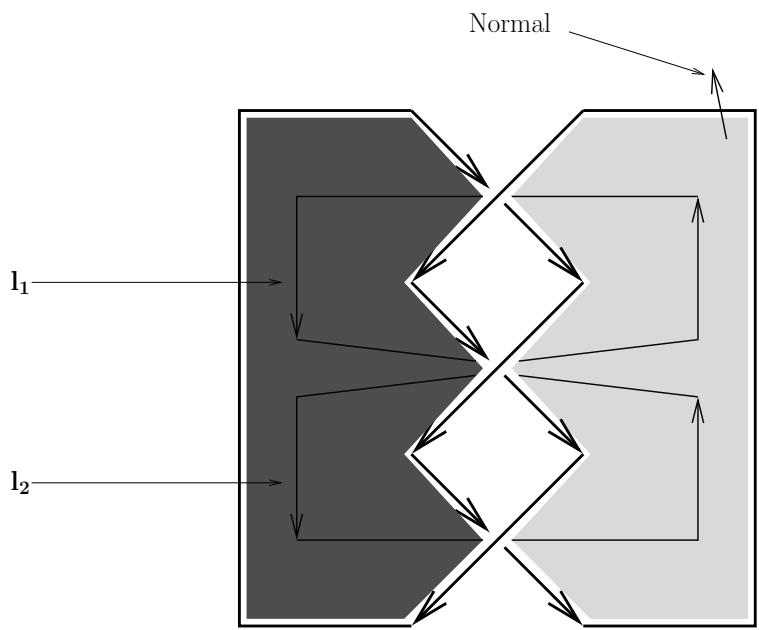}}}
    \caption{Generators of the first homology group of a Seifert surface
of the trefoil}\label{Fi:seifercircles}
\end{figure}

\begin{def.}\label{def.:seifertorient} Given a knot endowed with
an orientation and one of its Seifert surfaces, the orientation of
the normal to this Seifert surface is induced by the orientation
of the knot and the right-hand rule $($see Figure
\ref{Fi:seifercircles}$)$.
\end{def.}

\bigbreak

\bigbreak

We now proceed to define ``Signs at crossings'' and ``Linking
numbers'' in order to define a Seifert matrix of a diagram.

\bigbreak

\begin{figure}[h!]
    \psfrag{...}{\LARGE $\mathbf{\vdots}$}
    \psfrag{2i1+1}{\LARGE $\mathbf{2i\sb{1}+1}$}
    \psfrag{2i2+1}{\LARGE $\mathbf{2i\sb{2}+1}$}
    \psfrag{2i3+1}{\LARGE $\mathbf{2i\sb{3}+1}$}
    \psfrag{a4}{\LARGE $4b-3a$}
    \psfrag{+1}{\LARGE $+1$}
    \psfrag{-1}{\LARGE $-1$}
    \centerline{\scalebox{0.6}{\includegraphics{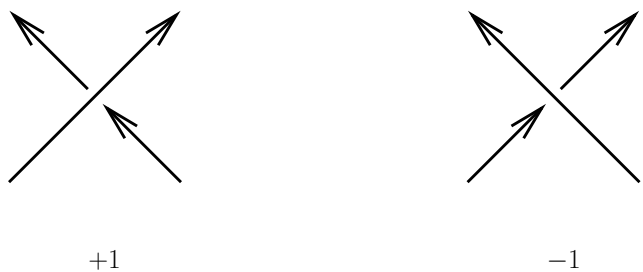}}}
    \caption{Signs at crossings}\label{Fi:link}
\end{figure}

\begin{def.}[Signs at crossings] We define signs at crossings of
an oriented knot diagram as shown in Figure \ref{Fi:link}.
\end{def.}

\bigbreak

\begin{def.} [Linking Number] Suppose $l\sb{1}$ and $l\sb{2}$ are
two knots and consider any diagram where both $l\sb{1}$ and
$l\sb{2}$ are depicted. The linking number, $lk(l\sb{1}, l\sb{2})$
of $l\sb{1}$ and $l\sb{2}$ is half the sum of the signs at the
crossings where one strand is from $l\sb{1}$ and the other one is
from $l\sb{2}$. The linking number is a topological invariant
$($\cite{BZ, lhKauffman00}$)$.
\end{def.}

\bigbreak

\begin{def.} [Seifert Matrix] Given a Seifert surface, $F$, consider a
set of generators of the first homology group of $F$, say
$\{l\sb{1}, \dots , l\sb{g} \}$. For any $i \in \{ 1, \dots , g
\}$, let $l\sb{i}\sp{+}$ denote the closed curve obtained by
pushing slightly $l\sb{i}$ in the direction of the normal to the
Seifert surface. Finally, let
\[
l\sb{ij}:=lk (l\sb{i}, l\sb{j}\sp{+})
\]
The Seifert matrix of the given Seifert surface is the square
$g\times g$ matrix whose $i, j$ entry is $l\sb{ij}$. We will
usually denote a Seifert matrix by the letter $S$.
\end{def.}

\bigbreak

\begin{figure}[h!]
    \psfrag{l1}{\LARGE $\mathbf{l\sb{1}}$}
    \psfrag{l1+}{\LARGE $\mathbf{l\sb{1}\sp{+}}$}
    \psfrag{l2}{\LARGE $\mathbf{l\sb{2}}$}
    \psfrag{l2+}{\LARGE $\mathbf{l\sb{2}\sp{+}}$}
    \psfrag{l11}{\LARGE $\mathbf{lk(l\sb{1}, l\sb{1}\sp{+})=-1}$}
    \psfrag{l12}{\LARGE $\mathbf{0=lk(l\sb{1}, l\sb{2}\sp{+})}$}
    \psfrag{l21}{\LARGE $\mathbf{lk(l\sb{2}, l\sb{1}\sp{+})=1}$}
    \psfrag{l22}{\LARGE $\mathbf{-1=lk(l\sb{2}, l\sb{2}\sp{+})}$}
    \centerline{\scalebox{0.5}{\includegraphics{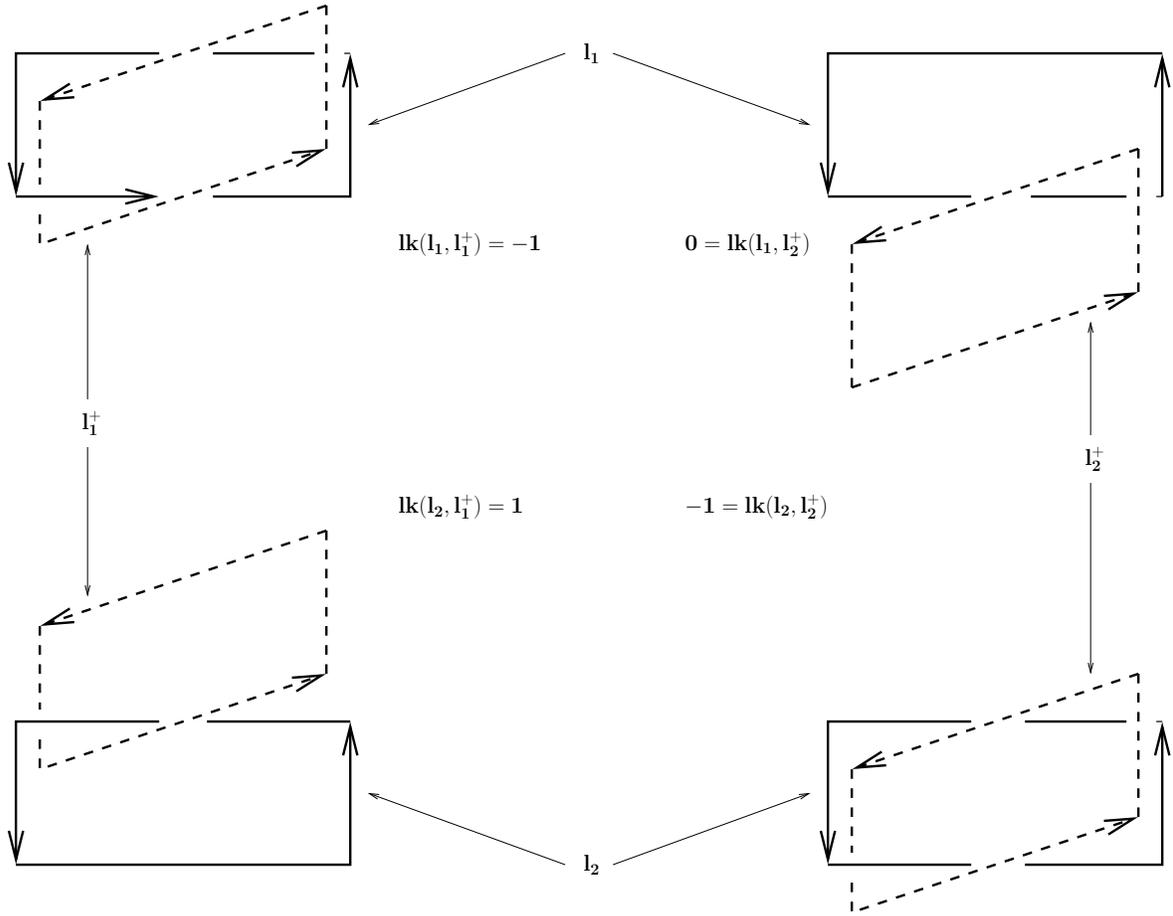}}}
    \caption{A pictorial representation of a Seifert matrix of the trefoil}
\label{Fi:seiferlinking}
\end{figure}

We calculate a Seifert matrix for the diagram of our previous
example, the trefoil (see also Figures \ref{Fi:seifercircles} and
\ref{Fi:seiferlinking}):

\begin{equation*}
S=\left(
\begin{matrix}
-1 & 0\\
1  & -1
\end{matrix}\right)
\end{equation*}

\bigbreak

\bigbreak

\bigbreak

A Seifert matrix, per se, is not a topological invariant. On the
other hand via the Seifert matrix one obtains a topological
invariant. In fact

\bigbreak

\begin{thm} Consider a knot $K$ along with a Seifert matrix of
$K$, say $S$. Let $t$ be an indeterminate. Then the matrix
$tS-S\sp{T}$ is a presentation matrix for the so-called Alexander
module of $K$ which is a topological invariant of $K$ $(S\sp{T}$
denotes the transpose of $S)$.
\end{thm} Proof: Omitted. See \cite{Lickorish}. $\hfill \blacksquare$

\bigbreak

Given a presentation matrix of a module, its elementary ideals are
defined as follows.

\bigbreak

\begin{def.}[\cite{CF, Lickorish}] Consider an $m\times n$ presentation matrix, $A$, of a
module $M$ over a commutative ring $R$. The $r$-th elementary
ideal ${\cal I}\sb{r}$ of $M$ is the ideal of $R$ generated by all
the $(m-r+1)\times (m-r+1)$ minors of $A$.
\end{def.}

\bigbreak

Elementary ideals of a module and the generators (modulo
multiplication by units) of the principal ideals contained in them
are invariants of the module at issue (\cite{CF}). In this way the
following are also topological invariants of knots.

\bigbreak

\begin{def.}[\cite{CF, Lickorish}]
The $r$-th Alexander ideal of an oriented link $L$ is the $r$-th
elementary ideal of the $\mathbb{Z}[t, t\sp{-1}]$ Alexander
module. The $r$-th Alexander polynomial of $L$ is a generator of
the smallest principal ideal of the Alexander module that contains
the $r$-th Alexander ideal. The first Alexander polynomial is
called the Alexander polynomial. We remark that the Alexander
polynomials are unique up to multiplication by units, $\pm
t\sp{n}$.
\end{def.}

\bigbreak

Resuming the study of the trefoil, we recall that a Seifert matrix
for a diagram of the trefoil read:

\begin{equation*}S= \left(
\begin{matrix}
-1 & 0\\
1 & -1
\end{matrix}
\right)
\end{equation*}

and so a presentation matrix for the Alexander module of the
trefoil is

\begin{equation*}tS-S\sp{T}= \left(
\begin{matrix}
-t+1 & -1\\
t & -t+1
\end{matrix}
\right)
\end{equation*}

Therefore, the Alexander polynomial of the trefoil is:

\[
\det \bigr( tS-S\sp{T}\bigl) = t\sp{2}-t+1
\]

\subsection{Simplifying Seifert surfaces}\label{subsect:simpleseif}

\noindent

\bigbreak

The purpose of this Subsection is to suggest a way of simplifying
the calculation of the Alexander polynomial (and if possible the
other Alexander polynomials) by simplifying the Seifert surface.
Given a knot diagram of the knot under study we calculate the
associated Seifert surface by the algorithm in Proposition
\ref{prop:seifertalgo}. We then proceed to simplifying this
Seifert surface by reducing it via deformations to the form ``disc
plus ribbons'' (see Definition \ref{def:discribbons} below). We
note that in this way each ribbon gives rise to a generator for
the first homology group of the Seifert surface (see below). In
the next Subsections we will use the pretzel knots to illustrate
the benefits of these ideas.

\bigbreak

\begin{def.}[Disc plus ribbons]\label{def:discribbons} By {\bf disc}
we mean the standard disc on the plane or one of its isotopes in
$3$-space. By {\bf ribbon} we mean a rectangle $($or one of its
isotopes in $3$-space$)$, although we will prefer rectangles such
that two parallel sides are significantly shorter than the other
sides. A {\bf disc with ribbons} or {\bf disc plus ribbons} is
then a {\bf disc} with a set of {\bf ribbons} attached to it along
the boundaries. Each {\bf ribbon} has both short sides attached to
the {\bf disc}; the long sides of the {\bf ribbons} are not
attached to anything $($see Figures \ref{Fi:pooo8a} and
\ref{Fi:seif10p5374}$)$. We remark that the ribbons may be
twisted.
\end{def.}

\bigbreak

\begin{def.}[Seifert surface in standard form]\label{def:seifstandard}
 A Seifert surface is
said in {\bf standard form} if it is in a {\bf disc plus ribbons}
form. Figures \ref{Fi:pooo8a} and \ref{Fi:seif10p5374} provide
illustrative examples.
\end{def.}

\bigbreak

In general a Seifert surface consists of several discs with
ribbons connecting them. Our strategy in transforming this Seifert
surface into a Seifert surface in standard form will be to choose
one disc and to merge the other discs into this one by shrinking
some ribbons. The preferred disc will, in general, be one with
more ribbons connecting to it.

\bigbreak

Once a Seifert surface (or for that matter, any surface) is in
standard form it is immediate to obtain a set of generators of the
first homology group of the surface. Each ribbon will give rise to
a generator in the following way. Consider this ribbon augmented
by a second ribbon embedded in the disc of the standard form of
the Seifert surface we are working with. This second ribbon
connects the short sides of the first ribbon, these connecting
short sides of both ribbons being of identical size. The ensemble
of these two ribbons constitutes a closed ribbon which retracts to
a circle embedded in three dimensions. This circle is the
generator of the first homology group we associate with the ribbon
we first considered. The retract mentioned above is realized by
shrinking the ensemble of the two ribbons along the short side.
See dotted lines in  Figure \ref{Fi:pooo9a}.

\bigbreak

\begin{def.}[Standard generators of first homology group of a
standard Seifert surface]\label{def:seifstandardgen}Given a
Seifert surface in standard form, we call the generators of the
first homology group as described in the preceding paragraph, the
{\bf standard generators} of this {\bf standard Seifert surface}.
\end{def.}

\bigbreak

\subsection{Calculations: pretzel knots on three tassels, each with an odd number of
crossings}\label{subsect:pooo}

\bigbreak

\noindent

\begin{figure}[h!]
    \psfrag{...}{\LARGE $\mathbf{\vdots}$}
    \psfrag{2i1+1}{\LARGE $\mathbf{2i\sb{1}+1}$}
    \psfrag{2i2+1}{\LARGE $\mathbf{2i\sb{2}+1}$}
    \psfrag{2i3+1}{\LARGE $\mathbf{2i\sb{3}+1}$}
    \psfrag{a4}{\LARGE $4b-3a$}
    \psfrag{1}{\LARGE $I$}
    \psfrag{3}{\LARGE $III$}
    \centerline{\scalebox{0.7}{\includegraphics{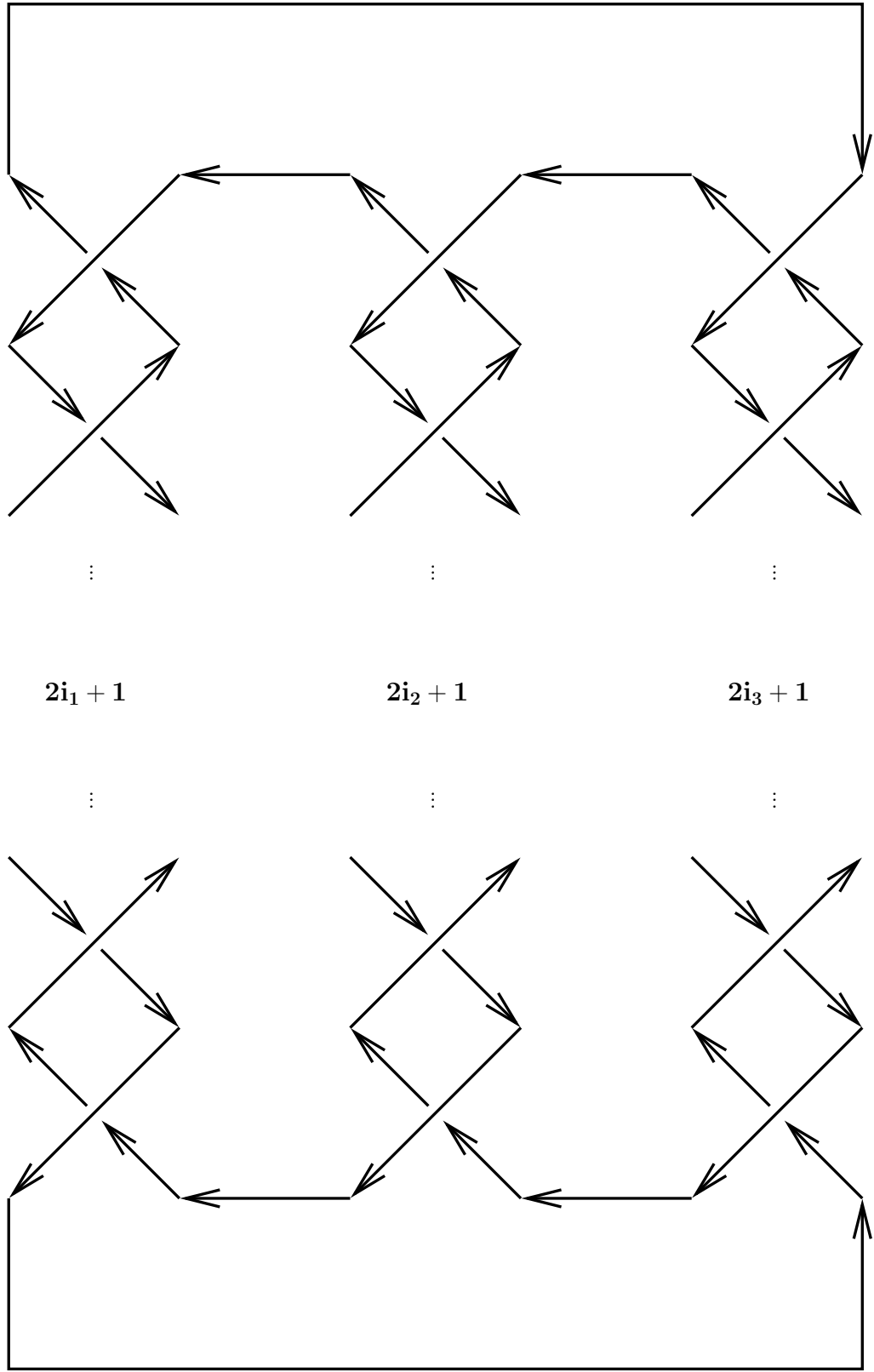}}}
    \caption{Pretzel knot on three tassels, each with an odd number of
crossings, $P(2i\sb{1}+1, 2i\sb{2}+1, 2i\sb{3}+1)$}\label{Fi:pooo}
\end{figure}

\begin{figure}[h!]
    \psfrag{...}{\LARGE $\mathbf{\vdots}$}
    \psfrag{2i1+1}{\LARGE $\mathbf{2i\sb{1}+1}$}
    \psfrag{2i2+1}{\LARGE $\mathbf{2i\sb{2}+1}$}
    \psfrag{2i3+1}{\LARGE $\mathbf{2i\sb{3}+1}$}
    \psfrag{a4}{\LARGE $4b-3a$}
    \psfrag{1}{\LARGE $I$}
    \psfrag{3}{\LARGE $III$}
    \centerline{\scalebox{0.7}{\includegraphics{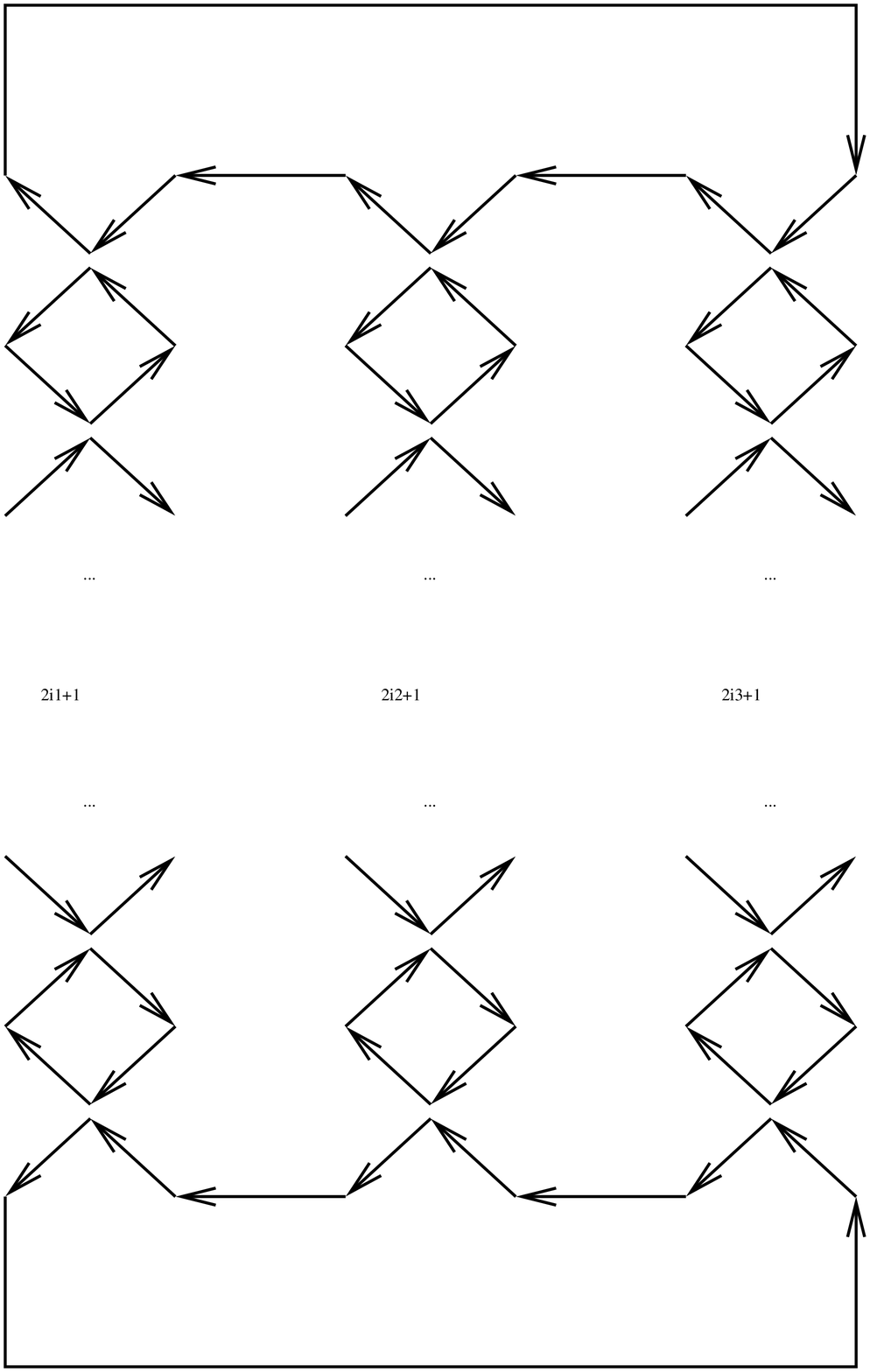}}}
    \caption{Towards constructing a Seifert surface of the
$P(2i\sb{1}+1, 2i\sb{2}+1, 2i\sb{3}+1)$ (1)}\label{Fi:pooo1}
\end{figure}

\begin{figure}[h!]
    \psfrag{...}{\LARGE $\mathbf{\vdots}$}
    \psfrag{2i1+1}{\LARGE $\mathbf{2i\sb{1}+1}$}
    \psfrag{2i2+1}{\LARGE $\mathbf{2i\sb{2}+1}$}
    \psfrag{2i3+1}{\LARGE $\mathbf{2i\sb{3}+1}$}
    \psfrag{a4}{\LARGE $4b-3a$}
    \psfrag{1}{\LARGE $I$}
    \psfrag{3}{\LARGE $III$}
    \centerline{\scalebox{0.7}{\includegraphics{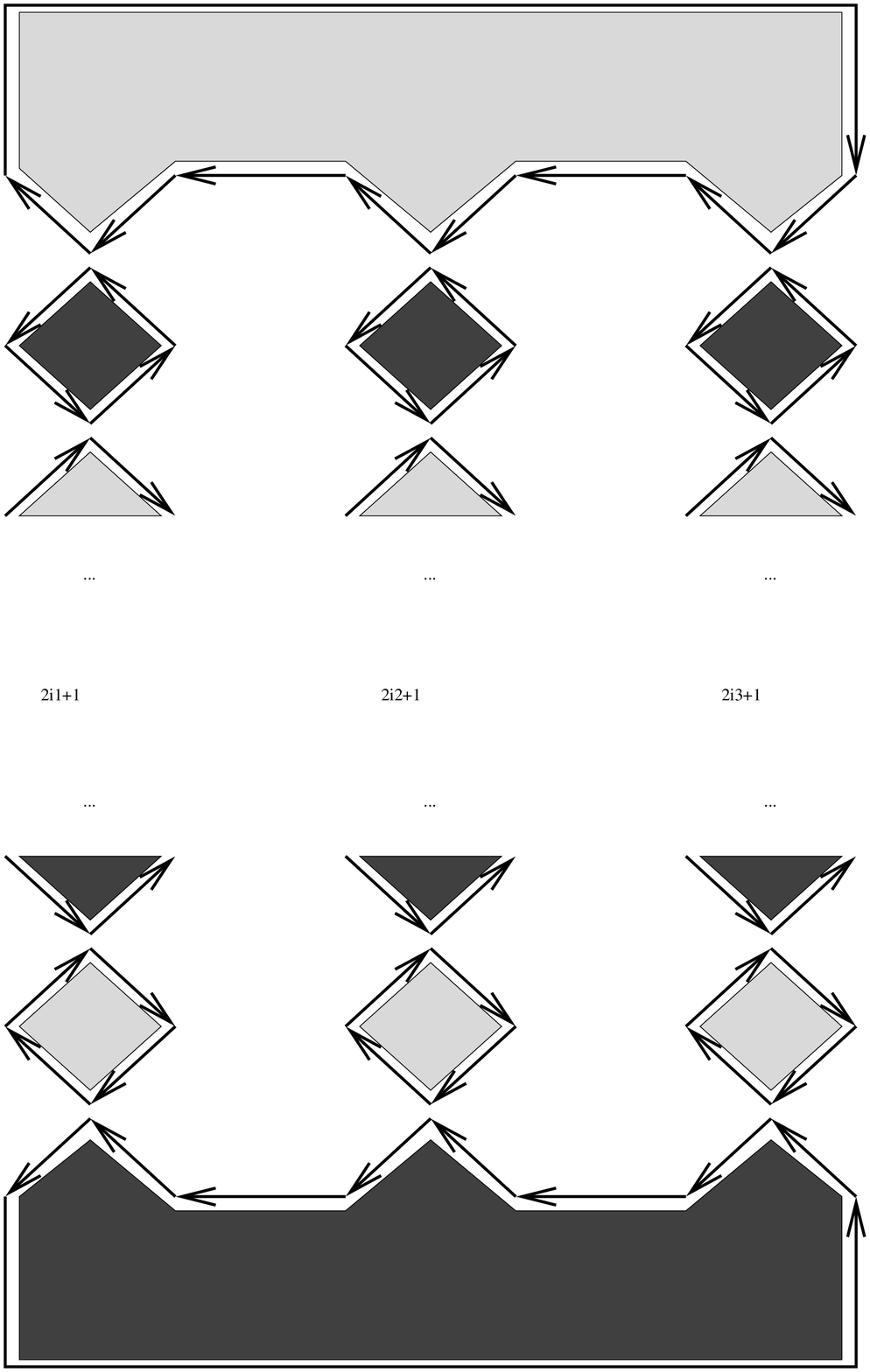}}}
    \caption{Towards constructing a Seifert surface of the
$P(2i\sb{1}+1, 2i\sb{2}+1, 2i\sb{3}+1)$ (2)}\label{Fi:pooo2}
\end{figure}

\begin{figure}[h!]
    \psfrag{...}{\LARGE $\mathbf{\vdots}$}
    \psfrag{2i1+1}{\LARGE $\mathbf{2i\sb{1}+1}$}
    \psfrag{2i2+1}{\LARGE $\mathbf{2i\sb{2}+1}$}
    \psfrag{2i3+1}{\LARGE $\mathbf{2i\sb{3}+1}$}
    \psfrag{a4}{\LARGE $4b-3a$}
    \psfrag{1}{\LARGE $I$}
    \psfrag{3}{\LARGE $III$}
    \centerline{\scalebox{0.7}{\includegraphics{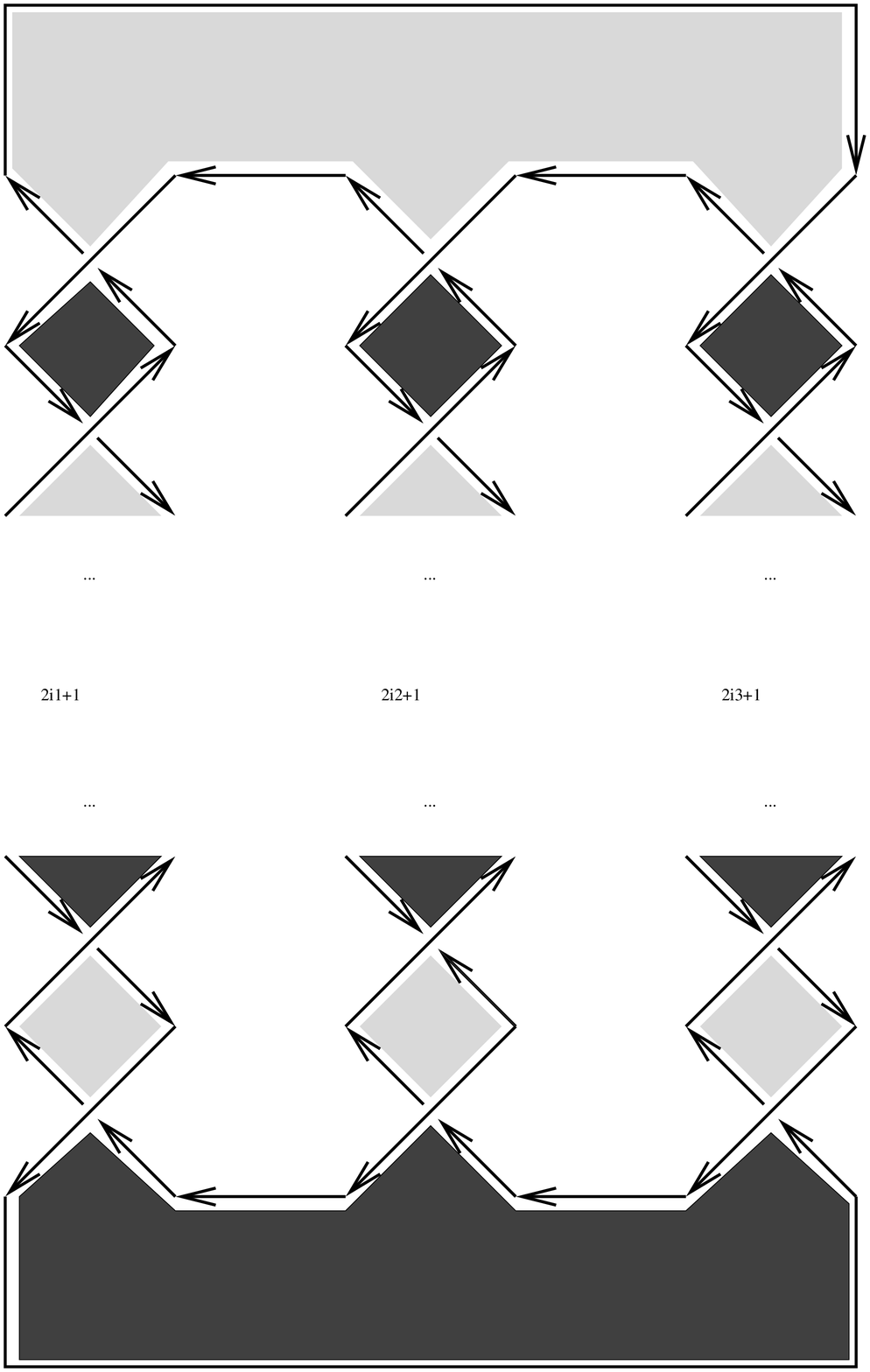}}}
    \caption{A Seifert surface of the $P(2i\sb{1}+1, 2i\sb{2}+1,
2i\sb{3}+1)$ (1)}\label{Fi:pooo3}
\end{figure}

\begin{figure}[h!]
    \psfrag{...}{\LARGE $\mathbf{\vdots}$}
    \psfrag{2i1+1}{\LARGE $\mathbf{2i\sb{1}+1}$}
    \psfrag{2i2+1}{\LARGE $\mathbf{2i\sb{2}+1}$}
    \psfrag{2i3+1}{\LARGE $\mathbf{2i\sb{3}+1}$}
    \psfrag{a4}{\LARGE $4b-3a$}
    \psfrag{1}{\LARGE $I$}
    \psfrag{3}{\LARGE $III$}
    \centerline{\scalebox{0.38}{\includegraphics{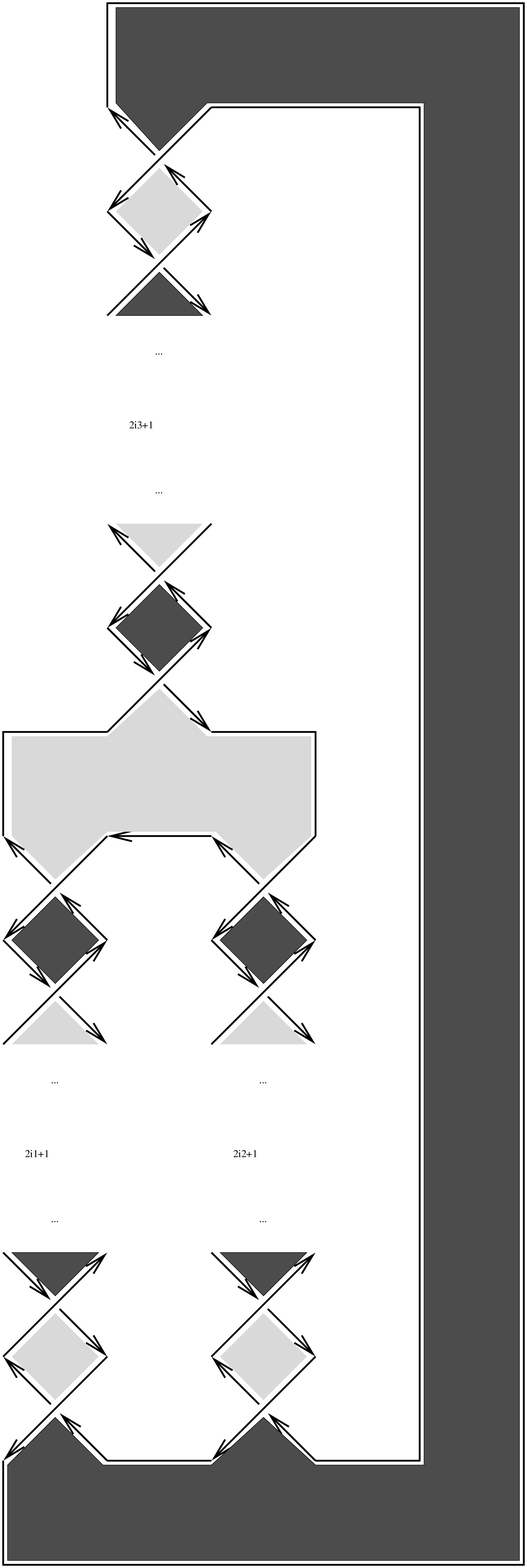}}}
    \caption{A Seifert surface of $P(2i\sb{1}+1, 2i\sb{2}+1, 2i\sb{3}+1)$ (2)}\label{Fi:pooo4a}
\end{figure}

\begin{figure}[h!]
    \psfrag{...}{\LARGE $\mathbf{\vdots}$}
    \psfrag{2i1+1}{\LARGE $\mathbf{2i\sb{1}+1}$}
    \psfrag{2i2+1}{\LARGE $\mathbf{2i\sb{2}+1}$}
    \psfrag{2i3+1}{\LARGE $\mathbf{2i\sb{3}+1}$}
    \psfrag{a4}{\LARGE $4b-3a$}
    \psfrag{1}{\LARGE $I$}
    \psfrag{3}{\LARGE $III$}
    \centerline{\scalebox{0.34}{\includegraphics{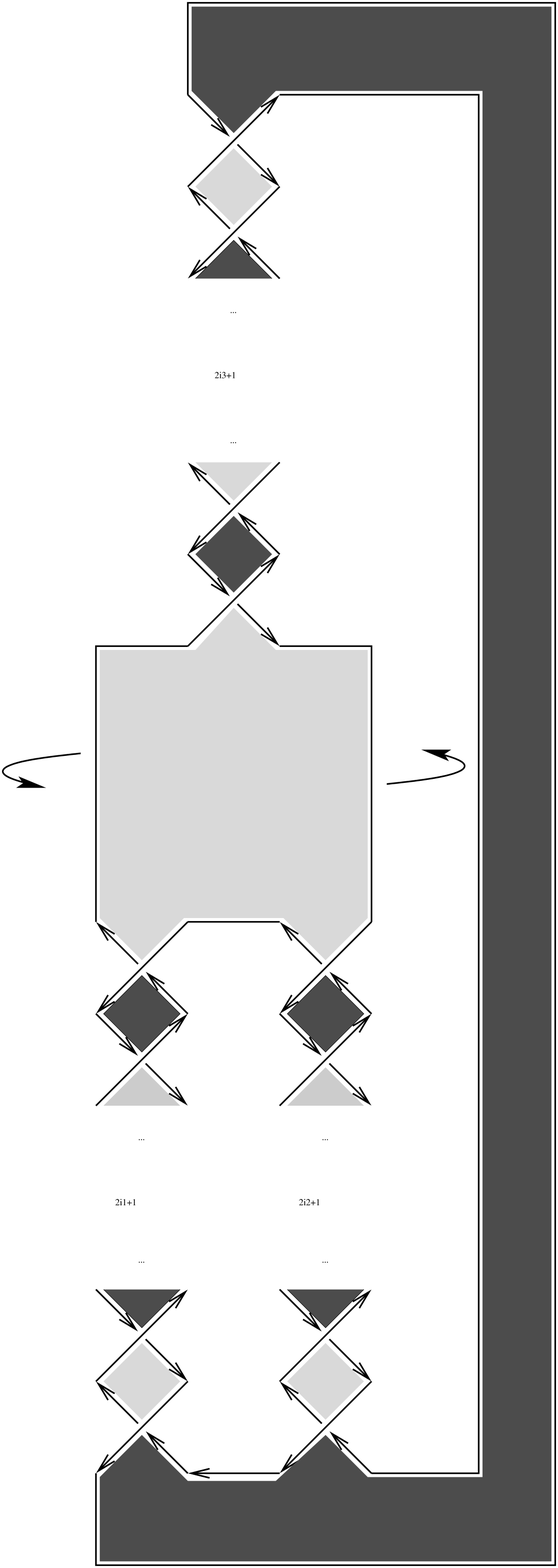}}}
    \caption{A Seifert surface of $P(2i\sb{1}+1, 2i\sb{2}+1,
2i\sb{3}+1)$ (3)}\label{Fi:pooo5a}
\end{figure}

\begin{figure}[h!]
    \psfrag{...}{\LARGE $\mathbf{\vdots}$}
    \psfrag{2i1+1}{\LARGE $\mathbf{2(i\sb{1}+i\sb{3})+2}$}
    \psfrag{2i2+1}{\LARGE $\mathbf{2(i\sb{2}+i\sb{3})+2}$}
    \psfrag{2i3+1}{\LARGE $\mathbf{2i\sb{3}+1}$}
    \psfrag{a4}{\LARGE $4b-3a$}
    \psfrag{1}{\LARGE $I$}
    \psfrag{3}{\LARGE $III$}
    \centerline{\scalebox{0.34}{\includegraphics{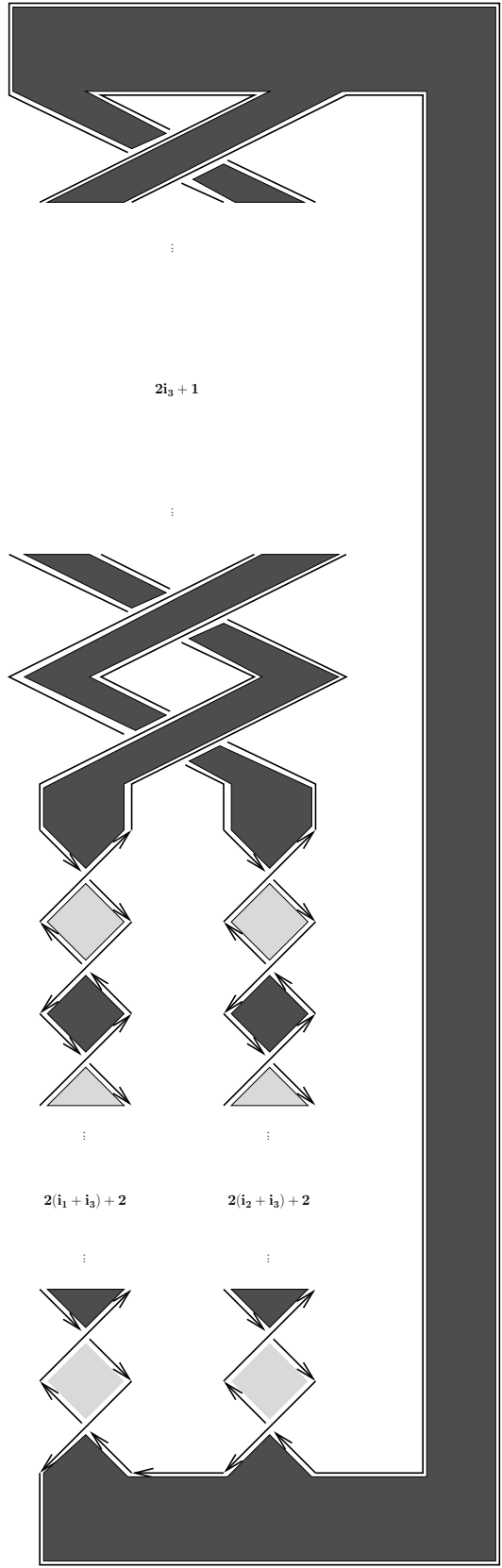}}}
    \caption{A Seifert surface of $P(2i\sb{1}+1, 2i\sb{2}+1, 2i\sb{3}+1)$
 - the standard form}\label{Fi:pooo8a}
\end{figure}

As a working example for the ideas described in the preceding
Subsection, we consider a pretzel knot on three tassels, each with
an odd number of crossings, $P(2i\sb{1}+1, 2i\sb{2}+1,
2i\sb{3}+1)$, see figure \ref{Fi:pooo}. In Figures \ref{Fi:pooo1}
and \ref{Fi:pooo2} we use the algorithm described in Proposition
\ref{prop:seifertalgo} to obtain a Seifert surface in Figure
\ref{Fi:pooo3}. We move the third tassel up (Figure
\ref{Fi:pooo4a}) and we rotate the disc in the center of Figure
\ref{Fi:pooo5a} to undo the twisting in the third tassel. As a
consequence we merge this disc into the other one giving rise to a
Seifert surface in standard form (Figure \ref{Fi:pooo8a}). This
produces also a braiding of the first two tassels and an increase
in their crossings of $2i\sb{3}+1$.

\bigbreak

\begin{figure}[h!]
    \psfrag{...}{\LARGE $\mathbf{\vdots}$}
    \psfrag{2i1+1}{\LARGE $\mathbf{2(i\sb{1}+i\sb{3})+2}$}
    \psfrag{2i2+1}{\LARGE $\mathbf{2(i\sb{2}+i\sb{3})+2}$}
    \psfrag{2i3+1}{\LARGE $\mathbf{2i\sb{3}+1}$}
    \psfrag{a4}{\LARGE $4b-3a$}
    \psfrag{l1}{\huge $l\sb{1}$}
    \psfrag{l2}{\huge $l\sb{2}$}
    \centerline{\scalebox{0.34}{\includegraphics{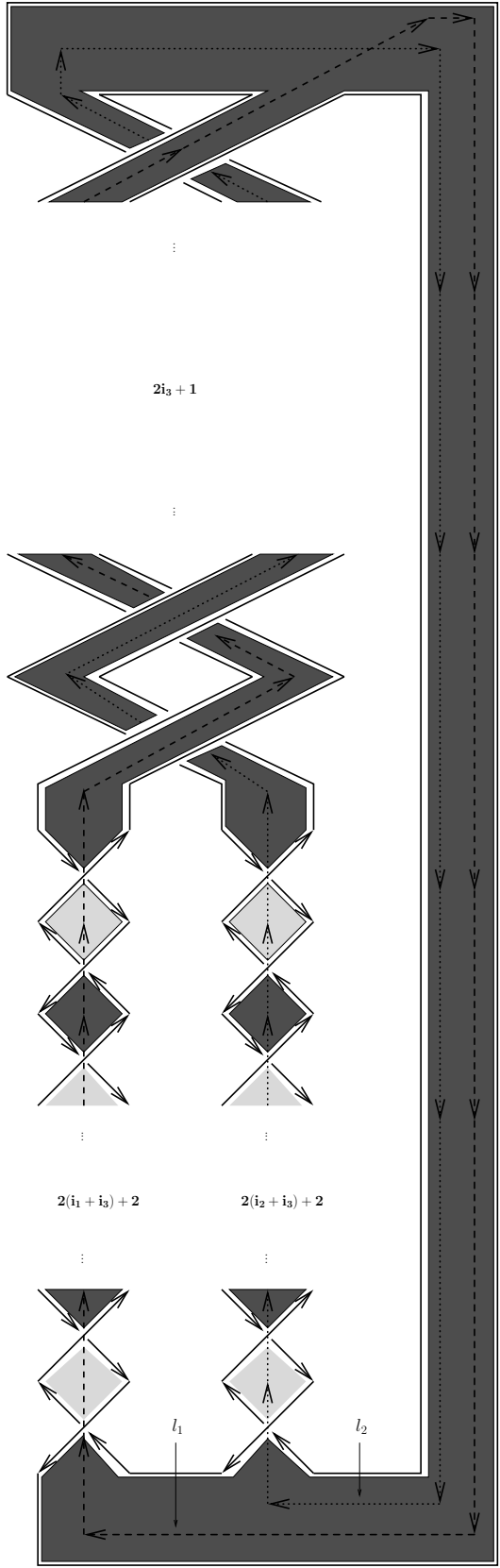}}}
    \caption{The standard form with the standard generators for the homology}\label{Fi:pooo9a}
\end{figure}

\begin{figure}[h!]
    \psfrag{...}{\LARGE $\mathbf{\vdots}$}
    \psfrag{2i1+1}{\LARGE $\mathbf{2(i\sb{1}+i\sb{3})+2}$}
    \psfrag{2i2+1}{\LARGE $\mathbf{2(i\sb{2}+i\sb{3})+2}$}
    \psfrag{2i3+1}{\LARGE $\mathbf{2i\sb{3}+1}$}
    \psfrag{a4}{\LARGE $4b-3a$}
    \psfrag{l1}{\huge $\mathbf{l\sb{1}}$}
    \psfrag{l2}{\huge $\mathbf{l\sb{2}}$}
    \psfrag{l1+}{\huge $\mathbf{l\sb{1}\sp{+}}$}
    \psfrag{l2+}{\huge $\mathbf{l\sb{2}\sp{+}}$}
    \psfrag{lk11}{\huge $\mathbf{lk(l\sb{1}, l\sb{1}\sp{+})=(i\sb{1}+i\sb{3})+1}$}
    \psfrag{lk22}{\huge $\mathbf{lk(l\sb{2}, l\sb{2}\sp{+})=(i\sb{2}+i\sb{3})+1}$}
    \centerline{\scalebox{0.6}{\includegraphics{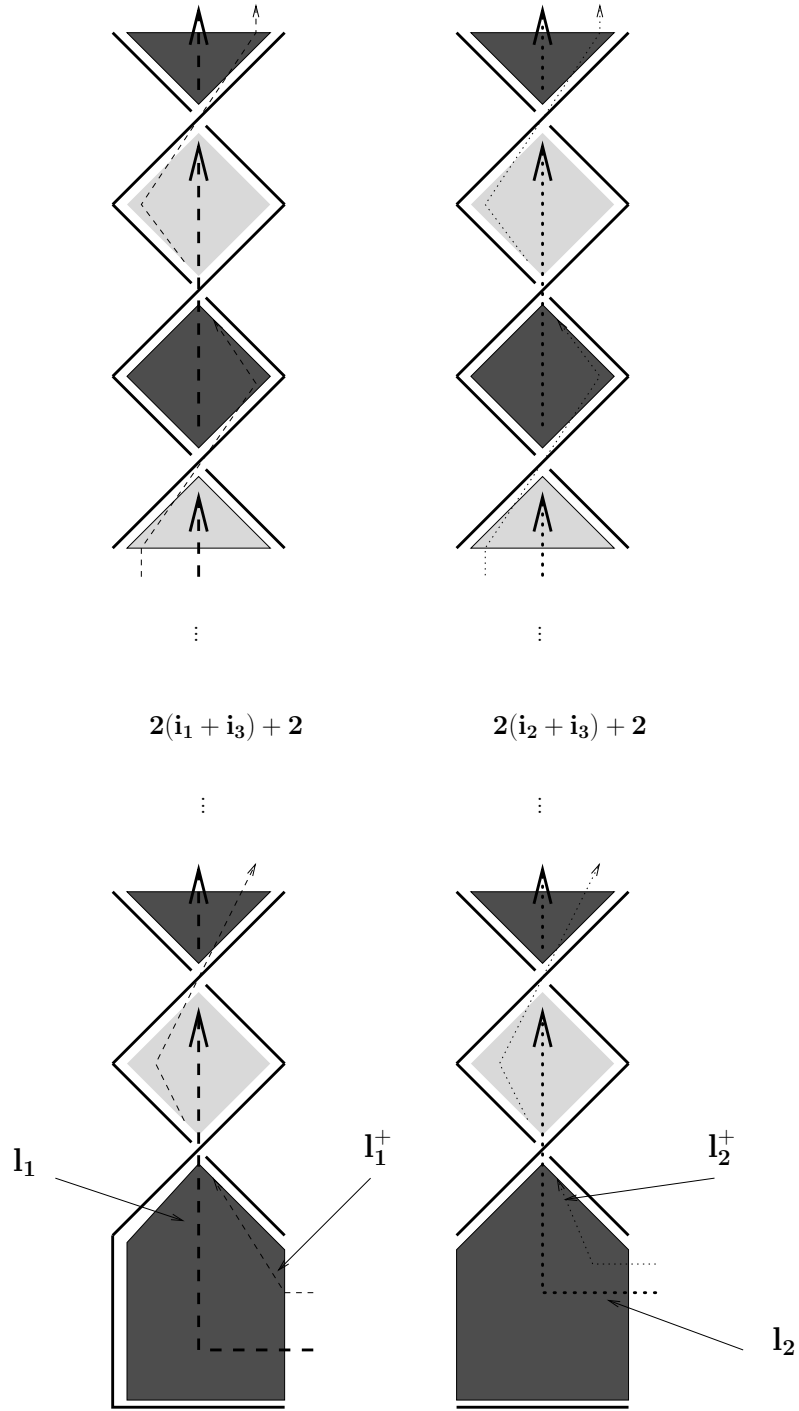}}}
    \caption{The $lk(l\sb{i}, l\sb{i}\sp{+})$ linking numbers}\label{Fi:pooolk11}
\end{figure}

\begin{figure}[h!]
    \psfrag{...}{\LARGE $\mathbf{\vdots}$}
    \psfrag{2i1+1}{\LARGE $\mathbf{2(i\sb{1}+i\sb{3})+2}$}
    \psfrag{2i2+1}{\LARGE $\mathbf{2(i\sb{2}+i\sb{3})+2}$}
    \psfrag{2i3+1}{\LARGE $\mathbf{2i\sb{3}+1}$}
    \psfrag{a4}{\LARGE $4b-3a$}
    \psfrag{l1}{\huge $\mathbf{l\sb{1}}$}
    \psfrag{l2}{\huge $\mathbf{l\sb{2}}$}
    \psfrag{l1+}{\huge $\mathbf{l\sb{1}\sp{+}}$}
    \psfrag{l2+}{\huge $\mathbf{l\sb{2}\sp{+}}$}
    \psfrag{l12}{\huge $\mathbf{lk(l\sb{1}, l\sb{2}\sp{+})=i\sb{3}+1}$}
    \psfrag{l21}{\huge $\mathbf{lk(l\sb{2}, l\sb{1}\sp{+})=i\sb{3}}$}
    \centerline{\scalebox{0.5}{\includegraphics{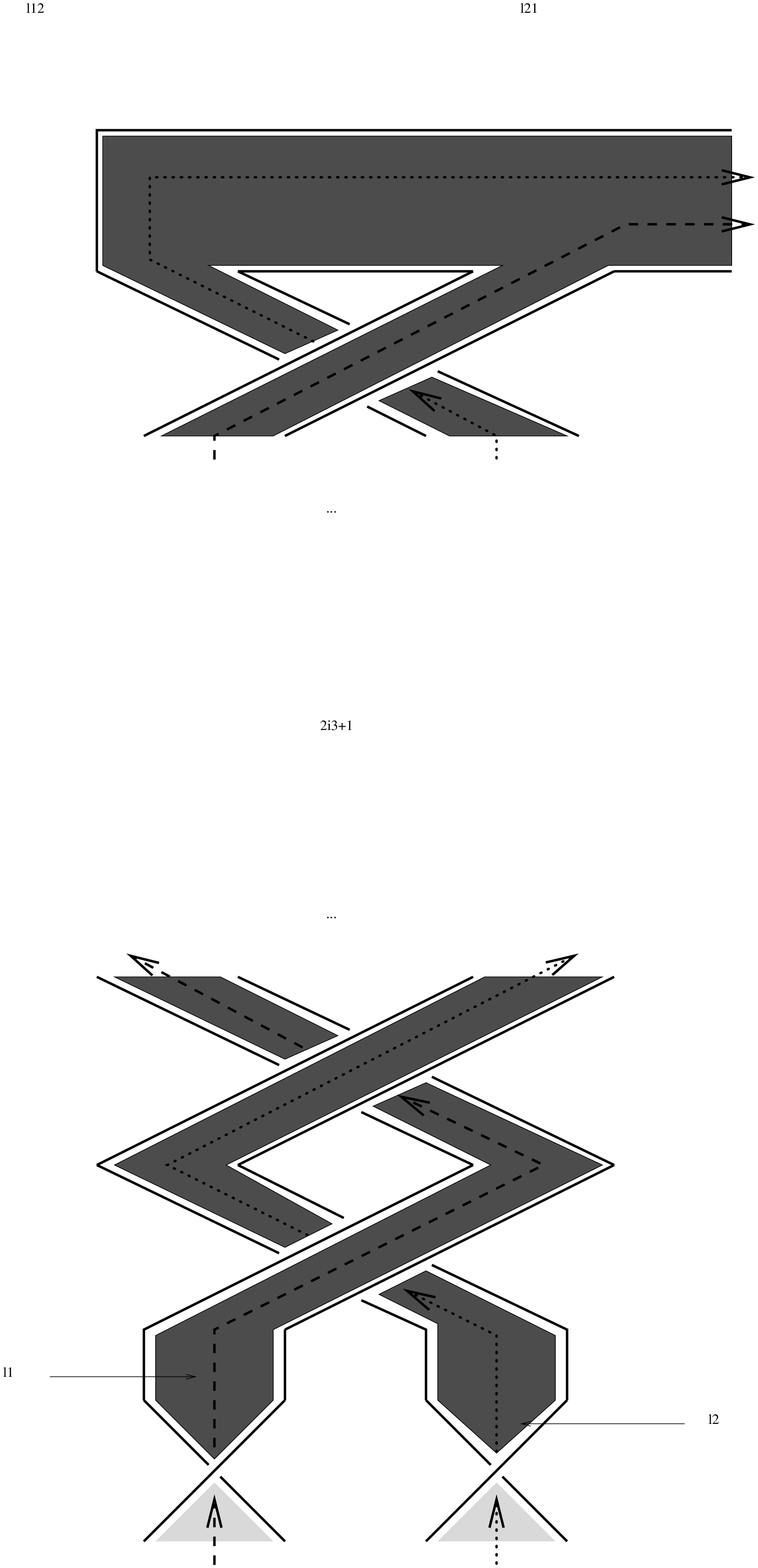}}}
    \caption{The $lk(l\sb{i}, l\sb{j}\sp{+})$ linking numbers}\label{Fi:pooolk12}
\end{figure}

With the Seifert surface in standard position, we easily identify
the standard generators of the first homology group of this
Seifert surface (see Figure \ref{Fi:pooo9a}): we associate one
generator of the first homology group of the Seifert surface to
each ribbon (and vice-versa). We then calculate the linking
numbers between generators and their translates as illustrated in
Figures \ref{Fi:pooolk11} and \ref{Fi:pooolk12}. Note that the
relevant crossings between generator and translate of generator
for this calculation are those shown in Figure \ref{Fi:pooolk11}.
As for Figure  \ref{Fi:pooolk12}, there is one relevant crossing
missing. In fact, in this case the two dotted lines will go over
one another one more time, after they leave the top right part of
this Figure. For $lk(l\sb{1}, l\sb{2}\sp{+})$, it is
$l\sb{2}\sp{+}$ that goes over $l\sb{1}$; for $lk(l\sb{2},
l\sb{1}\sp{+})$ it is $l\sb{1}\sp{+}$ that goes over $l\sb{2}$.

\bigbreak

The Seifert matrix is then:

\begin{equation*}S= \left(
\begin{matrix}
i\sb{1}+i\sb{3}+1 & i\sb{3}+1\\
i\sb{3} & i\sb{2}+i\sb{3}+1
\end{matrix}
\right)
\end{equation*}
and the presentation matrix of the Alexander module for the
pretzel knot $P(2i\sb{1}+1, 2i\sb{2}+1, 2i\sb{3}+1)$ is

\begin{equation*}tS-S\sp{T}= \left(
\begin{matrix}
 (i\sb{1}+i\sb{3}+1)(t-1) & (i\sb{3}+1 ) t -i\sb{3}\\
i\sb{3}t-(i\sb{3}+1) \bigr)  & (i\sb{2}+i\sb{3}+1) (t-1)
\end{matrix}
\right)
\end{equation*}
which is equivalent to the following matrix:

\begin{equation*}
\left(
\begin{matrix}
 i\sb{1}t-(i\sb{1}+1) & (i\sb{3}+1 ) t -i\sb{3}\\
i\sb{2}-(i\sb{2}+1)t  & (i\sb{2}+i\sb{3}+1) (t-1)
\end{matrix}
\right)
\end{equation*}

\bigbreak

It then follows that the Alexander polynomial of $P(2i\sb{1}+1,
2i\sb{2}+1, 2i\sb{3}+1)$ is

\[
\det \bigr( tS-S\sp{T}\bigl)
=i\sb{1}t-(i\sb{1}+1)(i\sb{2}+i\sb{3}+1) (t-1)- \bigl(
i\sb{2}-(i\sb{2}+1)t \bigr) \bigl(  (i\sb{3}+1 ) t -i\sb{3}
 \bigr)
\]

\bigbreak

We record here two results implicit in the preceding calculations
which will be useful in subsequent calculations.

\bigbreak

\begin{prop}[The $lk(l\sb{i}, l\sb{i}\sp{+})$ linking numbers]\label{prop : lk(i,i)}
Let $L\sb{i}$ be the ribbon which gives rise to the $l\sb{i}$
generator of the first homology group of the surface. The
contribution to $lk(l\sb{i}, l\sb{i}\sp{+})$ of the ribbon
$L\sb{i}$ as depicted in Figure \ref{Fi:pooolk11} is one half the
number of the crossings in the Figure. If each of these crossings
were reversed then it would be negative one half the number of
these crossings.
\end{prop}Proof: Omitted. $\hfill \blacksquare$

\bigbreak

\bigbreak

\begin{prop}[The $lk(l\sb{i}, l\sb{j}\sp{+})$ linking numbers]\label{prop : lk(i,j)}
Let $L\sb{i}$ and $L\sb{j}$ be the ribbons which give rise to the
$l\sb{i}$ and $l\sb{j}$ generators of the first homology group of
the surface $(i\neq j)$. The contribution to $lk(l\sb{i},
l\sb{j}\sp{+})$ and to $lk(l\sb{j}, l\sb{i}\sp{+})$ from the
braiding of $L\sb{i}$ and $L\sb{j}$ as depicted in Figure
\ref{Fi:pooolk12} is one half the number of crossings of these two
ribbons, $L\sb{i}$ and $L\sb{j}$. If each of these crossings were
reversed then it would be negative one half the number of these
crossings.
\end{prop}Proof: Omitted. $\hfill \blacksquare$

\bigbreak

 \vskip 20pt

 \clearpage

 \vskip 20pt

 \clearpage

 \vskip 20pt

 \clearpage

 \vskip 20pt

 \clearpage

 \vskip 20pt

\begin{figure}[h!]
    \psfrag{...}{\LARGE $\mathbf{\vdots}$}
    \psfrag{2i1+1}{\LARGE $\mathbf{2i\sb{1}+1}$}
    \psfrag{2i2+1}{\LARGE $\mathbf{2i\sb{2}+1}$}
    \psfrag{2i3+1}{\LARGE $\mathbf{2i\sb{3}}$}
    \psfrag{a4}{\LARGE $4b-3a$}
    \psfrag{1}{\LARGE $I$}
    \psfrag{3}{\LARGE $III$}
    \centerline{\scalebox{0.33}{\includegraphics{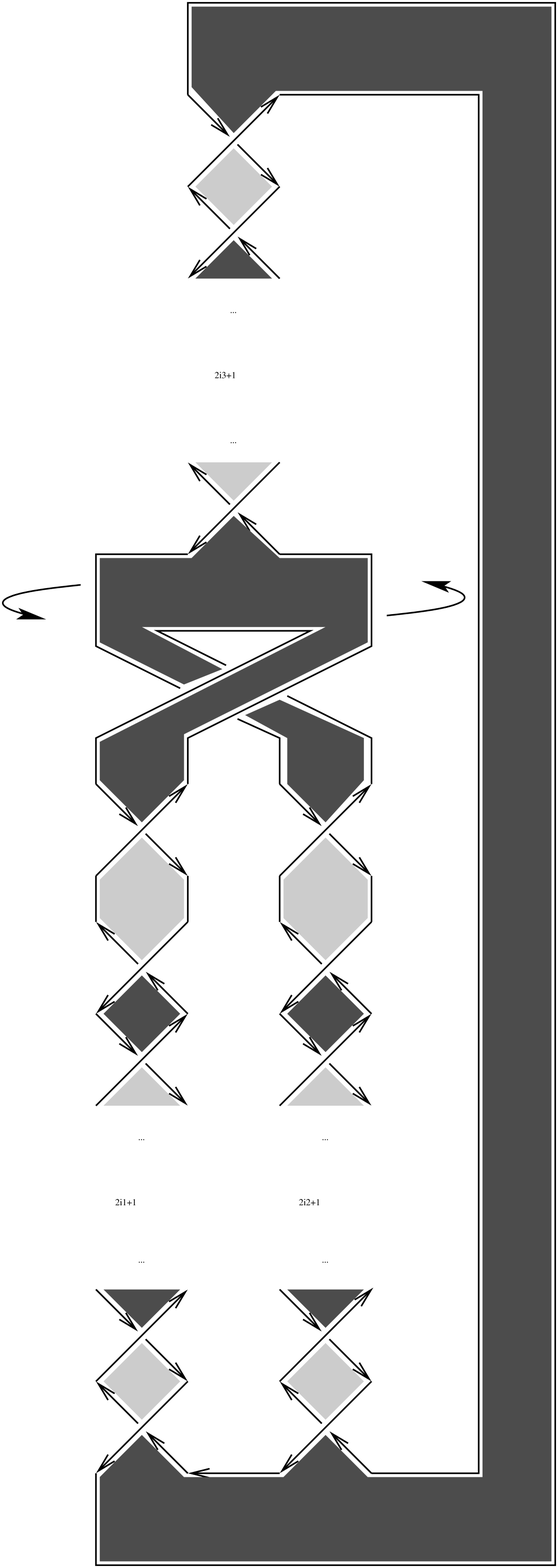}}}
    \caption{A Seifert surface of $P(2i\sb{1}+1, 2i\sb{2}+1,
2i\sb{3}+1)$ (4)}\label{Fi:pooo6a}
\end{figure}

 \clearpage

 \vskip 20pt

\begin{figure}[h!]
    \psfrag{...}{\LARGE $\mathbf{\vdots}$}
    \psfrag{2i1+1}{\LARGE $\mathbf{2i\sb{1}+1}$}
    \psfrag{2i2+1}{\LARGE $\mathbf{2i\sb{2}+1}$}
    \psfrag{2i3+1}{\LARGE $\mathbf{2i\sb{3}-1}$}
    \psfrag{a4}{\LARGE $4b-3a$}
    \psfrag{1}{\LARGE $I$}
    \psfrag{3}{\LARGE $III$}
    \centerline{\scalebox{0.32}{\includegraphics{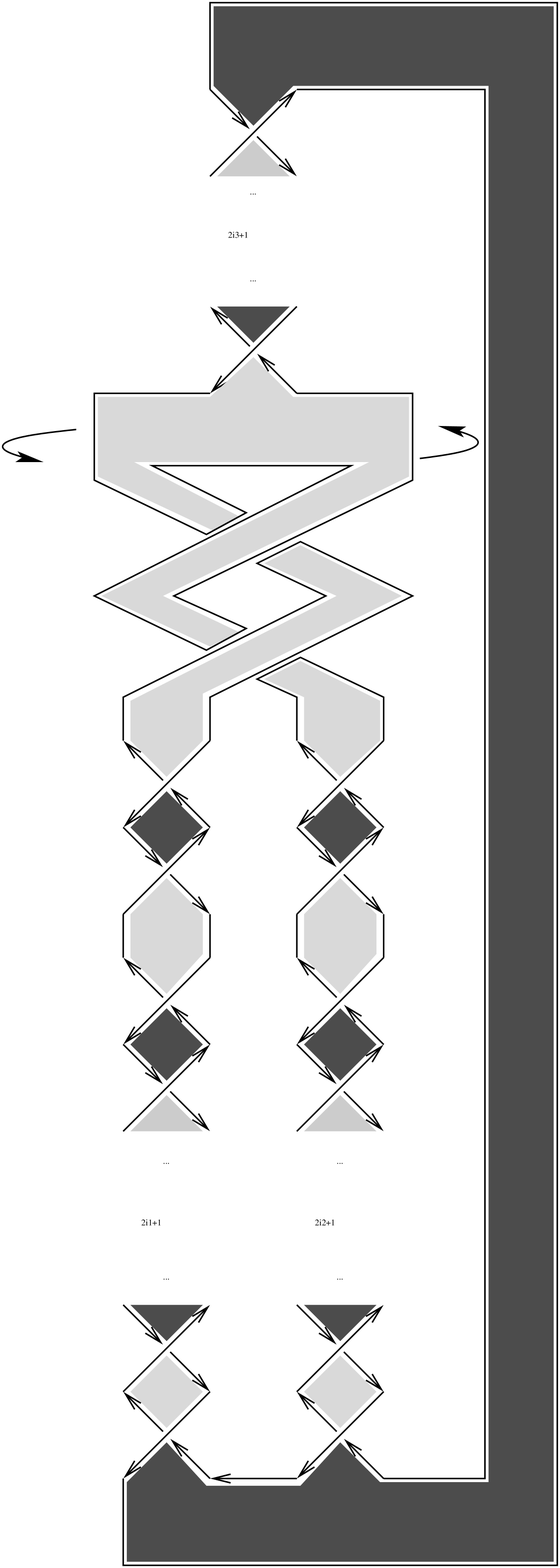}}}
    \caption{A Seifert surface of $P(2i\sb{1}+1, 2i\sb{2}+1,
2i\sb{3}+1)$ (5)}\label{Fi:pooo7a}
\end{figure}

 \clearpage

 \vskip 20pt

 \clearpage

 \vskip 20pt

 \clearpage

 \vskip 20pt

 \clearpage

 \vskip 20pt

\bigbreak

\subsection{Calculations: pretzel knots on $N$ tassels, each with an odd number of
crossings}\label{subsect:poo...o}

\noindent

We will now consider a pretzel knot on $N$ tassels, each with an
odd number of crossings. The next Proposition will distinguish two
types in this class of knots.

\begin{prop} Consider a pretzel knot on $N$ tassels, each with an
odd number of crossings. If $N$ is odd this is a $1$-component
knot. If $N$ is even this is a $2$-component knot.
\end{prop}Proof: Since each tassel has an odd number of crossings,
then starting from the top right (respect., left) side of the
tassel and going along the knot down this tassel we end up at the
bottom left (respect., right). Analogously when going along the
knot and up each tassel. So, suppose the orientation of the knot
takes us down from the top left of the $N$-th tassel to the bottom
right of this tassel, which is the bottom left of the $(N-1)$-th
tassel. Then we will be taken up along the $(N-1)$-th tassel to
its top left. Iterating this procedure we will get to the first
tassel after going through each one only once. If $N$ is odd then
we will get to the bottom right of the first tassel which is the
bottom left of the $N$-th tassel. A second trip along the knot
will take us to the the top right of the first tassel which is the
top left of the $N$-th tassel. In this way we went along the
entire knot and came back to the starting point for odd $N$. If
$N$ is even, then after going once along the knot, going through
each tassel only once, we will get to the top right of the first
tassel, which is the top left of the $N$-th tassel. Then we will
have reached the starting point without having gone through the
entire knot. If we now start at the top right of the $N$-th
tassel, for even $N$, we will go through the second component of
this knot. This completes the proof. $\hfill \blacksquare$

\bigbreak

\begin{cor} A pretzel knot with an even number of tassels, each
with an odd number of crossings has four possible distinct
orientations of its two components.
\end{cor}Proof: Omitted. $\hfill\blacksquare$

\bigbreak

\begin{def.} We will always choose an orientation of the pretzel
knot such that, in the first tassel, the over-arc of the first
crossing is oriented downwards and the over-arc of the second
crossing is oriented upwards.
\end{def.}

\bigbreak

Consider now a pretzel knot on $N$ tassels, each tassel with an
odd number of crossings, $P(2i\sb{1}+1, 2i\sb{2}+1, 2i\sb{3}+1,
\dots , 2i\sb{N}+1)$. Starting from the defining diagram, we move
the $N$-th tassel upwards as we did in the $N=3$ instance. Then we
undo the twists in this tassel by producing a braiding of the
remaining tassels and increasing the number of crossings in each
of the remaining tassels by the original number of crossings in
the $N$-th tassel (see case $N=2$). After this procedure we obtain
a Seifert surface in standard form. Each of the $N-1$ braided
ribbons represents a standard generators of the first homology
group of this Seifert surface. Moreover, the linking numbers
between standard generators and their translates are analogous to
those calculated for the $N=3$ case. As a matter of fact, in the
braiding of these $(N-1)$ ribbons, the braiding of two distinct of
them is like the braiding of ribbons (tassels) $1$ and $2$ in case
$N=3$ (just remove the other ribbons (tassels) in case $N>3$). In
this way the linking numbers are, for $j, k = 1, \dots , N-1$,

\begin{equation*}
lk(l\sb{j}, l\sb{k}\sp{+})=
\begin{cases}
i\sb{N}+1,  & \text{ if } j<k\\
i\sb{j}+i\sb{N}+1, & \text{ if } j=k\\
i\sb{N},  & \text{ if } j>k
\end{cases}
\end{equation*}

\bigbreak

The Seifert matrix is then the following $(N-1)\times (N-1)$
matrix

\bigbreak

\begin{equation*}S= \left(
\begin{matrix}
i\sb{1}+i\sb{N}+1 & i\sb{N}+1  & i\sb{N}+1  & \dots  & i\sb{N}+1\\
i\sb{N} & i\sb{2}+i\sb{N}+1 &  i\sb{N}+1  & \dots  & i\sb{N}+1\\
i\sb{N} & i\sb{N} &  i\sb{3}+i\sb{N}+1  & \dots  & i\sb{N}+1\\
\vdots        &      \vdots    &       \vdots         & \ddots &
\vdots     \\
i\sb{N} & i\sb{N} & \dots  & i\sb{N} &  i\sb{N-1}+i\sb{N}+1
\end{matrix}
\right)
\end{equation*}

\bigbreak

The presentation matrix of the Alexander module of $P(2i\sb{1}+1,
2i\sb{2}+1, 2i\sb{3}+1, \dots , 2i\sb{N}+1)$ is then

\bigbreak

\scalebox{.86}{$tS-S\sp{T}=$}

\vskip 10pt

\scalebox{.81}{$ \left(
\begin{matrix}
\bigl( i\sb{1}+i\sb{N}+1\bigr) (t-1) & \bigl( i\sb{N}+1\bigr) t-i\sb{N}   & \bigl( i\sb{N}+1\bigr) t-i\sb{N}  & \dots  & \bigl( i\sb{N}+1\bigr) t-i\sb{N} & \bigl( i\sb{N}+1\bigr) t-i\sb{N}\\
i\sb{N}t-\bigl( i\sb{N}+1\bigr)  & \bigl( i\sb{2}+i\sb{N}+1\bigr) (t-1) &  \bigl( i\sb{N}+1\bigr) t-i\sb{N}  & \dots  & \bigl( i\sb{N}+1\bigr) t-i\sb{N} & \bigl( i\sb{N}+1\bigr) t-i\sb{N}\\
i\sb{N}t-\bigl( i\sb{N}+1\bigr)   & i\sb{N}t-\bigl( i\sb{N}+1\bigr)   &  \bigl( i\sb{3}+i\sb{N}+1\bigr) (t-1)  & \dots  & \bigl( i\sb{N}+1\bigr) t-i\sb{N} & \bigl( i\sb{N}+1\bigr) t-i\sb{N}\\
\vdots        &      \vdots    &       \vdots         & \ddots  &
\vdots  &  \vdots      \\
i\sb{N}t-\bigl( i\sb{N}+1\bigr)   & i\sb{N}t-\bigl(
i\sb{N}+1\bigr)   & i\sb{N}t-\bigl( i\sb{N}+1\bigr)   & \dots   &  \bigl( i\sb{N-2}+i\sb{N}+1\bigr) (t-1)    & \bigl( i\sb{N}+1\bigr) t-i\sb{N} \\
i\sb{N}t-\bigl( i\sb{N}+1\bigr)   & i\sb{N}t-\bigl(
i\sb{N}+1\bigr)   & i\sb{N}t-\bigl( i\sb{N}+1\bigr)   & \dots  &
i\sb{N}t-\bigl( i\sb{N}+1\bigr)   & \bigl(
i\sb{N-1}+i\sb{N}+1\bigr) (t-1)
\end{matrix}
\right) $}

\bigbreak

By subtracting the second column from the first column, the third
from the second, ... , the $N$-th from the $(N-1)$-th we obtain
the equivalent matrix

\bigbreak

\begin{equation*} \left(
\begin{matrix}
i\sb{1}t-(i\sb{1}+1) & 0   & 0  & \dots  & 0   &   (i\sb{N}+1) t-i\sb{N}\\
i\sb{2}-(i\sb{2}+1)t   & i\sb{2}t-(i\sb{2}+1)  &  0  & \dots  &   0   &   (i\sb{N}+1) t-i\sb{N}\\
0   &   i\sb{3}-(i\sb{3}+1)t &  i\sb{3} t-(i\sb{3}+1)  & \dots  &   0   &   (i\sb{N}+1) t-i\sb{N}\\
\vdots        &      \vdots    &       \vdots         & \ddots &   \vdots   &   \vdots     \\
0   & 0   &   0   & \dots  & i\sb{N-2}t-(i\sb{N-2}+1)  &
(i\sb{N}+1) t-i\sb{N}  \\
0   & 0   &   0   & \dots  & i\sb{N-1}-( i\sb{N-1 }+1) t   &
(i\sb{N-1}+i\sb{N}+1) (t-1)
\end{matrix}
\right)
\end{equation*}

\bigbreak

By doing Laplace expansion over the last column we obtain the
following expression for $\det \bigr( tS-S\sp{T}\bigl) $, the
Alexander Polynomial of $P(2i\sb{1}+1, 2i\sb{2}+1, \dots ,
2i\sb{N}+1)$:

\bigbreak

\[
\bigl( (i\sb{N}+1)t-i\sb{N}\bigr) \cdot
\sum\sb{k=1}\sp{N-2}(-1)\sp{k}\cdot \prod\sb{j=1}\sp{k-1}\bigl(
i\sb{j}t-(i\sb{j}+1)   \bigr) \cdot \prod\sb{j=k+1}\sp{N-2}\bigl(
i\sb{j}-(i\sb{j}+1)t  \bigr) +
\]
\[
+ (-1)\sp{N-1}( i\sb{N-1}+i\sb{N}+1) (t-1) \cdot
\prod\sb{k=1}\sp{N-2}\bigl( i\sb{k}t-(i\sb{k}+1) \bigr) +
\]

\bigbreak

where a $(-1)\sp{N-1}$ factor has been omitted throughout.

\subsection{Calculations: pretzel knots on four tassels, exactly one with an even number of
crossings}\label{subsect:pooe}

\noindent

In this Subsection we work out the calculations for a Pretzel knot
on four tassels with exactly one tassel with an even number of
crossings. This is intended to pave the way for the calculations
for a pretzel knot with an arbitrary even number of tassels and
with exactly one tassel with an even number of crossings which
will done in the next Subsection.

\bigbreak

\begin{figure}[h!]
    \centerline{\scalebox{0.6}{\includegraphics{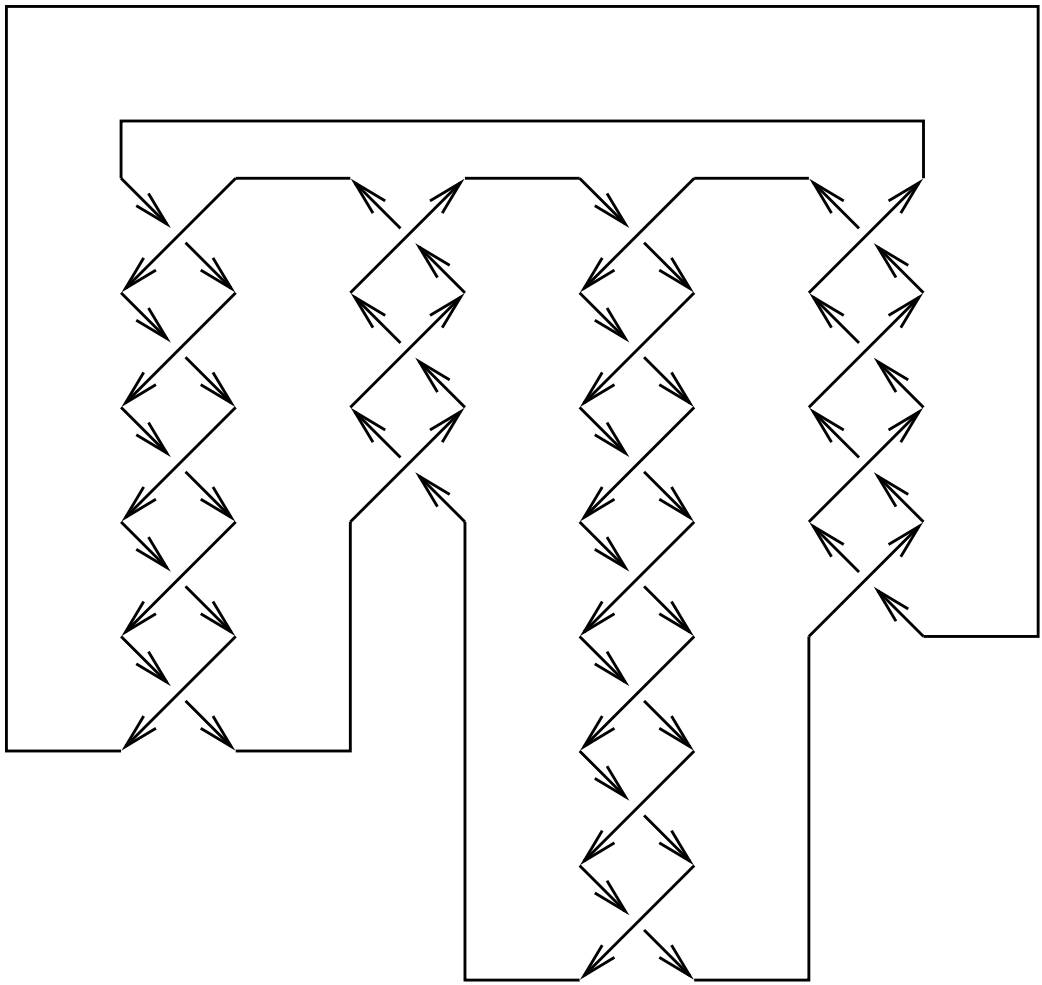}}}
    \caption{A diagram of $P(5, 3, 7, 4)$}\label{Fi:p5374}
\end{figure}

\begin{figure}[h!]
    \centerline{\scalebox{0.6}{\includegraphics{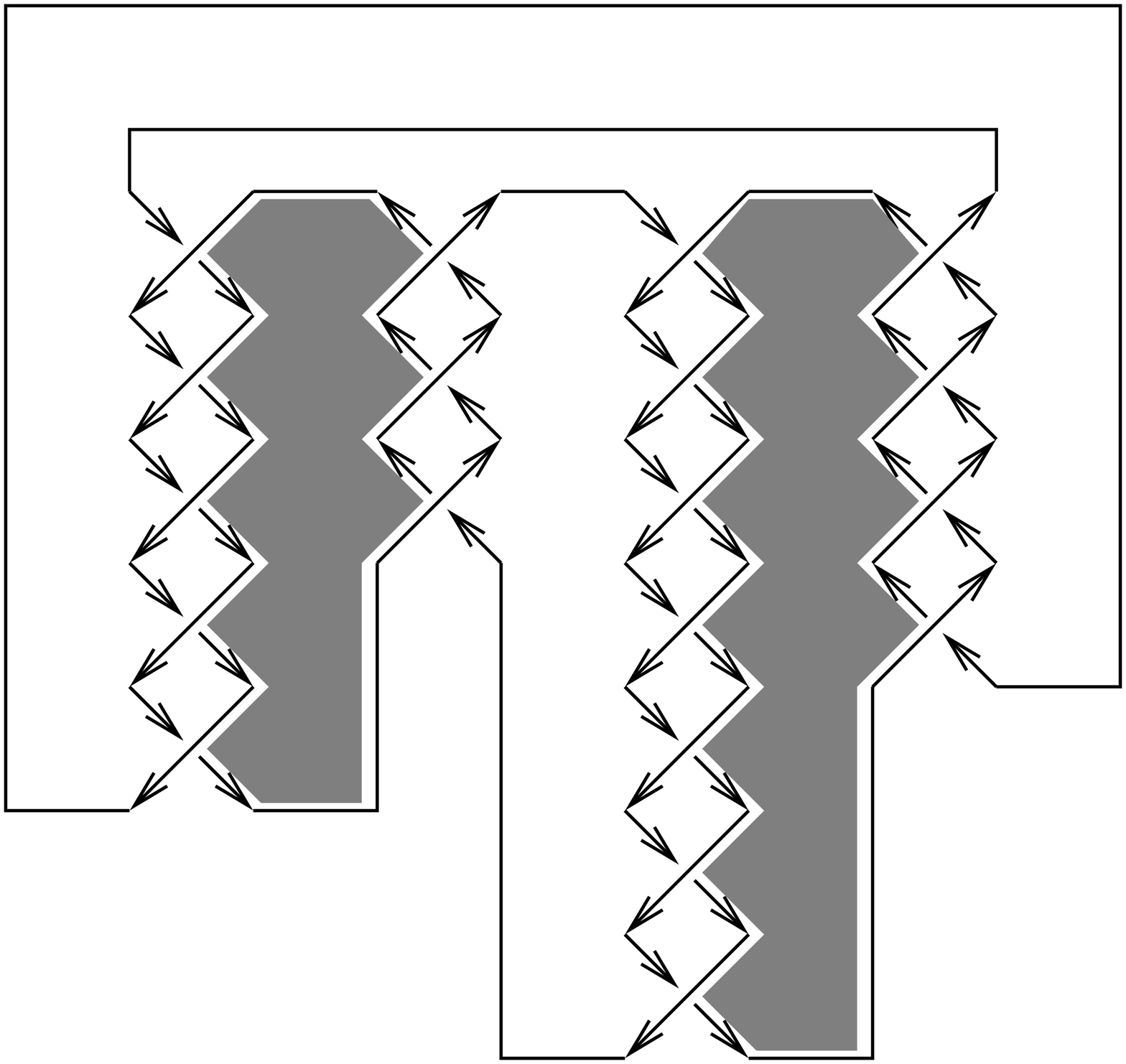}}}
    \caption{A Seifert surface for $P(5, 3, 7, 4)$ (1)}\label{Fi:seifp5374}
\end{figure}

We choose pretzel knot $P(5, 3, 7, 4)$ to do our calculations (see
Figure \ref{Fi:p5374}). A Seifert surface for this pretzel knot is
obtained as before using the algorithm in Proposition
\ref{prop:seifertalgo} (see Figure \ref{Fi:seifp5374}). Only one
side of the Seifert surface will be shaded in order not to blur
the interpretation of the Figures later on. We remark that on the
shaded side of the Seifert surface the normal points to the
reader. Our strategy here will be first to transfer the twists in
the ribbons from left to right, by deformations; the twists will
ultimately accumulate in the ribbons to the far right. This will,
however, produce a braiding of the ribbons. Then we will shrink
some of these ribbons in order to merge the center discs into the
larger disc. We will finally obtain a Seifert surface in standard
form. We now specify each step of this procedure.

\bigbreak

\begin{figure}[h!]
    \centerline{\scalebox{0.4}{\includegraphics{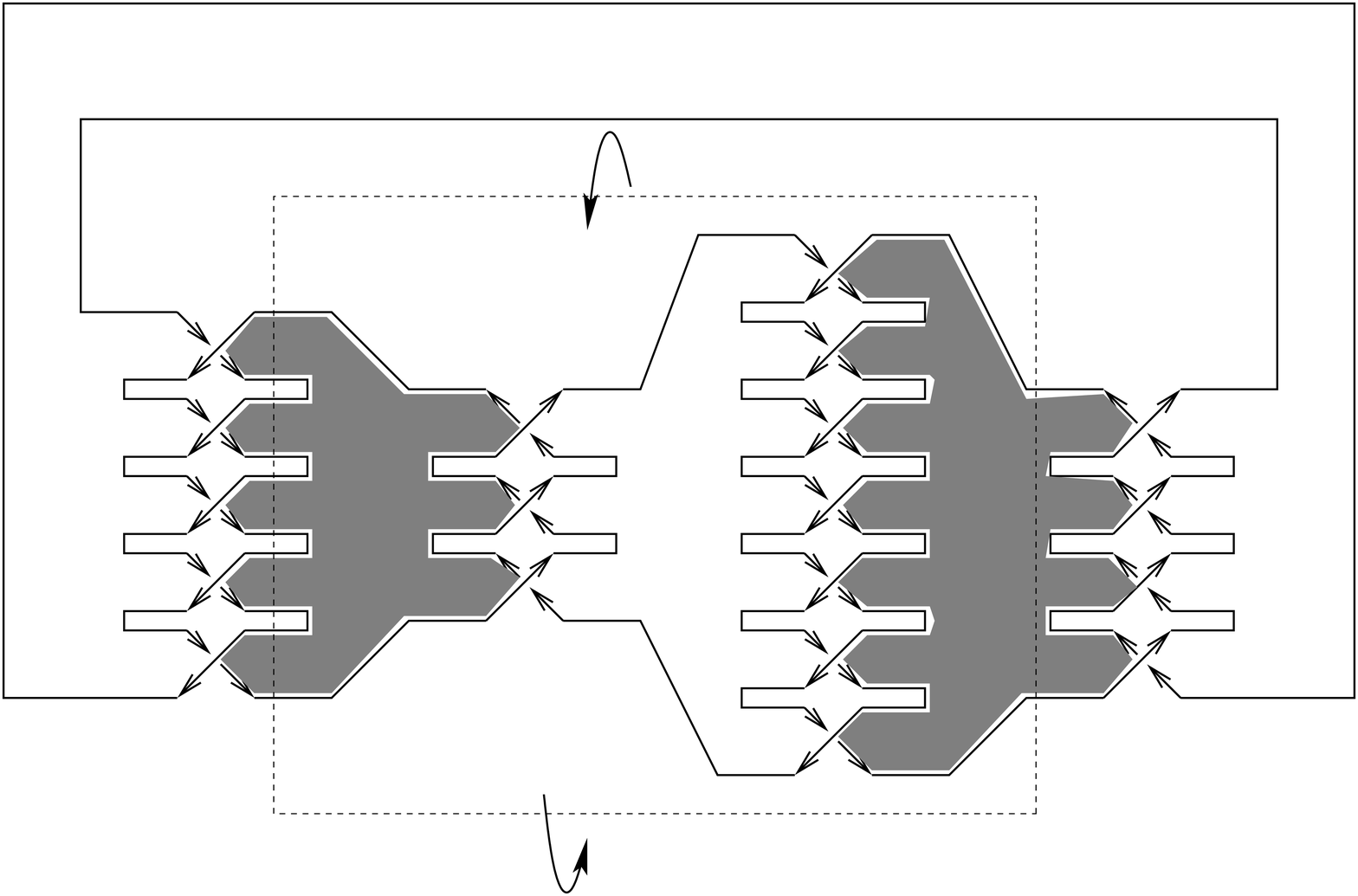}}}
    \caption{A Seifert surface for $P(5, 3, 7, 4)$ (2)}\label{Fi:seif2p5374}
\end{figure}

\begin{figure}[h!]
    \centerline{\scalebox{0.6}{\includegraphics{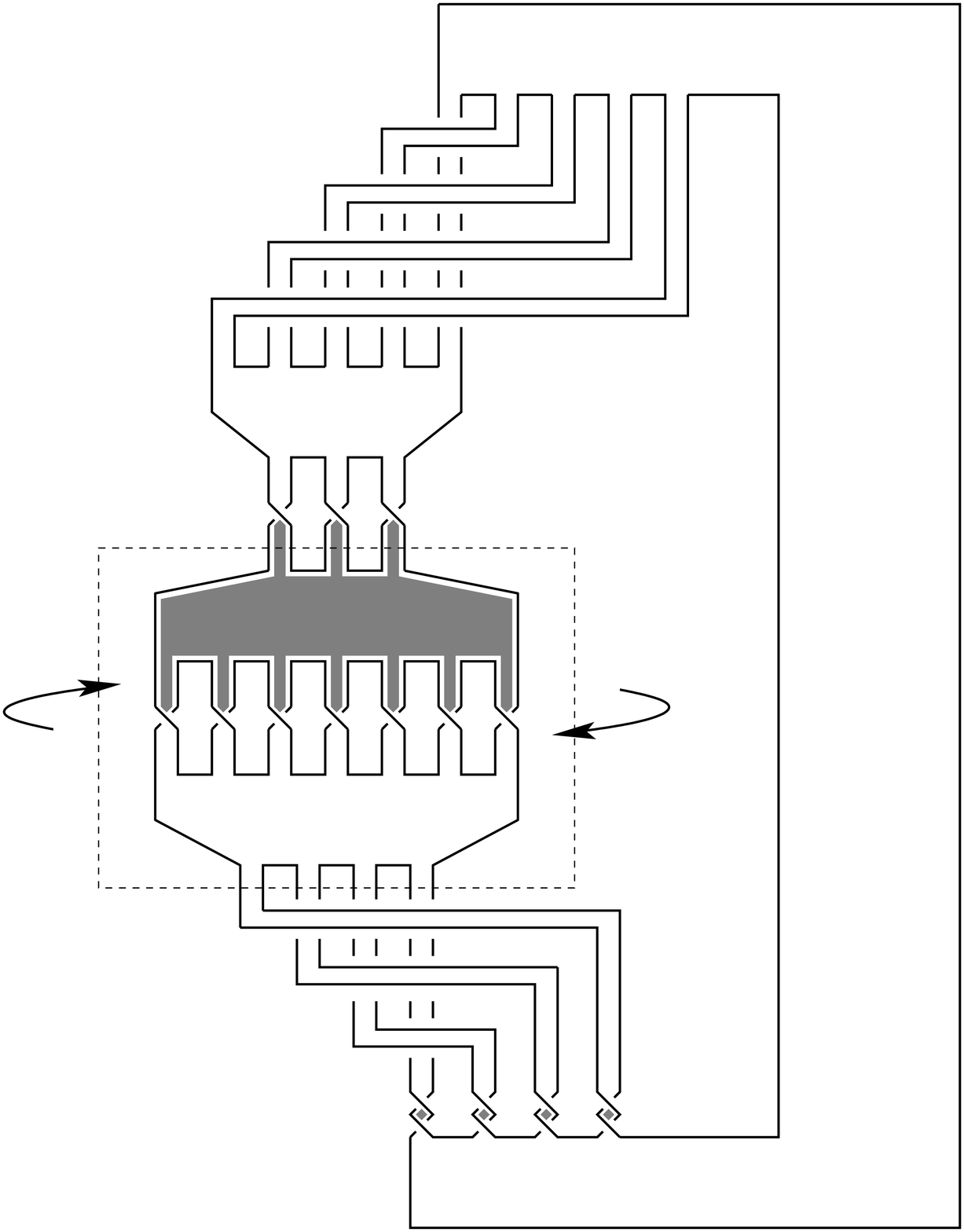}}}
    \caption{A Seifert surface for $P(5, 3, 7, 4)$ (3)}\label{Fi:seif3p5374}
\end{figure}

\begin{figure}[h!]
    \centerline{\scalebox{0.49}{\includegraphics{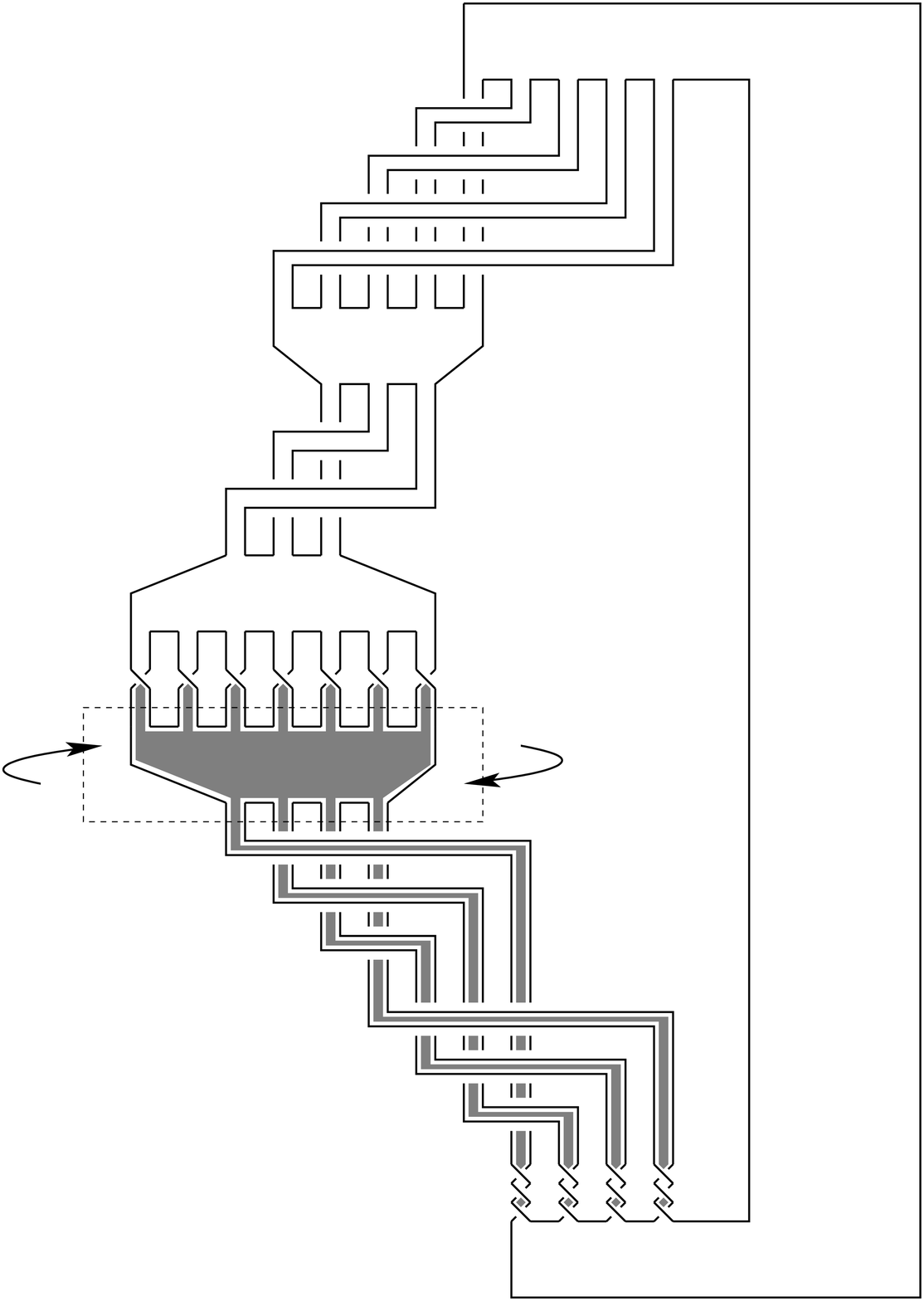}}}
    \caption{A Seifert surface for $P(5, 3, 7, 4)$ (4)}\label{Fi:seif4p5374}
\end{figure}

\begin{figure}[h!]
    \psfrag{A}{\huge $A$}
    \psfrag{B}{\huge $B$}
    \psfrag{C}{\huge $C$}
    \centerline{\scalebox{0.34}{\includegraphics{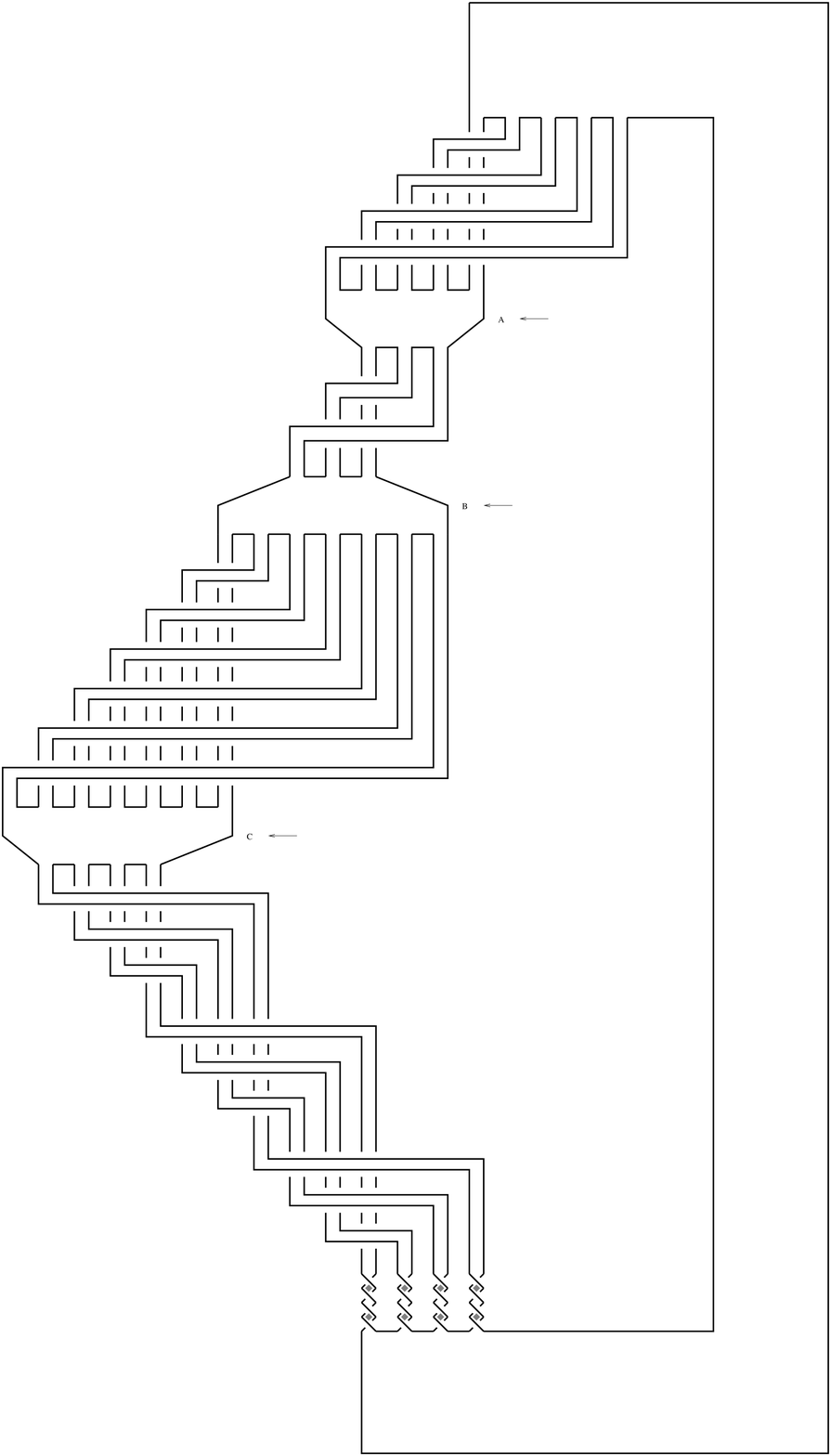}}}
    \caption{A Seifert surface for $P(5, 3, 7, 4)$ (5)}\label{Fi:seif5p5374}
\end{figure}

Figure \ref{Fi:seif2p5374} is obtained by stretching the central
portion of Figure \ref{Fi:seifp5374} so that the ribbons
connecting the different discs are clearly identified. Then the
twists in the ribbons to the left in Figure \ref{Fi:seif2p5374}
are transferred to the right by a $\frac{\pi}{2}$ rotation  of the
portion of the diagram boxed by the dotted rectangle. The result
of this rotation is shown in Figure \ref{Fi:seif3p5374}; this
Figure is further obtained from Figure \ref{Fi:seif2p5374} by
rotation of the whole by a $\frac{\pi}{2}$ clockwise on the plane
of the page. The upper twists in ribbons in Figure
\ref{Fi:seif3p5374} are transferred down by $\frac{\pi}{2}$
rotation of the portion of the diagram boxed by the dotted
rectangle is indicated. The result of this rotation is shown in
Figure \ref{Fi:seif4p5374}. The upper twists in ribbons in Figure
\ref{Fi:seif4p5374} are transferred down using the same technique.
The result is shown in Figure \ref{Fi:seif5p5374} where the twists
in ribbons are all in the lower portion of this Figure. As of this
Figure we proceed to merge the smaller discs into the larger one
mainly by shrinking specific ribbons.

\bigbreak

\begin{figure}[h!]
    \psfrag{x}{\huge $x$}
    \psfrag{A}{\huge $A$}
    \psfrag{B}{\huge $B$}
    \psfrag{C}{\huge $C$}
    \centerline{\scalebox{0.25}{\includegraphics{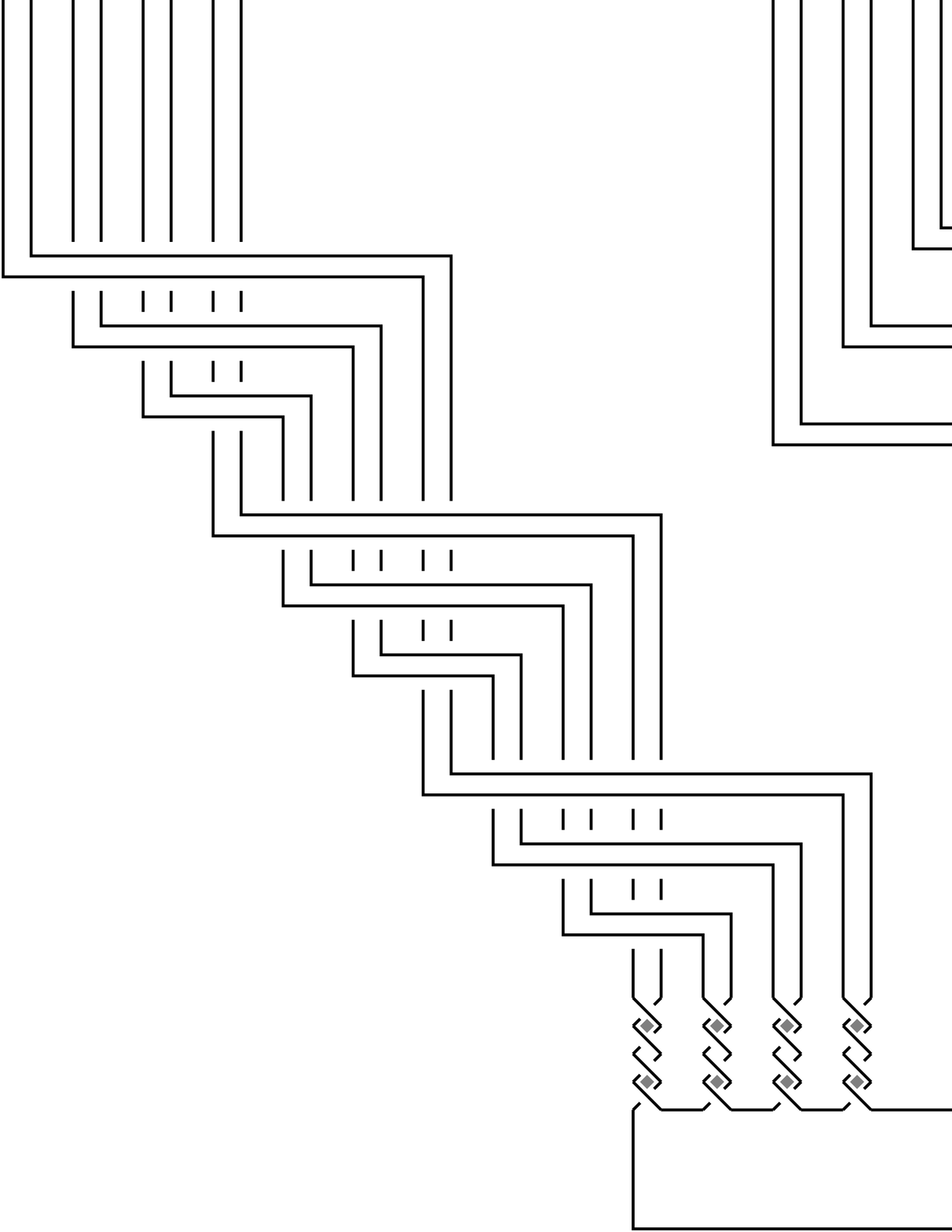}}}
    \caption{A Seifert surface for $P(5, 3, 7, 4)$ (6)}\label{Fi:seif6p5374}
\end{figure}

\begin{figure}[h!]
    \psfrag{x}{\huge $x$}
    \psfrag{A}{\huge $A$}
    \psfrag{B}{\huge $B$}
    \psfrag{C}{\huge $C$}
    \centerline{\scalebox{0.25}{\includegraphics{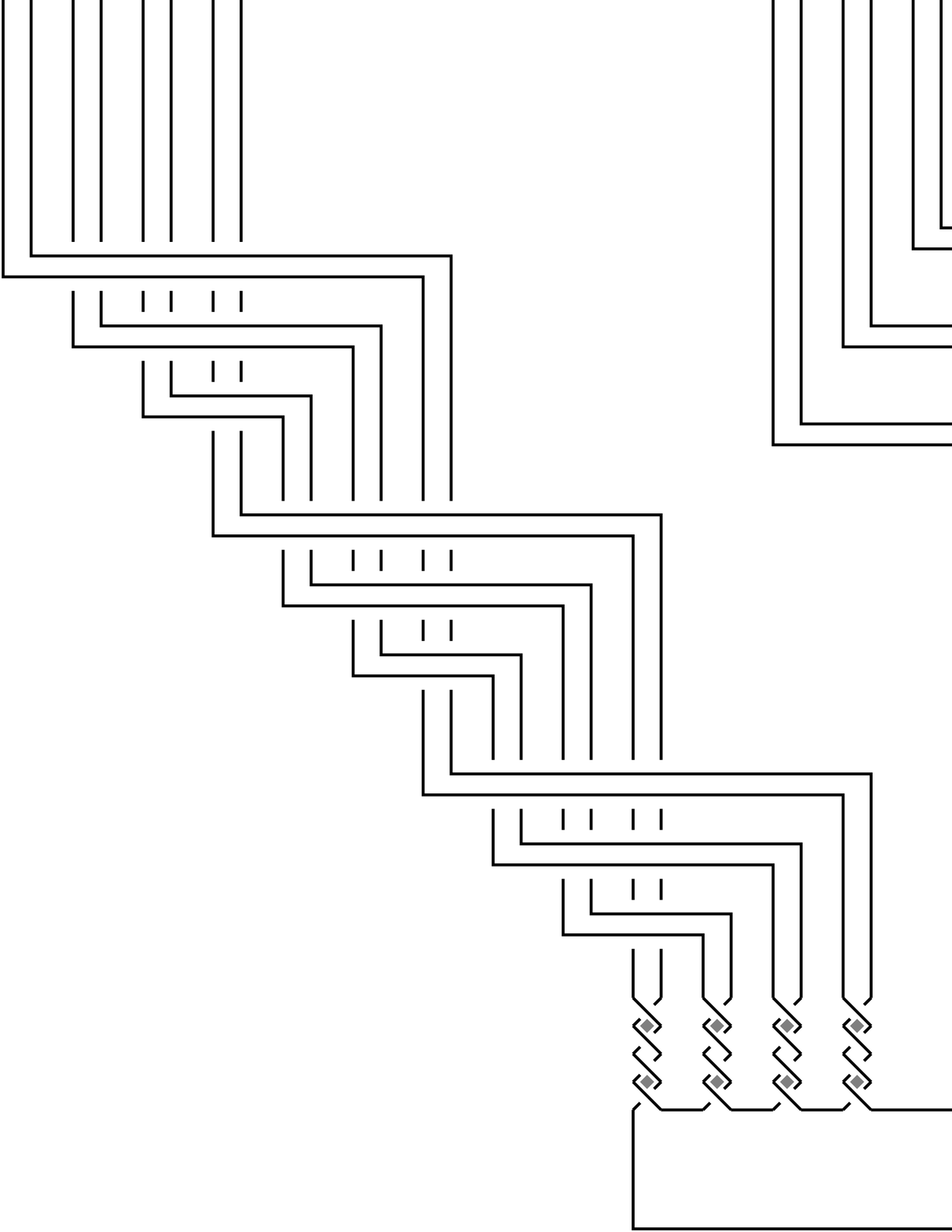}}}
    \caption{A Seifert surface for $P(5, 3, 7, 4)$ (7)}\label{Fi:seif7p5374}
\end{figure}

\begin{figure}[h!]
    \psfrag{x}{\huge $x$}
    \psfrag{A}{\huge $A$}
    \psfrag{B}{\huge $B$}
    \psfrag{C}{\huge $C$}
    \centerline{\scalebox{0.28}{\includegraphics{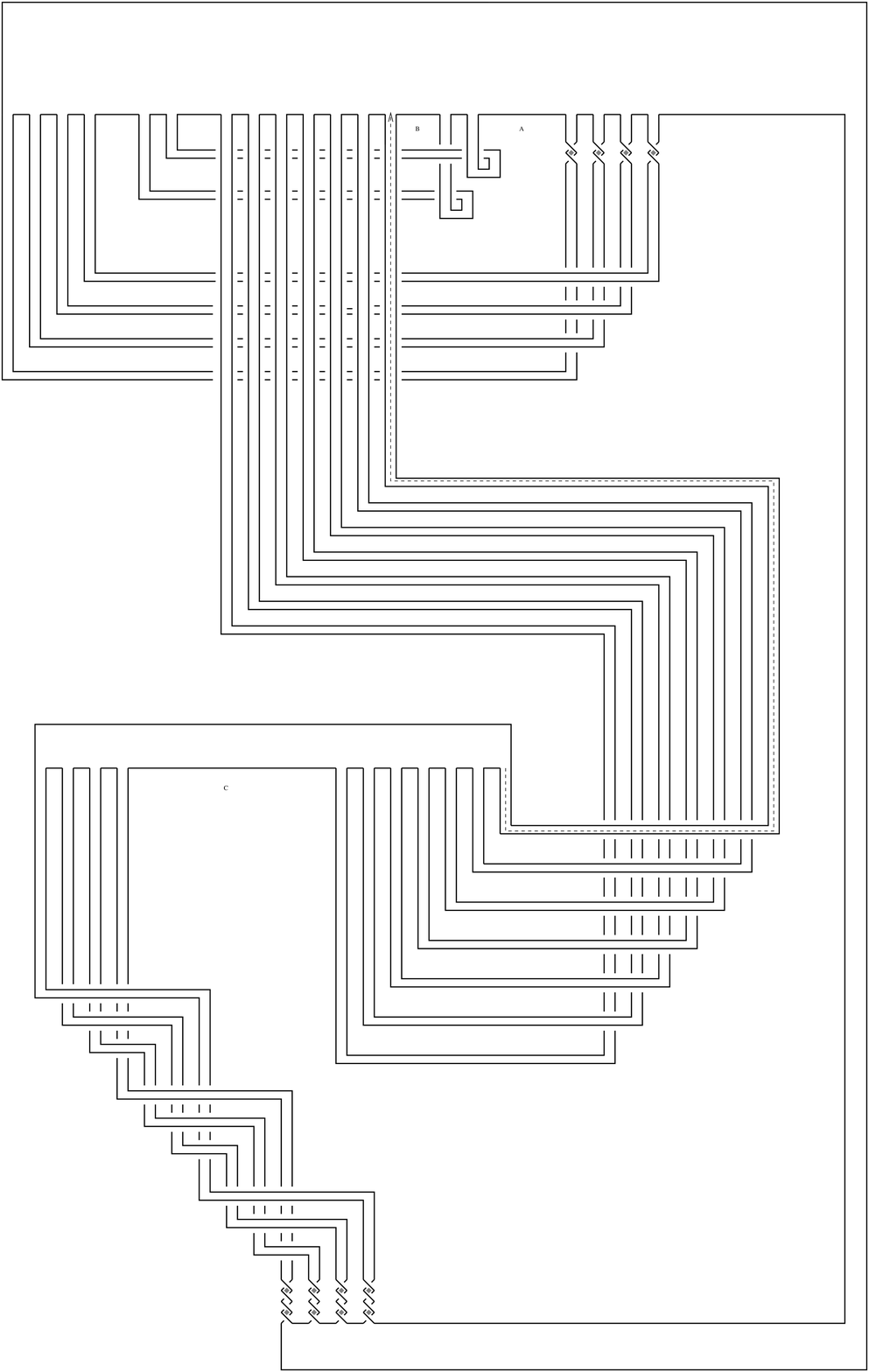}}}
    \caption{A Seifert surface for $P(5, 3, 7, 4)$ (8)}\label{Fi:seif8p5374}
\end{figure}

\begin{figure}[h!]
    \psfrag{x}{\huge $x$}
    \psfrag{A}{\huge $A$}
    \psfrag{B}{\huge $B$}
    \psfrag{C}{\huge $C$}
    \centerline{\scalebox{0.23}{\includegraphics{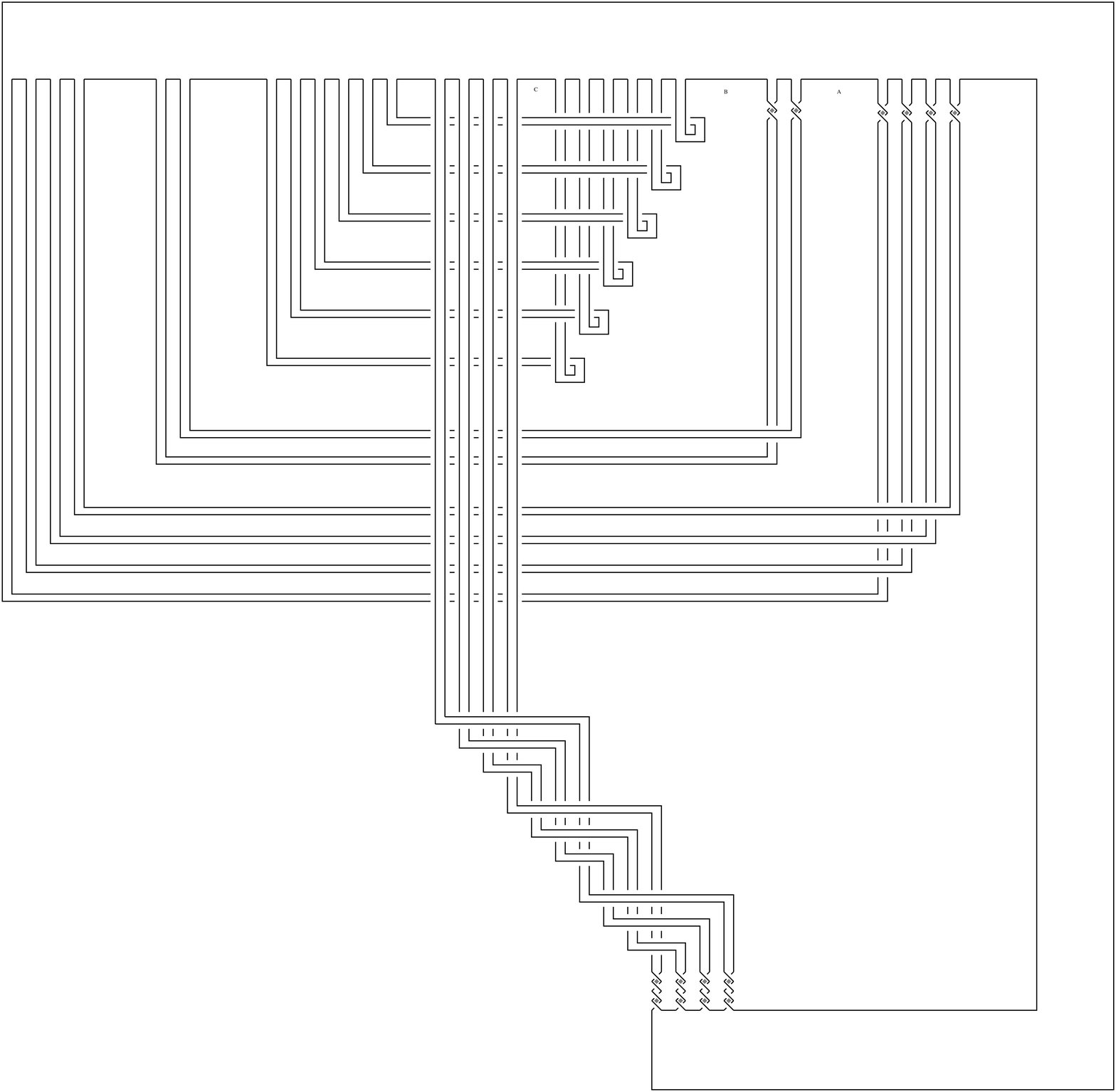}}}
    \caption{A Seifert surface for $P(5, 3, 7, 4)$ (9)}\label{Fi:seif9p5374}
\end{figure}

\begin{figure}[h!]
    \psfrag{x}{\huge $x$}
    \psfrag{A}{\huge $A$}
    \psfrag{B}{\huge $B$}
    \psfrag{C}{\huge $C$}
    \psfrag{x}{\huge $x$}
    \psfrag{1}{\huge $1$}
    \psfrag{2}{\huge $2$}
    \psfrag{3}{\huge $3$}
    \psfrag{4}{\huge $4$}
    \psfrag{5}{\huge $5$}
    \psfrag{6}{\huge $6$}
    \psfrag{7}{\huge $7$}
    \psfrag{8}{\huge $8$}
    \psfrag{9}{\huge $9$}
    \psfrag{10}{\huge $10$}
    \psfrag{11}{\huge $11$}
    \psfrag{12}{\huge $12$}
    \psfrag{13}{\huge $13$}
    \psfrag{14}{\huge $14$}
    \psfrag{15}{\huge $15$}
    \psfrag{16}{\huge $16$}
    \centerline{\scalebox{0.23}{\includegraphics{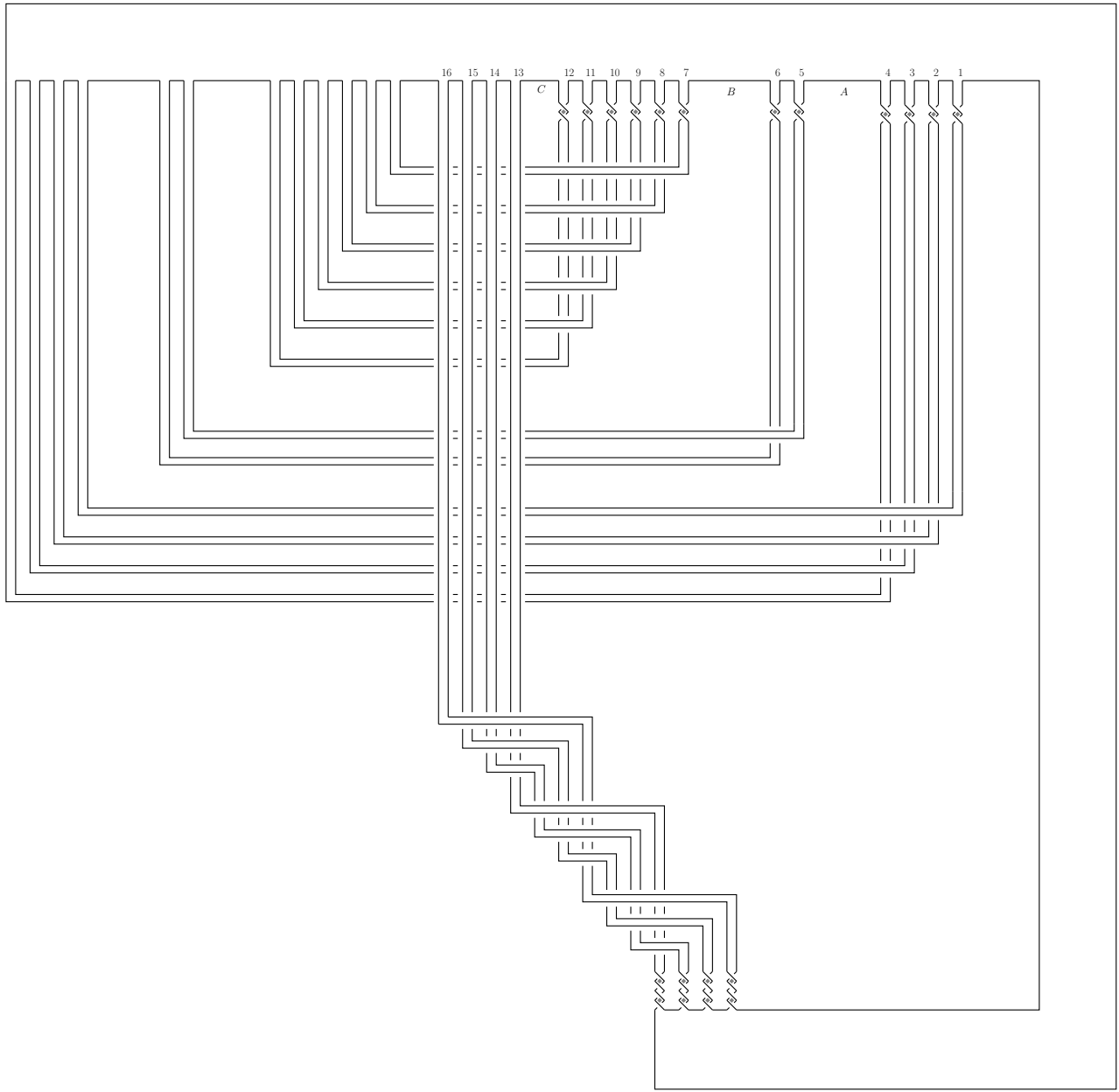}}}
    \caption{A Seifert surface for $P(5, 3, 7, 4)$ - the standard form}\label{Fi:seif10p5374}
\end{figure}

In Figure \ref{Fi:seif5p5374} we identify three vertices of the
diagram by the letters $A$, $B$, and $C$. These are the vertices
to the immediate left of each of these letters. The arrows to the
right of these letters indicate that a deformation of the surface
will be performed by, also, pushing on the corresponding vertices.
The same letters in Figure \ref{Fi:seif6p5374} indicate the new
location of the corresponding points after the deformation.
Furthermore, in this Figure, there is a dotted and oriented arc
inside a ribbon, in the upper part of the diagram. This indicates
that the corresponding ribbon will be shrunk in the indicated
direction thus merging the disc about point $A$ into the larger
disc. The result of this operation is shown in Figure
\ref{Fi:seif7p5374}. In this Figure, the dotted and oriented arc
inside a ribbon indicates the same operation on the corresponding
ribbon. The result of this operation is shown in Figure
\ref{Fi:seif8p5374}. Without further remarks, in Figure
\ref{Fi:seif9p5374} we show the result of the merging of the last
small disc into the larger one. After some deforming we obtain a
standard form in Figure \ref{Fi:seif10p5374}. Here we enumerate
the ribbons as indicated in the top of the diagram. We thus obtain
an enumeration of the standard generators of the first  homology
group of this surface (see Definition \ref{def:seifstandardgen}).
We orient these generators counterclockwise.

\bigbreak

Consider again the Seifert surface in standard form depicted in
Figure \ref{Fi:seif10p5374}. Since each ribbon is twisted, each
standard generator links with its translate. Thus, every diagonal
element of the Seifert matrix will be non-null. Except for the
four generators with numbers $13, 14, 15$ and  $16$ which
potentially link with every other generator, the remaining
generators exhibit linking within given subsets of generators. It
is further clear that these subsets of generators are directly
related to definite tassels in $P(5, 3, 7, 4)$. Specifically, the
four generators with numbers $13, 14, 15, 16$, stem from the only
tassel with the even number of crossings (four crossings). The
subset of generators $7, 8, 9, 10, 11, 12$ stems from the tassel
with seven crossings; the subset of generators $5, 6$ stems from
the tassel with three crossings; and the subset of generators $1,
2, 3, 4$ stems from the tassel with five crossings. In the general
case of a pretzel knot with (even) $N$ tassels exactly one with an
even number of crossings, there will be, after the procedure
described before, a subset of standard generators for the first
homology group with as many elements as there are crossings in the
tassel these generators stem from (we call this tassel the
reference tassel); the remaining tassels will generate subsets of
standard generators whose cardinality is the number of crossings
of the tassel minus one. The cause for this ``minus one''
difference has to do with the shrinking of one of the ribbons in
the passages from Figure \ref{Fi:seif6p5374} to Figure
\ref{Fi:seif7p5374}, \ref{Fi:seif7p5374} to \ref{Fi:seif8p5374},
and \ref{Fi:seif8p5374} to \ref{Fi:seif9p5374}. It is clear that
the reference tassel should be chosen to be the one with least
number of crossings since it is the generators stemming from this
tassel that link with each of the other generators.

\bigbreak

Finally, invoking Propositions \ref{prop : lk(i,i)}
 and \ref{prop : lk(i,j)} we obtain the following Seifert matrix.



\begin{equation*}
\setcounter{MaxMatrixCols}{16} S=\left(
\begin{matrix}
-1 &  0   &          0        &    0     &   0     &   0       & 0
&     0   &          0        &    0     &   0     &   0       & 1
   &   1     &   1       &   1   \\
1   &    -1 &  0   &          0        &    0     &   0     &   0
& 0 &     0   &         0
    &   0     &   0       &   1   &   1     &   1       &   1   \\
1   &    1 &  -1   &          0        &    0     &   0     &   0
& 0 &     0   &            0
    &   0     &   0       &   1   &   1     &   1       &   1   \\
1   &    1 &  1   &          -1        &    0     &   0     &   0
& 0 &     0   &            0
    &   0     &   0       &  1   &   1     &   1       &   1   \\
0   &    0 &  0   &          0        &    -1     &   0     &   0
& 0 &     0   &            0
    &   0     &   0       &  1   &   1     &   1       &   1   \\
    0     &   0       &   0   &   0     &   1       &   -1     &
0 &  0   &          0        &    0     &   0     &   0       &1
 &   1     &   1       &   1   \\
    0     &   0       &   0   &   0     &   0       &   0     &
-1 &  0   &          0        &    0     &   0     &   0       &
1   &   1     &   1       &   1   \\
    0     &   0       &   0   &   0     &   0       &   0     &
1 &  -1   &          0        &    0     &   0     &   0       &
1   &   1     &   1       &   1   \\
    0     &   0       &   0   &   0     &   0       &   0     &
   1       &    1     &   -1  &  0   &         0     &   0       &
1   &   1     &   1       &   1   \\
    0     &   0       &   0   &   0     &   0       &   0   &     1
           &    1     &   1     &   -1       & 0 &     0    &
1   &   1     &   1       &   1   \\
    0     &   0       &   0   &   0     &   0       &   0     &  1
&    1     &   1     &   1       & -1  &     0   &
1   &   1     &   1       &   1   \\
    0     &   0       &   0      &   0     &
0 &  0   &          1       &    1     &   1     &   1       & 1
&     -1   &  1   &   1     &   1       &   1   \\
 0   &          0        &    0     &   0     &   0       & 0 &
 0     &   0       &   0   &   0     &   0       &   0     &
-2 &     -2   &          -2        &    -2    \\
 0   &          0        &    0     &   0     &   0       & 0 &
 0     &   0       &   0   &   0     &   0       &   0     &
-1 &     -2   &          -2        &    -2   \\
 0   &          0        &    0     &   0     &   0       & 0 &
 0     &   0       &   0   &   0     &   0       &   0     &
-1 &     -1   &          -2        &    -2   \\
 0   &          0        &    0     &   0     &   0       & 0 &
 0     &   0       &   0   &   0     &   0       &   0     &
-1 &     -1   &          -1        &    -2
\end{matrix}
\right)
\end{equation*}

\bigbreak

\bigbreak

\bigbreak

\bigbreak

\clearpage

\clearpage

\clearpage

\clearpage

\clearpage

\clearpage

\clearpage

We now calculate the presentation matrix for the Alexander module
of $P(5, 3, 7, 4)$ from the Seifert matrix above. We will further
use elementary transformations on matrices to pass on to a more
convenient presentation matrix of this module, in order to
calculate the Alexander polynomial.

\bigbreak

\scalebox{.75}{ $ \setcounter{MaxMatrixCols}{16} tS-S\sp{T}=\left(
\begin{matrix}
1-t &  -1   &          -1        &    -1     &   0     &   0 & 0 &
0   &          0       &   0     &   0       &   0   &
t     &   t      &   t & t\\
t   &    1-t &  -1   &          -1        &    0     &   0     & 0
& 0 &     0   &         0
    &   0     &   0       &   t     &   t      &   t & t\\
t   &    t &  1-t   &          -1        &    0     &   0     &
0 & 0 &     0   &            0
    &   0     &   0        & t     &   t      &   t & t \\
t   &    t &  t   &          1-t        &    0     &   0     &   0
& 0 &     0   &            0
    &   0     &   0       &  t     &   t      &   t & t \\
0   &    0 &  0   &          0        &    1-t     &   -1     &
0 & 0 &     0   &            0
    &   0     &   0       &  t     &   t      &   t & t\\
    0     &   0       &   0   &   0     &   t       &   1-t     &
0 &  0   &          0        &    0     &   0     &   0
    & t     &   t      &   t & t\\
    0     &   0       &   0   &   0     &   0       &   0     &
1-t &  -1   &          -1        &    -1     &   -1     &   -1 &
t     &   t      &   t & t\\
    0     &   0       &   0   &   0     &   0       &   0     &
t &  1-t   &         -1        &   -1     & -1     &  -1       &
t     &   t      &   t & t\\
    0     &   0       &   0   &   0     &   0       &   0     &
   t       &    t     &   1-t  &  -1   &       -1     & -1
    & t     &   t      &   t & t\\
    0     &   0       &   0   &   0     &   0       &   0   & t
 &    t     &   t     &   1-t       & -1 &     -1    &
t     &   t      &   t & t\\
    0     &   0       &   0   &   0     &   0       &   0     &  t
&    t     &   t     &   t       & 1-t  &     -1   &
t     &   t      &   t & t\\
    0     &   0       &   0      &   0     &
0 &  0   &        t       &    t     &   t     &   t       & t
&     1-t   &   t     &   t      &   t & t\\
 -1   &          -1        &    -1     &   -1     &   -1       & -1 &
 -1     &   -1       &   -1   &   -1     &   -1       &   -1     &
2-2t &     1-2t   &          1-2t        &    1-2t    \\
 -1   &          -1        &    -1     &   -1     &   -1       & -1 &
 -1     &   -1       &   -1   &   -1     &   -1       &   -1     &
2-t &     2-2t   &          1-2t        &    1-2t   \\
 -1   &          -1        &    -1     &   -1     &   -1       & -1 &
 -1     &   -1       &   -1   &   -1     &   -1       &   -1     &
2-t &     2-t   &          2-2t        &    1-2t   \\
 -1   &          -1        &    -1     &   -1     &   -1       & -1 &
 -1     &   -1       &   -1   &   -1     &   -1       &   -1     &
2-t &     2-t   &          2-t        &    2-2t
\end{matrix}
\right)
$
}

\bigbreak

The next matrix is the result of subtracting, in the preceding
matrix, the $14$-th column from the $13$-th column, the $15$-th
from the $14$-th, and the $16$-th from the $15$-th.

\bigbreak

\scalebox{.85}{ $ \setcounter{MaxMatrixCols}{16} tS-S\sp{T}
\rightarrow \left(
\begin{matrix}
1-t &  -1   &          -1        &    -1     &   0     &   0 & 0 &
0   &          0       &   0     &   0       &   0   &
0     &   0      &   0 & t\\
t   &    1-t &  -1   &          -1        &    0     &   0     & 0
& 0 &     0   &         0
    &   0     &   0       &   0     &   0      &   0 & t\\
t   &    t &  1-t   &          -1        &    0     &   0     & 0
& 0 &     0   &            0
    &   0     &   0        & 0     &   0      &   0 & t \\
t   &    t &  t   &          1-t        &    0     &   0     &   0
& 0 &     0   &            0
    &   0     &   0       &  0     &   0      &   0 & t \\
0   &    0 &  0   &          0        &    1-t     &   -1     & 0
& 0 &     0   &            0
    &   0     &   0       &  0     &   0      &   0 & t\\
    0     &   0       &   0   &   0     &   t       &   1-t     &
0 &  0   &          0        &    0     &   0     &   0
    & 0     &   0      &   0 & t\\
    0     &   0       &   0   &   0     &   0       &   0     &
1-t &  -1   &          -1        &    -1     &   -1     &   -1 &
0     &   0      &   0 & t\\
    0     &   0       &   0   &   0     &   0       &   0     &
t &  1-t   &         -1        &   -1     & -1     &  -1       &
0     &   0      &   0 & t\\
    0     &   0       &   0   &   0     &   0       &   0     &
   t       &    t     &   1-t  &  -1   &       -1     & -1
    & 0     &   0      &   0 & t\\
    0     &   0       &   0   &   0     &   0       &   0   & t
 &    t     &   t     &   1-t       & -1 &     -1    &
0     &   0      &   0 & t\\
    0     &   0       &   0   &   0     &   0       &   0     &  t
&    t     &   t     &   t       & 1-t  &     -1   &
0     &   0      &   0 & t\\
    0     &   0       &   0      &   0     &
0 &  0   &        t       &    t     &   t     &   t       & t
&     1-t   &   0     &   0      &   0 & t\\
 -1   &          -1        &    -1     &   -1     &   -1       & -1 &
 -1     &   -1       &   -1   &   -1     &   -1       &   -1     &
1 &     0   &          0        &    1-2t    \\
 -1   &          -1        &    -1     &   -1     &   -1       & -1 &
 -1     &   -1       &   -1   &   -1     &   -1       &   -1     &
t &     1   &          0        &    1-2t   \\
 -1   &          -1        &    -1     &   -1     &   -1       & -1 &
 -1     &   -1       &   -1   &   -1     &   -1       &   -1     &
0 &     t   &          1        &    1-2t   \\
 -1   &          -1        &    -1     &   -1     &   -1       & -1 &
 -1     &   -1       &   -1   &   -1     &   -1       &   -1     &
0 &     0   &          t        &    2-2t
\end{matrix}
\right) $ }

\bigbreak

The next matrix is the result of subtracting the $16$-th row from
the $13$-th, the $14$-th, and the $15$-th rows of the preceding
matrix.

\bigbreak

\scalebox{.82}{ $ \setcounter{MaxMatrixCols}{16}  tS-S\sp{T}
\rightarrow \left(
\begin{matrix}
1-t &  -1   &          -1        &    -1     &   0     &   0 & 0 &
0   &          0       &   0     &   0       &   0   &
0     &   0      &   0 & t\\
t   &    1-t &  -1   &          -1        &    0     &   0     & 0
& 0 &     0   &         0
    &   0     &   0       &   0     &   0      &   0 & t\\
t   &    t &  1-t   &          -1        &    0     &   0     & 0
& 0 &     0   &            0
    &   0     &   0        & 0     &   0      &   0 & t \\
t   &    t &  t   &          1-t        &    0     &   0     &   0
& 0 &     0   &            0
    &   0     &   0       &  0     &   0      &   0 & t \\
0   &    0 &  0   &          0        &    1-t     &   -1     & 0
& 0 &     0   &            0
    &   0     &   0       &  0     &   0      &   0 & t\\
    0     &   0       &   0   &   0     &   t       &   1-t     &
0 &  0   &          0        &    0     &   0     &   0
    & 0     &   0      &   0 & t\\
    0     &   0       &   0   &   0     &   0       &   0     &
1-t &  -1   &          -1        &    -1     &   -1     &   -1 &
0     &   0      &   0 & t\\
    0     &   0       &   0   &   0     &   0       &   0     &
t &  1-t   &         -1        &   -1     & -1     &  -1       &
0     &   0      &   0 & t\\
    0     &   0       &   0   &   0     &   0       &   0     &
   t       &    t     &   1-t  &  -1   &       -1     & -1
    & 0     &   0      &   0 & t\\
    0     &   0       &   0   &   0     &   0       &   0   & t
 &    t     &   t     &   1-t       & -1 &     -1    &
0     &   0      &   0 & t\\
    0     &   0       &   0   &   0     &   0       &   0     &  t
&    t     &   t     &   t       & 1-t  &     -1   &
0     &   0      &   0 & t\\
    0     &   0       &   0      &   0     &
0 &  0   &        t       &    t     &   t     &   t       & t
&     1-t   &   0     &   0      &   0 & t\\
 0   &          0        &    0     &   0     &   0       & 0 &
 0     &   0       &   0   &   0     &   0       &   0     &
1 &     0   &          -t        &    -1    \\
 0   &          0        &    0     &   0     &   0       & 0 &
 0     &   0       &   0   &   0     &   0       &   0     &
t &     1   &          -t        &    -1    \\
 0   &          0        &    0     &   0     &   0       & 0 &
 0     &   0       &   0   &   0     &   0       &   0     &
0 &     t   &          1-t        &    -1    \\
 -1   &          -1        &    -1     &   -1     &   -1       & -1 &
 -1     &   -1       &   -1   &   -1     &   -1       &   -1     &
0 &     0   &          t        &    2-2t
\end{matrix}
\right) $ }

\bigbreak

The next matrix is the result of subtracting the second row from
the first, the third from the second, and the fourth from the
third; the sixth from the fifth; the eighth from the seventh, the
ninth from the eighth, the tenth from the ninth, the eleventh from
the tenth, and the twelfth from the eleventh, in the preceding
matrix.

\bigbreak

\scalebox{.76}{ $ \setcounter{MaxMatrixCols}{16}  tS-S\sp{T}
\rightarrow \left(
\begin{matrix}
1-2t &  t-2   &          0        &    0     &   0     &   0 & 0 &
0   &          0       &   0     &   0       &   0   &
0     &   0      &   0 & 0\\
0   &    1-2t &  t-2   &          0        &    0     &   0     &
0 & 0 &     0   &         0
    &   0     &   0       &   0     &   0      &   0 & 0\\
0   &    0 &  1-2t   &          t-2        &    0     &   0     &
0 & 0 &     0   &            0
    &   0     &   0        & 0     &   0      &   0 & 0 \\
t   &    t &  t   &          1-t        &    0     &   0     &   0
& 0 &     0   &            0
    &   0     &   0       &  0     &   0      &   0 & t \\
0   &    0 &  0   &          0        &    1-2t     &   t-2     &
0 & 0 &     0   &            0
    &   0     &   0       &  0     &   0      &   0 & 0\\
    0     &   0       &   0   &   0     &   t       &   1-t     &
0 &  0   &          0        &    0     &   0     &   0
    & 0     &   0      &   0 & t\\
    0     &   0       &   0   &   0     &   0       &   0     &
1-2t &  t-2   &          0        &    0     &   0     &   0 &
0     &   0      &   0 & 0\\
    0     &   0       &   0   &   0     &   0       &   0     &
0 &  1-2t   &         t-2        &   0     & 0     &  0       &
0     &   0      &   0 & 0\\
    0     &   0       &   0   &   0     &   0       &   0     &
   0       &    0     &   1-2t  &  t-2   &       0     & 0
    & 0     &   0      &   0 & 0\\
    0     &   0       &   0   &   0     &   0       &   0   & 0
 &    0     &   0     &   1-2t       & t-2 &     0    &
0     &   0      &   0 & 0\\
    0     &   0       &   0   &   0     &   0       &   0     &  0
&    0     &   0     &   0       & 1-2t  &     t-2   &
0     &   0      &   0 & 0\\
    0     &   0       &   0      &   0     &
0 &  0   &        t      &    t     &   t     &   t       & t
&     1-t   &   0     &   0      &   0 & t\\
 0   &          0        &    0     &   0     &   0       & 0 &
 0     &   0       &   0   &   0     &   0       &   0     &
1 &     0   &          -t        &    -1    \\
 0   &          0        &    0     &   0     &   0       & 0 &
 0     &   0       &   0   &   0     &   0       &   0     &
t &     1   &          -t        &    -1    \\
 0   &          0        &    0     &   0     &   0       & 0 &
 0     &   0       &   0   &   0     &   0       &   0     &
0 &     t   &          1-t        &    -1    \\
 -1   &          -1        &    -1     &   -1     &   -1       & -1 &
 -1     &   -1       &   -1   &   -1     &   -1       &   -1     &
0 &     0   &          t        &    2-2t
\end{matrix}
\right) $ }

\bigbreak

The next matrix is the result of subtracting the second column
from the first, the third from the second, and the fourth from the
third; the sixth from the fifth; the eighth from the seventh, the
ninth from the eighth, the tenth from the ninth, the eleventh from
the tenth, and the twelfth from the eleventh, in the preceding
matrix.

\bigbreak

\scalebox{.76}{ $ \setcounter{MaxMatrixCols}{16}  tS-S\sp{T}
\rightarrow \left(
\begin{matrix}
3-3t &  t-2   &          0        &    0     &   0     &   0 & 0 &
0   &          0       &   0     &   0       &   0   &
0     &   0      &   0 & 0\\
2t-1   &    3-3t &  t-2   &          0        &    0     &   0 & 0
& 0 &     0   &         0
    &   0     &   0       &   0     &   0      &   0 & 0\\
0   &    2t-1 &  3-3t   &          t-2        &    0     &   0 & 0
& 0 &     0   &            0
    &   0     &   0        & 0     &   0      &   0 & 0 \\
0   &    0 &  2t-1   &          1-t        &    0     &   0     &
0 & 0 &     0   &            0
    &   0     &   0       &  0     &   0      &   0 & t \\
0   &    0 &  0   &          0        &    3-3t     &   t-2     &
0 & 0 &     0   &            0
    &   0     &   0       &  0     &   0      &   0 & 0\\
    0     &   0       &   0   &   0     &   2t-1       &   1-t     &
0 &  0   &          0        &    0     &   0     &   0
    & 0     &   0      &   0 & t\\
    0     &   0       &   0   &   0     &   0       &   0     &
3-3t &  t-2   &          0        &    0     &   0     &   0 &
0     &   0      &   0 & 0\\
    0     &   0       &   0   &   0     &   0       &   0     &
2t-1 &  3-3t   &         t-2        &   0     & 0     &  0       &
0     &   0      &   0 & 0\\
    0     &   0       &   0   &   0     &   0       &   0     &
   0       &    2t-1     &   3-3t  &  t-2   &       0     & 0
    & 0     &   0      &   0 & 0\\
    0     &   0       &   0   &   0     &   0       &   0   & 0
 &    0     &   2t-1     &   3-3t       & t-2 &     0    &
0     &   0      &   0 & 0\\
    0     &   0       &   0   &   0     &   0       &   0     &  0
&    0     &   0     &   2t-1       & 3-3t  &     t-2   &
0     &   0      &   0 & 0\\
    0     &   0       &   0      &   0     &
0 &  0   &        0      &    0     &   0     &   0       & 2t-1
&     1-t   &   0     &   0      &   0 & t\\
 0   &          0        &    0     &   0     &   0       & 0 &
 0     &   0       &   0   &   0     &   0       &   0     &
1 &     0   &          -t        &    -1    \\
 0   &          0        &    0     &   0     &   0       & 0 &
 0     &   0       &   0   &   0     &   0       &   0     &
t &     1   &          -t        &    -1    \\
 0   &          0        &    0     &   0     &   0       & 0 &
 0     &   0       &   0   &   0     &   0       &   0     &
0 &     t   &          1-t        &    -1    \\
 0   &          0        &    0     &   -1     &   0       & -1 &
 0     &   0       &   0   &   0     &   0       &   -1     &
0 &     0   &          t        &    2-2t
\end{matrix}
\right) $ }

\bigbreak

Using Laplace's expansion over the last row and conveniently
relocating the last column to compute the minors, we obtain

\bigbreak

\begin{multline*}
(-1)(-1)\det \left(
\begin{smallmatrix}
3-3t & t-2 & 0 & 0\\
2t-1 & 3-3t & t-2 & 0\\
0 & 2t-1 & 3-3t & 0\\
0  &  0  &  2t-1 & t
\end{smallmatrix}
\right)  \det \left(
\begin{smallmatrix}
3-3t & t-2 \\
2t-1 & 1-t
\end{smallmatrix}
\right)  \det  \left(
\begin{smallmatrix}
3-3t & t-2 & 0 & 0 & 0 & 0\\
2t-1 & 3-3t & t-2 & 0 & 0 & 0\\
0 & 2t-1 & 3-3t & t-2 & 0 & 0\\
0  &  0  &  2t-1 & 3-3t & t-2 &0\\
0  &  0  &   0  &  2t-1 & 3-3t & t-2\\
0  &  0  &   0  &  0    &  2t-1 & 1-t\\
\end{smallmatrix}
\right)   \det \left(
\begin{smallmatrix}
1 & 0 & -t\\
t & 1 & -t\\
0 & t & 1-t
\end{smallmatrix}
\right) \\
+(-1)(-1)\det \left(
\begin{smallmatrix}
3-3t & t-2 & 0 & 0\\
2t-1 & 3-3t & t-2 & 0\\
0 & 2t-1 & 3-3t & t-2\\
0  &  0  &  2t-1 & 1-t
\end{smallmatrix}
\right)  \det \left(
\begin{smallmatrix}
3-3t & 0 \\
2t-1 & t
\end{smallmatrix}
\right)  \det  \left(
\begin{smallmatrix}
3-3t & t-2 & 0 & 0 & 0 & 0\\
2t-1 & 3-3t & t-2 & 0 & 0 & 0\\
0 & 2t-1 & 3-3t & t-2 & 0 & 0\\
0  &  0  &  2t-1 & 3-3t & t-2 &0\\
0  &  0  &   0  &  2t-1 & 3-3t & t-2\\
0  &  0  &   0  &  0    &  2t-1 & 1-t\\
\end{smallmatrix}
\right)   \det \left(
\begin{smallmatrix}
1 & 0 & -t\\
t & 1 & -t\\
0 & t & 1-t
\end{smallmatrix}
\right) \\
+(-1)(-1)\det \left(
\begin{smallmatrix}
3-3t & t-2 & 0 & 0\\
2t-1 & 3-3t & t-2 & 0\\
0 & 2t-1 & 3-3t & t-2\\
0  &  0  &  2t-1 & 1-t
\end{smallmatrix}
\right)  \det \left(
\begin{smallmatrix}
3-3t & t-2 \\
2t-1 & 1-t
\end{smallmatrix}
\right)  \det  \left(
\begin{smallmatrix}
3-3t & t-2 & 0 & 0 & 0 & 0\\
2t-1 & 3-3t & t-2 & 0 & 0 & 0\\
0 & 2t-1 & 3-3t & t-2 & 0 & 0\\
0  &  0  &  2t-1 & 3-3t & t-2 &0\\
0  &  0  &   0  &  2t-1 & 3-3t & 0\\
0  &  0  &   0  &  0    &  2t-1 & t\\
\end{smallmatrix}
\right)   \det \left(
\begin{smallmatrix}
1 & 0 & -t\\
t & 1 & -t\\
0 & t & 1-t
\end{smallmatrix}
\right) \\
+t(-1)\det \left(
\begin{smallmatrix}
3-3t & t-2 & 0 & 0\\
2t-1 & 3-3t & t-2 & 0\\
0 & 2t-1 & 3-3t & t-2\\
0  &  0  &  2t-1 & 1-t
\end{smallmatrix}
\right)  \det \left(
\begin{smallmatrix}
3-3t & t-2 \\
2t-1 & 1-t
\end{smallmatrix}
\right)  \det  \left(
\begin{smallmatrix}
3-3t & t-2 & 0 & 0 & 0 & 0\\
2t-1 & 3-3t & t-2 & 0 & 0 & 0\\
0 & 2t-1 & 3-3t & t-2 & 0 & 0\\
0  &  0  &  2t-1 & 3-3t & t-2 &0\\
0  &  0  &   0  &  2t-1 & 3-3t & t-2\\
0  &  0  &   0  &  0    &  2t-1 & 1-t\\
\end{smallmatrix}
\right)   \det \left(
\begin{smallmatrix}
1 & 0 & -1\\
t & 1 & -1\\
0 & t & -1
\end{smallmatrix}
\right) \\
+(2-2t)\det \left(
\begin{smallmatrix}
3-3t & t-2 & 0 & 0\\
2t-1 & 3-3t & t-2 & 0\\
0 & 2t-1 & 3-3t & t-2\\
0  &  0  &  2t-1 & 1-t
\end{smallmatrix}
\right)  \det \left(
\begin{smallmatrix}
3-3t & t-2 \\
2t-1 & 1-t
\end{smallmatrix}
\right)  \det  \left(
\begin{smallmatrix}
3-3t & t-2 & 0 & 0 & 0 & 0\\
2t-1 & 3-3t & t-2 & 0 & 0 & 0\\
0 & 2t-1 & 3-3t & t-2 & 0 & 0\\
0  &  0  &  2t-1 & 3-3t & t-2 &0\\
0  &  0  &   0  &  2t-1 & 3-3t & t-2\\
0  &  0  &   0  &  0    &  2t-1 & 1-t\\
\end{smallmatrix}
\right)   \det \left(
\begin{smallmatrix}
1 & 0 & -t\\
t & 1 & -t\\
0 & t & 1-t
\end{smallmatrix}
\right) =\\
= 2 + 192t + 972t\sp{2} - 12289t\sp{3} + 49274t\sp{4} -
120582t\sp{5} + 213765t\sp{6} - 295426t\sp{7} + 328185t\sp{8}-\\
\shoveleft{ - 295426t\sp{9} + 213765t\sp{10} - 120582t\sp{11} +
49274t\sp{12} - 12289t\sp{13} + 972t\sp{14} + 192t\sp{15} +
2t\sp{16}}
\end{multline*}

\bigbreak

From this  matrix it is now clear how to obtain the matrix for any
given even number of tassels such that one of the tassels has an
even number of crossings.

\subsection{Calculations: pretzel knots on $N$ tassels, exactly one with an even number of
crossings}\label{subsect:poo...oe}

\noindent

Let us assume that there are $N$ tassels, for a given even
positive integer $N$, and that the $N$-th tassel has $2i\sb{N}$
crossings $(i\sb{N}>1)$, and each of the remaining tassels have
$2i\sb{k}+1$ tassels, for $1\leq k \leq N-1$ $(i\sb{k}\, 's >0)$.
The presentation matrix of the Alexander module of the Pretzel
knot $P(2i\sb{1}+1, 2i\sb{2}+1, \dots , 2i\sb{N-1}+1, 2i\sb{N})$
is equivalent to the sum of the following two $\bigl(
\sum\sb{j=1}\sp{N}2i\sb{j}\bigr) \times \bigl(
\sum\sb{j=1}\sp{N}2i\sb{j}\bigr) $ matrices, $M\sb{1}$ and
$M\sb{2}$.

$M\sb{1}$ is a matrix whose last column and last row have non-zero
entries, all other entries being zero. The only non-zero entries
of the last row are $-1$'s at the $\sum\sb{j=1}\sp{k}2i\sb{j}$
positions, for $1 \leq k \leq N-1$. The only non-zero entries of
the last column are $t$'s at the $\sum\sb{j=1}\sp{k}2i\sb{j}$
positions, for $1 \leq k \leq N-1$.

$M\sb{2}$ is a matrix in block-diagonal form. Each of the blocks
is a $2i\sb{k} \times 2i\sb{k}$ matrix, for $1\leq k \leq N$. For
$1\leq k \leq N-1$, the $k$-th block is of the following type:

\bigbreak

\begin{equation*} B\sb{k} = \left(
\begin{matrix}
3-3t &  t-2   &          0        &    \cdots    &   \cdots     &
0
& 0 & 0   &          0       &   0     \\
2t-1   &    3-3t &  t-2   &          0        &    \cdots      &
\cdots & 0 & 0 &     0   &         0\\
0   &    2t-1 &  3-3t   &          t-2        &    0      &
\cdots & \cdots & 0 &     0   &            0\\
\vdots    &    \vdots &  \vdots & \ddots   &     \vdots      &
\vdots & \vdots & \vdots &     \vdots   &            \vdots\\
\vdots    &    \vdots &  \vdots &     \vdots  & \ddots       &
\vdots & \vdots & \vdots &     \vdots   &            \vdots\\
0   &   \cdots      &   \cdots     &
   0       &    2t-1     &   3-3t  &  t-2   &       0   &  0     &     0  \\
0   &   0     &   \cdots       &
    \cdots &    0     &   2t-1     &   3-3t       & t-2 &     0  &   0\\
0   &   0     &   0       &  \cdots       &
    \cdots &    0     &   2t-1     &   3-3t       & t-2 &     0  \\
   0   &   0     &   0       &   0     &
    \cdots
&    \cdots     &   0     &   2t-1       & 3-3t  &     t-2\\
 0      &   0     & 0 &  0   &        0
&    \cdots     &   \cdots     &   0       & 2t-1 & 1-t
\end{matrix}
\right)
\end{equation*}

\bigbreak

The $N$-th block is of the following type:

\bigbreak

\begin{equation*} B\sb{N} = \left(
\begin{matrix}
1 &  0   &          0        &    \cdots    &   \cdots     & 0
& 0 & 0   &          -t       &   -1     \\
t   &    1 &  0   &          0        &    \cdots      &
\cdots & 0 & 0 &     -t   &         -1\\
0   &    t &  1   &          0        &    0      &
\cdots & \cdots & 0 &     -t   &            -1\\
\vdots    &    \vdots &  \vdots & \ddots   &     \vdots      &
\vdots & \vdots & \vdots &     \vdots   &            \vdots\\
\vdots    &    \vdots &  \vdots &     \vdots  & \ddots       &
\vdots & \vdots & \vdots &     \vdots   &            \vdots\\
0   &   \cdots      &   \cdots     &
   0       &    t     &   1  &  0   &       0   &  -t     &     -1  \\
0   &   0     &   \cdots       &
    \cdots &    0     &   t     &   1       & 0 &     -t  &   -1\\
0   &   0     &   0       &  \cdots       &
    \cdots &    0     &   t     &   1       & -t &     -1  \\
   0   &   0     &   0       &   0     &
    \cdots
&    \cdots     &   0     &   t       & 1-t  &     -1\\
 0      &   0     & 0 &  0   &        0
&    \cdots     &   \cdots     &   0       & t & I-It
\end{matrix}
\right)
\end{equation*}
where $I$ is such that $N=2I$.

\bigbreak

The Alexander polynomial of $P(2i\sb{1}+1, 2i\sb{2}+1, \dots ,
2i\sb{N-1}+1, 2i\sb{N})$ is the determinant of $M\sb{1}+M\sb{2}$.
We compute it by doing Laplace's expansion on the last row of
$M\sb{1}+M\sb{2}$. In order to do that the following matrices will
be helpful. For $1\leq k \leq N-1$, $B'\sb{k}$ is a $2i\sb{k}
\times 2i\sb{k}$ matrix; it is obtained by replacing the last
column of $B\sb{k}$ by a column with a $t$ in the last entry and
otherwise zero. $B''\sb{N}$ is a $(2i\sb{N}-1) \times
(2i\sb{N}-1)$ matrix; it is obtained by removing the last column
and the last row from $B\sb{N}$. $B'\sb{N}$ is also a
$(2i\sb{N}-1) \times (2i\sb{N}-1)$ matrix; it is is obtained by
removing the column before the last one and the last row from
$B\sb{N}$.

\bigbreak

\begin{equation*} B'\sb{k} = \left(
\begin{matrix}
3-3t &  t-2   &          0        &    \cdots    &   \cdots     &
0
& 0 & 0   &          0       &   0     \\
2t-1   &    3-3t &  t-2   &          0        &    \cdots      &
\cdots & 0 & 0 &     0   &         0\\
0   &    2t-1 &  3-3t   &          t-2        &    0      &
\cdots & \cdots & 0 &     0   &            0\\
\vdots    &    \vdots &  \vdots & \ddots   &     \vdots      &
\vdots & \vdots & \vdots &     \vdots   &            \vdots\\
\vdots    &    \vdots &  \vdots &     \vdots  & \ddots       &
\vdots & \vdots & \vdots &     \vdots   &            \vdots\\
0   &   \cdots      &   \cdots     &
   0       &    2t-1     &   3-3t  &  t-2   &       0   &  0     &     0  \\
0   &   0     &   \cdots       &
    \cdots &    0     &   2t-1     &   3-3t       & t-2 &     0  &   0\\
0   &   0     &   0       &  \cdots       &
    \cdots &    0     &   2t-1     &   3-3t       & t-2 &     0  \\
   0   &   0     &   0       &   0     &
    \cdots
&    \cdots     &   0     &   2t-1       & 3-3t  &     0\\
 0      &   0     & 0 &  0   &        0
&    \cdots     &   \cdots     &   0       & 2t-1 & t
\end{matrix}
\right)
\end{equation*}

\bigbreak

\bigbreak

\begin{equation*} B'\sb{N} = \left(
\begin{matrix}
1 &  0   &          0        &    \cdots    &   \cdots     & 0
& 0 & 0   &            -1     \\
t   &    1 &  0   &          0        &    \cdots      &
\cdots & 0 & 0 &           -1\\
0   &    t &  1   &          0        &    0      &
\cdots & \cdots & 0 &             -1\\
\vdots    &    \vdots &  \vdots & \ddots   &     \vdots      &
\vdots & \vdots & \vdots &     \vdots   \\
\vdots    &    \vdots &  \vdots &     \vdots  & \ddots       &
\vdots & \vdots & \vdots &     \vdots   \\
0   &   \cdots      &   \cdots     &
   0       &    t     &   1  &  0   &       0     &     -1  \\
0   &   0     &   \cdots       &
    \cdots &    0     &   t     &   1       & 0  &   -1\\
0   &   0     &   0       &  \cdots       &
    \cdots &    0     &   t     &   1     &     -1  \\
   0   &   0     &   0       &   0     &
    \cdots
&    \cdots     &   0     &   t       &     -1
\end{matrix}
\right)
\end{equation*}

\bigbreak

\bigbreak

\begin{equation*} B''\sb{N} = \left(
\begin{matrix}
1 &  0   &          0        &    \cdots    &   \cdots     & 0
& 0 & 0   &          -t           \\
t   &    1 &  0   &          0        &    \cdots      &
\cdots & 0 & 0 &     -t   \\
0   &    t &  1   &          0        &    0      &
\cdots & \cdots & 0 &     -t   \\
\vdots    &    \vdots &  \vdots & \ddots   &     \vdots      &
\vdots & \vdots & \vdots &     \vdots   \\
\vdots    &    \vdots &  \vdots &     \vdots  & \ddots       &
\vdots & \vdots & \vdots &     \vdots   \\
0   &   \cdots      &   \cdots     &
   0       &    t     &   1  &  0   &       0   &  -t      \\
0   &   0     &   \cdots       &
    \cdots &    0     &   t     &   1       & 0 &     -t  \\
0   &   0     &   0       &  \cdots       &
    \cdots &    0     &   t     &   1       & -t  \\
   0   &   0     &   0       &   0     &
    \cdots
&    \cdots     &   0     &   t       & 1-t  \\
\end{matrix}
\right)
\end{equation*}

\bigbreak

The Alexander polynomial is then:

\[
\det B''\sb{N}\sum\sb{j=1}\sp{N-1}\det
B'\sb{j}\prod\sb{\substack{l=1\\l\neq j}}\sp{N-1}\det B\sb{l} -
t\det B'\sb{N}\prod\sb{l=1}\sp{N-1}\det B\sb{l}  + (I-It)\det
B''\sb{N}\prod\sb{l=1}\sp{N-1}\det B\sb{l}
\]

We remark that it is easy to see, by Laplace expansion on the last
column, that

\[
\det
B'\sb{N}=\sum\sb{j=1}\sp{2i\sb{N}-1}(-1)(-1)\sp{j+1}t\sp{2i\sb{N}-1-j}
\]

and

\[
\det
B''\sb{N}=\sum\sb{j=1}\sp{2i\sb{N}-2}(-t)(-1)\sp{j+1}t\sp{2i\sb{N}-1-j}
+(1-t)
\]

\end{document}